\documentclass[12pt]{article}
\newcommand{\free}[1] {\ensuremath{\hspace{0.2cm} (\kern -.7em{<}#1{>}\kern-.7em)\hspace{0.2cm}}}
\usepackage{latexsym}
\usepackage{amsmath}
\renewcommand{\wr}{\operatorname{wr}}

\begin{document}

\begin{titlepage}
\title {Restricted Lie algebras of polycyclic groups, II.}
\author {A. I. Lichtman,\\
\vspace{0.5cm}
 lichtman@uwp.edu}

\maketitle
\noindent
\begin{abstract} 

  This is the second paper   in the series of three.  We study restricted Lie algebras  of polycyclic groups and  obtain conditions for existence  of $p$-series  with associated  restricted Lie algebra abelian  or free abelian with rank equal to the Hirsch number of the group. We develop methods for constructing such series,  these  methods are  based on construction of filtrations and valuation functions in the group rings.

This paper continues author's work [9] and [10].

\end{abstract}
\end{titlepage}
\setcounter{page}{2}
\setcounter{section}{1}
\setcounter{equation}{0}
\section*{\center \S 1. Statement of the results. Notation.}

{\bf 1.1.}   Let $H$ be a group.  A series of normal subgroups 

  \begin{equation}H=H_1\supseteq H_2\supseteq\cdots  \end{equation}

 is a $p$-series in $H$ if $[H_i,H_j]\subseteq H_{i+j}$ and $H_i^p\subseteq H_{ip}$. We will denote by $L_p(H,H_i)$ the restricted Lie algebra associated to $p$-series $(1.1)$ and by $U_p(L_p(H,H_i))$ its universal $p$-envelope. The algebra 
\begin{equation}L_p(H,H_i)=\sum_{i=1}^{\infty}H_i/H_{i+1}  \end{equation}
is obtained by the classical construction of Lazard [7];  we give in section $2$ a brief description of this construction and of the main properties of it. The  properties of the algebra $L_p(H,H_i)$  and the properties  of the algebra $U_p(L_p(H,H
_i))$ depend not only on the group $H$ but also on choice of the $p$-series in this group.  

We studied in Lichtman $[10]$ the restricted Lie algebras of  polycyclic groups  and obtained necessary and sufficient conditions for existence  in  a poly-{infinite cyclic} group with Hirsch number $r$ a series $(1.1)$ with associated Lie algebra $L_p(H,H_i)$ free abelian of rank $r$.  We obtained also in $[10]$ some necessary conditions for existence of such series in an arbitrary polycyclic group. 

The  main results of this paper are Theorems I-XII. We begin by formulating Theorem V which will be proven in section $5$. \bigskip

 {\bf Theorem  V.} {\it Let $H$ be an infinite    polycyclic-by-finite group with Hirsh number $r$.  Assume that there exists    a $p$-series $(1.1)$ with unit intersection such that the corresponding restricted Lie algebra $L_p(H,H_i)$ is finitely generated.    Then there exists a 
torsion free  normal  subgroup $F$   with index a power of $p$ such that the ideal $L_p(F,F_i)$ associated to the $p$-series $F_i=F\bigcap H_i\,\, (i=1,2,\cdots)$ is a restricted free abelian  subalgebra of the center of  $L_p(H,H_i)$
   with  index  a power of $p$ and rank $1\leq r_1\leq r$.

Hence  the  center $Z$ of $L_p(H,H_i)$  has a finite index which is a power of $p$ and $L_p(H,H_i)$ is a nilpotent Lie algebra.} \bigskip

It is known that   finitely generated restricted  abelian Lie algebra $F$ has a representation

\begin{equation} F=T+T_0\end{equation}
where $T$ is either a free abelian subalgebra or $T=0$  and $T_0$ is finite; this follows also from classical theorems  about modules over a polynomial ring.  {\it We will call throughout the paper  the rank  of the free abelian 
subalgebra $T$   the rank of $F$}. \bigskip 

We recall also that an element $x$ of a restricted Lie algebra is nilpotent if there exists $p^n$ such that $x^{[p]^n}=0$, and that a finitely generated restricted Lie algebra without nilpotent elements is free; the last fact follows also from the representation $(1.3)$

Before formulating   Theorems VI-XII  we recall first the elementary  fact  that  if $(1.1)$ is a $p$-series in an arbitrary group $H$, and $H_0=\bigcap_{i=1}^{\infty}H_i$ then the Lie algebra $L_p(H,H_i)$ is isomorphic  to the restricted Lie algebra $L_p(\bar{H},\bar{H}_i)$ of the group $ \bar{H}=H/H_0$ associated to the series $\bar{H}_i=H_i/H_0\,\, (i=1,2,\cdots)$  (see section $2$). The group $\bar{H}$ is a  residually  finite $p$-group, and the series $\bar{H}_i\,\, (i=1,2,\cdots)$ has unit intersection; this shows that in the study of  the algebra $L_p(H,H_i)$ we can assume that series $(1.1)$ has a unit intersection and, we can consider only the case when  $H$ is a residually  finite $p$-group.  {\it We will use throughout the paper the notation $H\in res\, \mathcal{N}_p$ for  the fact that $H$ is a residually finite $p$-group.} 

Further, we will consider in this case the topology defined by    series $(1.1)$.  If 

\begin{equation}M_n(Z_pH)\,\, (n=1,2,\cdots)\end{equation}

 is the series of dimension subgroups in characteristic $p$ (see section $2$) then the  topology defined by this series will be called throughout the paper by  $p$-topology. We will prove that if $H$ is a polycyclic group with Hirsch number $r$,  $(1.1)$ is a $p$-series with unit intersection and the algebra $L_p(H,H_i)$ is abelian of  rank $r$ then the topology defined by this series is equivalent to the $p$-topology. This means that for every dimension subgroup  $M_n(Z_p H)$ a number $i(n)$ can be found such that $H_{i(n)}\subseteq M_n(Z_pH)$;  it is worth remarking that the  inclusion $H_i\supseteq M_i(Z_p H)\,\, (i=1,2,\cdots)$ holds in  an arbitrary group.

  The  necessary  and sufficient conditions for  existence in a polycyclic group  of a $p$-series   with associated restricted Lie algebra abelian of rank $r$ are obtained in the following   Theorems  VII and VI. \bigskip 

{\bf Theorem VII.} {\it  Let $H$ be a polycyclic group with Hirsch number $r$. There exists in $H$ a  $p$-series $(1.1)$ with unit intersection and associated restricted Lie algebra $L_p(H,H_i)$ abelian of rank $r$ if and only if $H$ contains normal subgroups $Q\supseteq N$ such that the quotient group $H/Q$ is a $p$-group, $Q/N$ is a  free abelian group,  $N$ is torsion free nilpotent, and
for every element $h\in H$ the following condition holds}

 \begin{equation} [h^{p^r},N]\subseteq N^{\prime}N^p\end{equation}   

{\it  or, equivalently, if  $R=gp(h,N)$  then the quotient group $\bar{R}=R/N^{\prime}N^p\in res \, \mathcal{N}_p$}

The necessity of the conditions of Theorem VII follow from statement vi) of Theorem VI which we will now formulate; we will show    in subsection $1.3.$ that the sufficiency of these  conditions follows from Theorem XII 
and Corollary $6.4.$ \bigskip

{\bf Theorem VI.} {\it Let $H$ be a polycyclic group with Hirsch number $r$. Assume that there exists a $p$-series $H_i\,\, (i=1,2,\cdots)$ with unit intersection such that $L_p(H,H_i)$ is  abelian of rank $r$. }

i)   {\it Let $U$ be an arbitrary subgroup of $H$ with Hirsch number $k$, $U_i=U\bigcap H_i\,\, (i=1,2,\cdots)$. Then $L_p(U,U_i)$ is an abelian algebra of rank $k$.} 

ii)  {\it  Let $U$ be a normal subgroup of $H$ with Hirsch number $k$,   $\bar{H}_i$ be the image of the subgroup $H_i$ in $H/U$. Then the  subgroup $\bigcap_{i=1}^{\infty} \bar{H}_i$ is finite and the algebra $L_p(\bar{H},\bar{H}_i)$ is abelian of rank $r-k$. In particular, if
 $\bar{H}=H/U$ contains no finite normal subgroups then it is a residually  
 \{finite $p$-group\}. }

iii) {\it  Let $U$ be a normal subgroup of $H$. If  $\bar{H}=H/U$ is a residually \{finite $p$-group\} then $\bigcap_{i=1}^{\infty}\bar{H}_i=1$}.

iv)   {\it Let $W$ be the unique maximal normal nilpotent-by-finite subgroup of $H$. Then $W$ is an extension of a torsion free nilpotent group by a finite $p$-group.  The quotient group $H/W$ is an extension of a free abelian group by a finite $p$-group.}

v) {\it The topology defined in $H$ by the $p$-series $H_i\,\, (i=1,2,\cdots)$  is equivalent to the $p$-topology. 

The topologies  defined in an arbitrary subgroup $U$   by the series  $U_i=U\bigcap H_i\,\, (i=1,2,\cdots)$ and    $U\bigcap M_n(H)\,\, (n=1,2,\cdots) $ are equivalent to the $p$-topology in $U$.}

vi)  {\it There exists an index $i_0$ such that if $i\geq i_0$ then the   subgroup  $Q=H_i$  contains  a torsion free nilpotent subgroup $N$  which is invariant in $H$,    $Q/N$ is free abelian, the algebra $L_p(Q,Q_i)$ is free abelian of rank $r$  and }  

\begin{equation}H/N^{\prime}N^p\in res\,\mathcal{N}_p\end{equation}

{\it Clearly, $H/Q$ is a finite $p$-group.}

vii)  {\it Let $F\supseteq S$ be two normal subgroups in $H$ such that $H/F$ and $F/S$ are residually \{finite $p$-groups\}. Then $H/S$ is a residually \{finite $p$-group\}}. 

viii)  {\it  Let 

 \begin{equation} H=H^{\ast}_1\supseteq H^{\ast}_2\supseteq\cdots\end{equation}   

be a series in $H$  with unit intersection and finitely generated associated graded Lie algebra $L_p(H,H^{\ast}_i)$. If the  topology defined by series $(1.7)$ is equivalent to the $p$-topology then the center of  $L_p(H,H^{\ast}_i)$ has rank $r$.  Moreover, there exists a number $k$ such that if  
$U=H_i \,\, (i\geq k)$ then the subalgebra $L_p(U,U_i^{\ast})\cong \sum_{i\geq k} H_i/H_{i+1}$ associated to the $p$-series $U_i^{\ast}=U\bigcap H_i^{\ast}\,\,  (i=1,2,\cdots)$ is a  central free abelian subalgebra of rank $r$. } \bigskip

 Condition $(1.5)$ in Theorem VII holds  in an arbitrary group $H$ which contains a normal subgroup $N$ such that  the quotient group $H/N^{\prime}N^p$ is a residually \{finite $p$-group\}. On the other hand, if $H$ is a polycyclic group with Hirsch number $r$ which contains a $p$-series $(1.1)$ with associated graded algebra $L_p(H,H_i)$ abelian of rank $r$ and $Q$ and $N$ are the normal subgroups obtained  in Theorems VII and VI we obtain from statement vi) of Theorem VI that $H/N^{\prime}N^p$ is a residually \{finite $p$-group\}.

Theorem VI will be proven in section $6$.  We  will show in subsection $1.3.$ that  Theorem VII follows  from Theorem VI and from Theorem XII which will be formulated in subsection $1.2.$

We prove in section $6$ Theorem $6.1.$ which gives an additional information about the normal subgroups $Q$ and $N$ which are obtained in Theorems VI and VII.

Theorem V implies that if the algebra $L_p(H,H_i)$ is finitely generated then the rank of the center of it  is less than or equal to the Hirsch number of $H$.
 We prove in section $7$   Theorems XI   which gives the necessary and sufficient conditions  for  the center $Z$ to have rank  equal to the Hirsch number of $H$.  \bigskip

{\bf 1.2.} Theorems VI and VII  are related to an existence of a $p$-series $H_i\,\, (i=1,2,\cdots)$ with associated graded Lie algebra abelian of rank $r$, where $r$ is the Hirsch number of $H$. We consider the question when the algebra $L_p(H,H_i)$ is {\it free abelian}  of rank $r$ in Theorem XII. We have already pointed out that we should  consider only the series with unit intersection; further,   it is easy to see (Lemma $2.1.$) that  if the algebra $L_p(H,H_i)$ is free abelian  then the group $H$ must be torsion free. 

    Theorem XII will be proven in section $8$. This theorem  provides a  sufficient condition for existence in a torsion free polycyclic group $H$ with Hirsch number $r$ of a $p$ series $(1.1)$ with restricted Lie algebra $L_p(H,H_i)$ free abelian of rank $r$.  Theorem XII  generalizes  author's results   [10]; these results were obtained for the poly-\{infinite cyclic\} groups. \bigskip

  {\bf Theorem XII.} {\it Let $H$ be a torsion free polycyclic group with Hirsch number $r$, $U$ be a normal subgroup with Hirsch number $k$ and torsion free quotient group $\bar{H}=H/U$.   }

{\it Assume that the following $3$ conditions hold.} 

1) {\it There exists in $U$ a $p$-series }
 \begin{equation} U=U_1\supseteq U_2\supseteq\cdots \end{equation} 

{ \it with associated restricted Lie algebra $L_p(U,U_i)$ free abelian of rank $k$.}

2) {\it  There exists in the group $\bar{H}=H/U$ a $p$-series} 

\begin{equation} \bar{H}=\bar{H}_1\supseteq \bar{H}_2\supseteq\cdots\end{equation} 
{\it with associated restricted  Lie algebra $L_p(\bar{H},\bar{H}_i)$ free abelian of rank $r-k$.} 

3) {\it For for every subgroup $R=gp(h,U)$ generated by $U$ and an element $h\in H$ the quotient group $\bar{R}=R/U^{\prime}U^p\in res\,\,\mathcal{N}_p$   or, equivalently, $[h^{p^k}, U]\subseteq U^{\prime}U^p$}. 

{\it Then there exists a $p$-series $(1.1)$ } { \it  with unit intersection and  associated restricted Lie algebra $L_p(H,H_i)$ free abelian of rank $r$ such that  } 

\begin{equation} \bar{H}_i=(H_iU)/U\,\, (i=1,2,\cdots) \end{equation}  

The  last statement of Theorem XII together with the classical results of Lazard   implies immediately the following   corollary. \bigskip 

{\bf Corollary 8.3.} {\it  The natural homomorphism $\phi\colon H\longrightarrow \bar{H}$ defines a homomorphism  of graded algebras}

\begin{equation} \tilde{\phi}\colon L_p(H,H_i)\longrightarrow L_p(\bar{H},\bar{H}_i)\end{equation}

Below are a few remarks on the conditions of Theorem XII.

First,  the series $H_i,U_i,\bar{H}_i\,\, (i=1,2,\cdots)$ in Theorem XII have unit intersection. This follows from Proposition $3.5.$ in Lichtman $[10 ]$; it follows also  from  statement ii) of  Theorem   VI.

Second, everyone of   conditions  $1)$ and $3)$ is necessary. The necessity of condition $1)$ follows from the fact that the subalgebra of a free abelian algebra is free abelian, and the rank of $L_p(U,U_i)$ must be $k$ via statement i) of Theorem VI. The necessity of condition $3)$ follows from statement vii) of Theorem VI.

Third, it  is easy to show that condition $3)$  does not follow from conditions $1$ and $2$ even if $H\in res\, \mathcal{N}_p$ (see, for instance, [10], section $10.3.$)  but it holds if the group $H$ is an extension of a torsion free nilpotent group by a finite $p$-group.  We have for this class of groups the following immediate corollary of Theorem XII. \bigskip

{\bf Corollary 1.1.} {\it Let    $H$ be a torsion free group which  contains    a nilpotent normal subgroup   whose Hirsch number is $r$ and index  is a power of $p$. Let $U$ be a normal subgroup of $H$ which satisfies  conditions $1)$ and $2)$ of Theorem XII. Then there exists in $H$ a $p$-series which satisfies all the conclusions of Theorem XII.} \bigskip

We have one more corollary of  Theorem XII.\bigskip

{\bf Corollary 1.2.} {\it Let $H$ be a torsion free polycyclic group with Hirsch number $r$ which contains a nilpotent  normal subgroup $U$ with Hirsch number $k$ such that  the quotient group $H/U^pU^{\prime}$ is a residually finite $p$-group. Assume that the quotient group $\bar{H}=H/U$ is torsion free and contains a $p$-series $(1.9)$ such that the algebra $L_p(\bar{H},\bar{H}_i)$ is free abelian of rank $r-k$. Then the  group $H$ contains a $p$-series $(1.1)$ with unit intersection such that the algebra $L_p(H,H_i)$ is free abelian of rank $r$.}\bigskip 

{\bf Proof.}  The group $U$ contains a series of normal subgroups of  length $k$ with infinite cyclic factors.  Theorem XII together with a straightforward induction argument yields  that $U$ contains a $p$-series $(1.8)$ with associated restricted Lie algebra free abelian of rank $k$.
We apply now Theorem XII to the group $H$ and its normal subgroup $U$ and the assertion follows. \bigskip

 We make in the proof of Theorem XII  an essential use of Theorems VI,VII, IX and X. \bigskip

{\bf Theorem X.} {\it Let $H$ be a torsion free   polycyclic group with Hirsch number $r$ which contains a $p$-series  $(1.1)$  
 with unit intersection. Assume  that the Lie algebra $L_p(H,H_i)$ is free abelian of rank $r$. Let $\Phi$ be a group of automorphisms of $H$ such that the order of every automorphism $\phi\in  \Phi$ on  the quotient group $H/H^{\prime}H^p$ is a power of $p$. }

 { \it Then there exists a $p$-series }

$$ H=H_1^{\ast}\supseteq H_2^{\ast}\supseteq\cdots$$

 {\it with unit intersection such that  all the subgroups $H_i^{\ast}\,\, (i=1,2,\cdots)$ are $\Phi$-invariant,   $\Phi$ centralizes all the factors $H_i/H_{i+1}\,\, (i=1,2,\cdots)$ and  the algebra $L_p(H,H_i^{\ast})$  is free abelian of rank $r$. }  \bigskip

 {\bf Theorem IX.} {\it Let $H$ be a torsion free polycyclic group with Hirsch number $r$. Assume that there exists a $p$-series $(1.1)$ with associated restricted Lie algebra $L_p(H,H_i)$ free abelian (abelian) of rank $r$. Then there exists a $p$-series of characteristic subgroups with associated restricted Lie algebra free abelian (abelian) of rank $r$}.\bigskip

Theorem IX is obtained as a corollary of the following more general result.\bigskip

{\bf Theorem VIII.} {\it Let $H$ be a  finitely generated torsion free group which has  a  $p$-series $H_i\,\, (i=1,2,\cdots)$ with unit intersection and with the associated restricted Lie algebra $L_p(H,H_i)$ free abelian (abelian) of finite rank. Assume that the topology defined by this $p$-series is equivalent to the $p$-topology. Then there exists a $p$-series

 \begin{equation}H=U_1\supseteq U_2\supseteq\cdots\end{equation} 

whose terms $U_i\,\, (i=1,2,\cdots)$ are characteristic subgroups and the Lie algebra $L_p(H,U_i)$ is free abelian (abelian) of finite rank.} \bigskip

{\bf 1.3.}   {\bf Derivation  of Theorem VII from  Theorem VI and  Proposition 7.5.} We have already observed that the necessity of the conditions of Theorem VII is in fact statement vi) of Theorem VI. 

To prove the sufficiency we consider the normal subgroup $Q$ which was obtained in statement vi) of Theorem VI.  Corollary $1.2.$ implies that  $Q$  contains a $p$-series with unit intersection 

\begin{equation} Q=Q_1\supseteq Q_2\supseteq\cdots\end{equation}

and with associated restricted Lie algebra $L_p(Q,Q_i)$ free abelian of rank $r$.  We see now that  sufficiency of the conditions of Theorem VII will  follow  from Theorem XII and the following fact which will be proven in section $7$.\bigskip

{\bf Proposition 7.5.}  {\it  Let $H$ be a polycyclic group with Hirsch number $r$, $U$ be a normal subgroup of finite index $(H\colon U)=p^n$. Assume that there exists $p$-series }

\begin{equation} U=U_1\supseteq U_2\supseteq\cdots\end{equation}

{\it  with unit intersection such that the algebra $L_p(U,U_i)$ is abelian of rank $r$. Then there exists a $p$-series $(1.1)$ such that the algebra $L_p(H,H_i)$ is abelian of rank $r$. } \bigskip

{\bf 1.4. Polycentral systems in rings.}   Our methods are based on  the results   of section $3$ on polycentral system in rings. Before formulating these results we define first the following  concept which will be used throughout the whole paper.

 Let $R$ be a ring, $t_1,t_2,\cdots, t_n$ be a system of elements in $R$. Assume that the element  $t_1$ is central, the  ideal $(t_1)$ generated by $t_1$ is residually nilpotent   and $t_1$ is  regular in $R$ or, equivalently, the graded ring associated to the filtration $(t_1)^i\,\, (i=1,2,\cdots)$ is isomorphic to the polynomial ring $(R/A)[t_1]$ (see Corollary $3.1.$); further, for every $1\leq i\leq n-1$ the  element $t_{i+1}$ is central and regular modulo the ideal $A_i=<t_1,t_2,\cdots,t_i>$ generated by $t_1,t_2,\cdots, t_i$ and the ideal $(t_{i+1})$ is residually nilpotent in $R/A_i$. If these conditions hold we  will say that this system is polycentral independent; if all the elements $t_i\,\, (i=1,2,\cdots,n)$
are central in $R$ we will say that the system is central independent.

We will also consider this situation in a more detailed way taking into account the number of central elements in the system $T$ and its subsystems. Let $T_1$ be a central independent system in $R$, $A_1$ be the ideal generated by $T_1$; for every $1\leq i\leq n-1$ let $T_{i+1}$ be a system of elements which is central and independent modulo the ideal $A_i$ generated by the system $T_1\bigcup T_2\bigcup\cdots\bigcup  T_i$. Clearly the system $T=<T_1, T_2,\cdots  T_m>$ is polycentral independent in $R$. We will use this notation in order to make clear that   $T$ is a polycentral independent  system which is composed  from the independent systems $T_1,T_2,\cdots, T_m$.  It is worth remarking that we do not assume that the subsystems  $T_i\,\, (i=2,3,\cdots,m)$ are central in $R$; we assume that  every subsystem $T_i$  must be  central (and independent) in the quotient rings $R/A_{i-1}\,\, (i=2,3,\cdots,m)$.   

We  order the elements of every $T_i\,\, (i=1,2,\cdots,m)$ in an arbitrary way and then extend these orders to an order in $T$ by assuming that the elements of $T_i$ preceed the elements of $T_{i+1}$. The standard monomials on $T$  are defined in the usual way; it is convenient to assume that $1$ is a standard monomial of degree zero. 
  We will construct  independent polycentral systems in some classes of group rings of torsion free polycyclic groups, and in other classes of rings but at this point we will only mention the following two  cases where these systems are obtained easily.\bigskip 

$1)$ Let  $L$ is a  nilpotent Lie algebra with a central series 

\begin{equation} L=L_1\supseteq L_2\supseteq\cdots\supseteq L_{k-1}\supseteq  L_k=0\end{equation}

We pick in  every $L_{k-i}\,\, (i=1,2,\cdots,k-1)$ a system of elements $T_i$ which forms a basis of the quotient space $L_{k-i}/L_{k-i+1}$;  
 in particular, the system $T_1$ is a basis of $L_{k-1}$  and is central in $L$. The system of elements $T_1\bigcup T_2\bigcup\cdots \bigcup T_{k-1}$ is an independent  polycentral system in the universal eneveloping algebra $U(L)$.\bigskip 

$2)$ Let $H$ be a finitely generated torsion free nilpotent group. Let 

\begin{equation}  H=H_1\supseteq H_2\supseteq\cdots\supseteq  H_{k-1}\supseteq H_k=0\end{equation}
be a central series in $H$ with torsion free factors $H_i/H_{i+1}\,\, (i=1,2,\cdots,k-1)$; we recall that   the factors of the upper central series of $H$ have this property. We pick  in $H_{k-i}\,\, (i=1,2,\cdots,k-1)$ a system of elements  $E_i$ which forms a basis of the free abelian group $H_{k-i}/H_{k-i+1}$. Let $E_i-1=\{e-1|e\in E_i\}\,\, (i=1,2,\cdots,k-1)$. 
The system 
$E_1-1,E_2-1,\cdots, E_{k-1}-1$ is a polycentral independent system in the group ring $KH$ of $H$ over an arbitrary field.  Further, we consider the group ring of $H$  over  the ring of integers or the ring of 
 $p$-adic integers $\Omega$; in this case the system 

\begin{equation} p,E_1-1, E_2-1,\cdots, E_{k-1}-1\end{equation} 

is an independent polycentral system in $\Omega H$. The same is true in a group ring $CH$  over an arbitrary ring $C$ of characteristic zero such that the powers of the ideal $(p)$ define a $p$-adic valuation in $C$.  

More generally, we consider a torsion free nilpotent group $H$ without elements of infinite $p$-height and construct  filtrations and  valuations in the group rings over a field of $K$ of characteristic $p$ or over the ring of the integers $C$.   These results are obtained in Theorem IV and Corollary $4.1.$ We prove also Theorem IV$^{\prime}$ which is an analog of these results when the characteristic of $K$ is zero. Theorem IV$^{\prime}$ provides also a new proof of Hall-Hartley Theorem about the residual nilpotence of the augmentation ideal $\omega(KH)$.  

Polycentral ideals were considered by J.Roseblade and by P. Smith  in group rings of polycyclic groups and in Noetherian rings (see Passman [13], section $11$), and by Passman in [14] in connection with the AR-property and the localization theory.  We consider here a different situation when the rings are in general non-noetherian and our results are related to the valuations defined by the polycentral systems.

Our applications of the polycentral systems are based on the following Theorem I and II  which will be proven in section $3$.\bigskip

{\bf Theorem I.} {\it Let $R$ be a ring, $T=<t_1,t_2,\cdots,t_n>$  be an independent  polycentral system in $R$.}

 {\it The ideal $A$ generated by the system $T$ is residually nilpotent. If $\tilde{R}$ is the completion of $R$ in the topology defined by this ideal and $\mathcal{X}$ is a system of coset representatives for the elements of the quotient ring $R/A$ then every element $x\in \tilde{R}$ has a unique representation }

\begin{equation} x=\sum_{n=0}^{\infty} \lambda_n\pi_n\end{equation}

{\it where $\lambda_n\in \mathcal{X}\,\, (n=0,1,\cdots) $,  
  $\pi_n$ are standard monomials on $T$, and $lim_{n\to\infty}v(\pi_n)=\infty$.}\bigskip
\bigskip

 Let $C$ be the set of integers. Here and throughout the paper we mean that a function  $\rho$ from  a ring $ R$ into the set $C\bigcup \infty$  is a pseudovaluation if  for every $x,y\in R$ we have

\begin{equation} \rho(x+y)\geq \mbox{min} \{\rho(x),\rho(y)\}\end{equation}
\begin{equation}  \rho(xy)\geq \rho(x)+\rho(y)\end{equation}

and $\rho(0)=\infty$.   A pseudovaluation is a valuation if relation $(1.20 )$  is an equation. {\it The pseudovaluations and valuations which will be considered throughout this paper are  discrete. We will assume also that  $\rho(x)=\infty$ only if $x=0$.} We refer the reader to Cohn's book [3], or Bourbaki [2], chapter VI,  for the  the basic concepts and properties of valuations and filtrations, and the graded rings associated with them.

The following  Theorem II  will be applied for construction of valuations and pseudovaluations in  rings $R$ with  independent polycentral systems,  and  in particular in group rings of polycyclic groups; we will use then  these pseudovaluations for    the study of $p$-series in groups. \bigskip

{\bf Theorem II.} {\it Let $R$ be a ring, $T=<t_1,t_2,\cdots,t_n> $ be a polycentral independent system in $R$ which is composed from the central systems $T_1,T_2,\cdots,T_k$, $A$ be the ideal generated by the system $T$. Let $f$ be a function on  on $T$ whose values are  natural numbers and  $f(t_1)>2 f(t_2)\,\, \mbox{for}\,\, t_1\in T_i,t_2\in T_{i+1}\,\,  (i=1,2,\cdots,k-1)$.

 Then there exists a pseudovaluation $v$ of $R$ such that }

\begin{equation} v(t)=f(t)\,\ \mbox{if}\,\, (t\in T) \end{equation}

{\it and the graded ring $gr_v(R)$ is isomorphic to the polynomial ring\\  $(R/A)[\tilde{t}_1,\tilde{t}_2,\cdots, \tilde{t}_n]$ over the zero degree component $R/A$, the topology defined in $R$ by this pseudovaluation is equivalent to the topology defined by the powers of the ideal $A$.  Furthermore,   $v$ is a unique pseudovaluation such that $v(t)=f(t)\,\, (t\in T)$  and the graded ring associated to it is isomorphic to $R[\tilde{t}_1,\tilde{t}_2,\cdots,\tilde{t}_n]$.}\bigskip

We have the following immediate corollary of Theorem II.\bigskip

{\bf Corollary 1.3.} {\it Let $R$ be a ring, $T=<t_1,t_2,\cdots, t_n>$ be a polycentral independent system in $R$ which is composed from the central systems $T_1,T_2,\cdots,T_k$, $A$ be the ideal generated by the system $T$. Let $M_i\,\, (i=1,2,\cdots,k)$ be a system of natural numbers such that $M_i>2M_{i+1}\,\, (i=1,2,\cdots,k-1)$.

 Then there exists a unique pseudovaluation $v$ of $R$ such that }

\begin{equation} v(t)=M_i\,\, \mbox{if}\,\, t\in T_i;\,\,  (i=1,2,\cdots,k) \end{equation}

{\it and the graded ring $gr_v(R)$ is isomorphic to the polynomial ring\\ $(R/A)[\tilde{t}_1,\tilde{t}_2,\cdots, \tilde{t}_n]$, the topology defined in $R$ by this pseudovaluation is equivalent to the topology defined by the powers of the ideal $A$.  }\bigskip

Theorem  II will be derived from Theorem III which is formulated and proven in section $3$.

Theorem II implies in particular that the completions $\tilde{R}_{\rho}$ and $\tilde{R}_v$ of $R$ in the $\rho$-topology and in the $v$-topology are homeomorphic. Let $\mathcal{X}$ be a system of coset representatives in the quotient ring $R/A$ for the ideal $A$. Theorems I and II  now  yield the following corollary of Theorem III.\bigskip

{\bf Corollary 3.1.} {\it An arbitrary element $x\in \tilde{R}_{\rho}\cong \tilde{R}_v$ has a unique representation $(1.18)$ 
where $\lambda_n\,\,   (n=1,2,\cdots)$ are elements from a system of coset representatives $\mathcal{X}$   of the quotient ring $R/A$ and $\pi_n\,\, (n=1,2,\cdots)$ are standard monomials on the system $T$. The length $l(\pi_n)$ and the value $v(\pi_n)$ run to infinity if $n\to\infty$.} \bigskip

The following Corollary $3.8$ of Theorems I and II will be proven in section $3$; this corollary provides a method for construction of pseudovaluations in a ring $R$ using  polycentral systems in a graded ring $gr(R)$.\bigskip

{\bf Corollary 3.9.} {\it  Let $R$ be a ring with a discrete pseudovaluation $\rho$, $gr_{\rho}(R)$  be the associated graded ring. Assume that there exists in $gr(R)$ an independent polycentral sustem  $T$, let $A$ be the ideal generated in $gr(R)$ by this system. Then  there exists in $R$ a discrete pseudovaluation $v$ such that the graded ring $gr_v(R)$ is isomorphic to a subring of the Laurent polynomial ring in the system of variables $T,t,t^{-1}$ over the ring 
$(gr(R))/A$.}\bigskip

 Theorems I-III can be applied for constructions of valuations in some classes of skew fields, mainly the universal field of fractions of the free ring $D_K<X>$ over a field $D$. \bigskip

{\bf 1.4.} Let $H$ be a torsion free polycyclic group with Hirsch number $r$ which satisfies conditions of Theorem VI, and hence has a $p$-series with unit intersection and associated graded Lie algebra $L_p(H,H_i)$ abelian of rank $r$. Theorem XII provides sufficient conditions for the existence of a $p$-series with associated algebra $L_p(H,H_i)$  free abelian of rank $r$. We will study the necessary and   sufficient conditions for the  existence in a polycyclic group of a  series with associated graded algebra  $L_p(H,H_i)$  free abelian    of rank $r$ in the  paper [11]; the results of the current paper reduce the   problem  to the case when the quotient group $H/H^{\prime}$  is a finite  $p$-group. Here we will  only prove Proposition $1.1.$ and will sketch an example which show that there are torsion free polycyclic group which   have $p$-series with unit intersection and associated graded algebra abelian of rank equal to the Hirsch rank $r$ of 
the group,  nevertheless   these groups do not have series with associated graded algebra {\it free} abelian of rank $r$.\bigskip 

{\bf Proposition 1.1.} {\it Let $H$ be a group. Assume that the quotient group $H/\gamma_p(H)$ of $H$ by the $p^{th}$ term of the low central series has exponent $p$.

 Then $H$ can not have a $p$-series $(1.1)$  with unit intersection and associated graded algebra free abelian.}\bigskip 

{\bf Proof.} Assume that there exists a $p$-series $(1.1)$ with algebra $L_p(H,H_i)$ free abelian.  Let $h$ be one of the elements with minimal weight, say $\omega(h)=k$. 

We consider now an arbitrary  commutator $u=[h_1,h_2,\cdots,h_l]$ with $l\geq p$. Since the algebra $L_p(H,H_i)$ is abelian  we obtain from Lemma $2.2.$ below   that the weight $w ([h_1,h_2]) $ of the commutator  $[h_1,h_2]$ is greater that $w(h_1)+w(h_2)$, and then that

\begin{equation} w([h_1,h_2,\cdots,h_l])>w(h_1)+w(h_2)+\cdots w(h_l)\geq lk\geq pk\end{equation}

Since the algebra $L_p(H,H_i)$ is free abelian the weight of $h^p$ must be 
$pk$  by Lemma $2.1$. On the other hand,  since $h^p$ is a product of commutators of  length greater than  or equal $p$ its weight  must be greater than $pk$ by $(1.23)$. This contradiction shows that $L_p(H,H_i)$  can not be free abelian, and the proof is complete.\bigskip 

Proposition $1.1.$ provides a method for construction groups which can not have $p$-series with unit intersection and associated graded algebra free abelian.
We will briefly sketch  here with a proof one  example, the proof and the computations will be given   in [11]; we will provide there more examples.\bigskip

{\bf Example.} Let $F$ be a free group with generators $f_1,f_2,f_3$, $N$ be a normal subgroup such that the quotient group $F/N$ has exponent $2$. Let $\bar{F}=F/N^{\prime}$. The group $\bar{F}$ is generated by the elements $\bar{f}_1,\bar{f}_2,\bar{f}_3$,  the images of $f_1,f_2,f_3$, and it is  an extension of the free abelian group $\bar{N}=N/N^{\prime}$  by the abelian group $\bar{F}/\bar{N}\cong F/N$ of exponent $2$ and rank $3$. The group $\bar{F}$  does have a $p$-series 

\begin{equation}\bar{F}=\bar{F}_1\supseteq \bar{F}_2\supseteq\cdots\end{equation}

with unit intersection and associated graded algebra free abelian of rank equal to the Hirsch rank of $\bar{F}$ which is equal to $17$,  this series can be constructed  by methods of Lichtman [10].  We consider then      the quotient group $H$ of $\bar{F}$ by the normal subgroup $U$ generated by the elements 

\begin{equation} \bar{f}_1^2[\bar{f}_3,\bar{f}_2], \bar{f}_2^2[\bar{f}_1,\bar{f}_3],\bar{f}_3^2[\bar{f}_2,\bar{f}_1]\end{equation}

We will prove  in [11] that the  group $H$ is torsion free; it is easy to verify that     $H$ is an extension of a  free abelian group by a group of exponent $2$, it is generated by the elements $a,b,c$ which are  the images 
of the elements $\bar{f}_1,\bar{f}_2,\bar{f}_3$,  subject to the relations 

\begin{equation} a^2=[b,c], b^2=[c,a],c^2=[a,b]\end{equation}

The group $H$ is a  residually \{finite $2$-group\};  
 relations $(1.26)$ imply that the quotient group $H/H^{\prime}$ has exponent $2$, we conclude now from   Proposition $1.1.$
that the group $H$ can not contain a $2$-series with unit intersection and the algebra $L_2(H,H_i)$ free abelian. 

To show that such series does not exist  if  $p\not=2$ we observe that the group $H$ is not a residually-\{finite $p$-group\} for $p\not=2$ so an arbitrary $p$-series will have in this case a non-unit intersection.\bigskip 

{\bf 1.5.} Our results on restricted Lie algebras of polycyclic groups are related to the cases when the algebra $L_p(H,H_i)$ is finitely generated and   we  deal almost entirely with the case when the rank of $L_p(H,H_i)$ coincides with the Hirsch number of $H$. We will construct in section $9$ examples of $p$-series in a free abelian group of rank $2$  with associated graded algebra free abelian of rank $1$; we will construct also examples of series with associated graded algebra infinitely generated.\bigskip

\setcounter{section}{2}
\setcounter{equation}{0}

\section*{\center \S 2.  Preliminaries.}

{\bf 2.1.} We give now a brief account  of the main concepts and results  about
the restricted Lie algebras of groups; these results are obtained in  Lazard's
article [7]. 

Let $H$ be a group, 

\begin{equation} H=H_1\supseteq H_2\supseteq\cdots\end{equation} 

be a $p$-series. The
restricted Lie algebra associated to  series $(2.1)$  is a Lie
algebra over the prime field $Z_p$ whose vector space is $\sum_{i=1}^\infty H_i/H_{i+1}$; an  element  $\tilde{h}=hH_{i+1}$ from $H_i/H_{i+1}$
is  called homogeneous of degree $i$ and the Lie operation for homogeneous
elements is defined in the following way.  If $x_1\in H_{i_1}/H_{i_1+1}$, $x_2\in H_{i_2}/H_{i_2+1}$  and   $x^{\ast}_\alpha$ is the
coset representative of $x_\alpha\,\, (\alpha =1,2)$  then  $[x_1,x_2]$  is the
element of  $ H_{i_1+i_2}/H_{i_1+i_2+1}$ which contains the commutator $[x^{\ast}_1, x^{\ast}_2]$.  This operation  extends by
distributivity on arbitrary elements from $\sum_{i=1}^\infty H_i/H_{i+1}$  and it defines the structure of graded Lie algebra  on the set $\sum_{i=1}^\infty H_i/H_{i+1}$.   If   $x\in H_i/H_{i+1}$ and $x^{\ast}$  is its
representative in $H$ then $x^{[p]}$ is defined as the homogeneous component of the element $(x^{\ast})^p$
and we obtain the structure of a restricted Lie algebra in $L_p(H,H_i)$.  We will denote by $\tilde{h}$ the homogeneous component of an element $h\in H$.

 Let $U=\bigcap_ {i=1}^{\infty}H_i$, $\bar{H}=H/U,\bar{H}_i=H_i/U\,\, (i=1,2,\cdots,)$.  The definition of the algebra $L_p(H,H_i)$ implies that there exists a natural isomorphism between the graded algebras $L_p(H,H_i)$ and  $L_p(\bar{H},\bar{H}_i)$; this fact reduces the study of the algebra $L_p(H,H_i)$ to the case when $\bigcap_{i=1}^{\infty}H_i=1$. 

If $h$ is an element of $H$ and $h\in H_i\backslash H_{i+1}$ then the weight $w(h)$ of $h$ is $i$; the weight of $1$ is $\infty$.

 Let $F$ be a  subgroup of $H$.  Then the  series $(2.1)$ induces in $F
$ a $p$-series  $F_i=F\bigcap H_i\,\, (i=1,2, \cdots)$.  We denote the associated Lie
algebra of $F$ by  $L_p(F,F_i)$ and we will use this notation throughout the paper. There is a natural imbedding of the graded Lie algebra $L_p(F,F_i)$ in $L_p(H,H_i)$; if $F$ is normal in $H$ then
$L_p(F,F_i)$ is an ideal in $L_p(H,H_i)$. Let  $\phi: H\longrightarrow H/F=G$ be an epimorphism  of groups, and let $G_i=\phi(H_i)=(H_iF)/F\,\, (i=1,2,\cdots)$. We will make a use of the following result which is Theorem $2.4.$ 
 in Lazard's article [7].\bigskip

{\bf  Proposition  2.1.} {\it Let $F$ be a normal subgroup of $H$, $G=H/F$. The epimorphism $\phi\colon H\longrightarrow G $ defines in a natural way an epimorphism $\tilde{\phi}: L_p(H,H_i)\longrightarrow L_p(G,G_i)$  which preserves the degrees of the homogeneous elements. The kernel of $\tilde{\phi}$ is  the ideal $L_p(F,F_i)$.}\bigskip

{\bf Corollary 2.1.} {\it Assume that the normal subgroup  $F$ in Proposition $2.1.$ has  finite index. Let $Q$ be the subgroup formed by   all the elements of $H$ whose images in $G$ belong to $\bigcap_{i=1}^{\infty} G_i$. Then $Q$ is a normal subgroup which contains $F$ , the index $(H:Q)$ is a power of $p$ and $L_p(Q,Q_i)=L_p(F,F_i)$.} \bigskip

{\bf Proof.} Let $G_0=\bigcap_{i=1}^{\infty} G_i$. Then  $Q$ is the inverse image of $G_0$ in $H$  and $Q\supseteq F$. The quotient algebras $L_p(H,H_i)/L_p(Q,Q_i)$ is isomorphic to $L_p(G, G_i)$, so $L_p(F,F_i)=L_p(Q,Q_i)$. 

The quotient group $\bar{H}=H/Q$  contains a $p$-series $\bar{H}_i=(H_iQ)/Q\,\, (i=1,2,\cdots)$ with unit intersection. This series has a finite length because the group $H/Q$ is finite. This implies that $H/Q$ is a finite $p$-group, and the proof is complete.\bigskip

The following result is  the first half of Theorem $2.2.$ in Lazard's article 
$[7]$. \bigskip

{\bf Proposition 2.2.} {\it Let $\phi$ be a homomorphism from a group $H$ into a group $G$. Let $H_i$ and $G_i$ be $p$-series in $G$ and $H$ respectively and assume that $\phi(H_i)\subseteq G_i\,\, (i=1,2,\cdots)$. Then $\phi$ defines in the following  way a homomorphism $\tilde{\phi}\colon L_p(H,H_i)\longrightarrow L_p(G,G_i) $ such that  if  $h\in H$ is an element of weight $k$ then}

$$\tilde{\phi}(\tilde{h})=\phi(h)+G_{k+1}$$

The following two facts   follow from this theorem immediately.\bigskip

{\bf Corollary 2.2.} {\it Assume that the conditions of Proposition $2. 2.$ hold. Let  $h$ be  an element of $H$ and $\tilde{h}$ be its homogeneous component. The element $\tilde{h}$ belongs to the kernel of $\tilde{\phi}$ iff the weight of $\phi(h)$ is greater than the weight  of  $h$.  If   $\tilde{\phi}(\tilde{h})\neq 0$ then $\tilde{\phi}(\tilde{h})=\widetilde{\phi(h)}$.  }\bigskip

{\bf Corollary 2.3.} {\it Let $\Phi$ be a group of automorphisms of $H$ such that $\phi(H_i)=H_i\,\, (i=1,2,\cdots)$. Then the group $\Phi$ defines in a natural way a group of automorphisms $\tilde{\Phi}$ of the  restricted Lie algebra $L_p(H,H_i)$; the homogeneous components  $H_i/H_{i+1}\,\, (i=1,2,\cdots)$ are $\tilde{\Phi}$-invariant.  }\bigskip

{\bf Lemma 2.1.}  {\it  Let $h$   be a non-unit  element of $H$, $\tilde{x}$ be its homogeneous component.  Assume that  $w(h)=k$. If $w(h^{p^n})=kp^n$ then 
$\widetilde{h^{p^n}}=\tilde{h}^{[p^n]}$; if $w(h)^{p^n}>p^nw(h)$ then $\tilde{h}^{[p]^n}=0$}.\bigskip

The proof is straightforward. \bigskip

{\bf Lemma 2.2.} i) { \it Let  $[h_1,h_2,\cdots,h_l]$ be a right normed commutator in $H$. If  its weight is equal to $w(h_1)+w(h_2)+\cdots+w(h_l)$ then the homogeneous component of this commutator is  $[\tilde{h}_1,\tilde{h}_2,\cdots,\tilde{h}_l]$; if its weight is greater than  $w(h_1)+w(h_2)+\cdots+w(h_l)$ then}  $[\tilde{h}_1,\tilde{h}_2,\cdots,\tilde{h}_l]=0$.

ii)  {\it Let $h_1,h_2,\cdots,h_l$ be elements of $H$ with the same weight 
$q$.  Assume that the weight of the element $h_1h_2\cdots h_l$ is also $q$. Then the homogeneous component of $h_1h_2\cdots h_l$ is} \begin{equation} 
\tilde{h}_1+\tilde{h}_2+\cdots + \tilde{h}_l\end{equation}

{\it If the weight of  $h_1h_2\cdots h_l$ is greater than $q$ then $\tilde{h}_1+\tilde{h}_2+\cdots+\tilde{h}_l=0$. }\bigskip

{\bf Proof.} All   the statements are known facts whose proofs are straightforward. For instance, the proof of  statement ii) is obtained in the following way.  If the condition of the assertion hold then the  image of $h_1h_2\cdots h_l$ in the quotient group $H_q/H_{q+1}$ is the sum of the images of the elements $h_i\,\, (i=1,2,\cdots,l)$, hence its homogeneous component is $(2.2)$. The rest of the statements are proven by similar arguments.\bigskip

{\bf Corolllary 2.4.} {\it Let $H$ be a nilpotent group. Then the restricted Lie algebra $L_p(H,H_i)$ is a nilpotent Lie algebra.}\bigskip

{\bf Proof.} Follows immediately from statement i) of Lemma $2.2.$\bigskip

{\bf 2.2.} We need an analog of Proposition $2.1.$ and Lemmas $2.1.$ and $2.2.$   for rings. Let $R$ be a ring with a non-negative pseudovaluation $v$, $R_i$ is the filtration defined by $v$, that is 

$$R_i=\{r\in R| v(r)\geq i\}\,\, (i=0,1,\cdots)$$
It is clear that the pseudovaluation function $v$ is completely defined if the filtration $R_i$ is given. Let $A$ be an ideal in $R$, we have in $A$ an induced filtration $A_i=A\bigcap R_i\,\, (i=0,1,\cdots)$, let $gr(R)$ and $gr(A)$ be the graded ring of $R$ and $A$ associated to the filtrations $R_i, A_i\,\, (i=0,1,\cdots)$ respectively. Further if $\bar{X}$ denotes the image of a subset $X\in R$ under the natural homomorphism $\phi\colon R\longrightarrow R/A$ we obtained in $\bar{R}$ a filtration $\bar{R}_i\,\, (i=0,1,\cdots)$ and a pseudovaluation $\bar{v}$ defined by this filtration. It is natural to say that  pseudovaluation $\bar{v}$ is obtained from $v$ by the homomorphism $\phi$. The following fact is a modified  version of Proposition $2$ in section $3.4.$ of Bourbaki [2]; its proof can be read off from [2] or obtained by a straightforward argument.\bigskip

{\bf  Proposition  2.3.} {\it The homomorphism $\phi\colon  R\longrightarrow R/A$   defines in $\bar{R}$  a filtration $\bar{R}_i\,\, (i=0,1,\cdots)$,  a pseudovaluation $\bar{v}$ defined by this filtration and   a homomorphism of graded rings $gr(\phi)\colon gr(R)\longrightarrow gr(\bar{R})$ with 
kernel $gr(A)$.}\bigskip

 {\bf Proposition 2.4.} {\it Let $R$ be a ring, $v$ and $w$ be two non-negative pseudovaluations  which define equivalent topologies in $R$. Assume that there exists a system of elements $T\subseteq R$ such that $v(t)=w(t)\,\, (t\in T)$ and the graded rings $gr_v(R)$ and $ gr_w(R)$ are isomorphic to the polynomial rings  $(R/A)[\tilde{T}_v]$ and $(R/A)[\tilde{T}_w]$ respectively, where $A$ is the ideal generated by $T$ and $\tilde{T}$ is the set of the homogeneous components of elements $t\in T$.  Then
$v(r)=w(r)$ for every $r\in R$.} \bigskip

{\bf Proof.} Let $\mathcal{X}$ be a system of coset representatives for the quotient ring $R/A$,   $\pi_{\alpha}\,\, (\alpha=1,2,\cdots,n)$ be distinct standard monomials with $v$-value $k$. Then 

\begin{equation} v(\sum_{\alpha=1}^n\lambda_{\alpha}\pi_{\alpha})=k \end{equation} 

 and 

\begin{equation}w(\sum_{\alpha=1}^n\lambda_{\alpha}\pi_{\alpha})=k\end{equation} 

 because of the condition $gr_v(R)\cong gr_w(R)\cong (R/A)[\tilde{T}]$; we obtain in particular that $w(\pi)=v(\pi)$ for an arbitrary standard monomial.

Now assume that there exists an element $r\in R$ such that $v(r)=k$ and $w(r)=l>k$. Let $A_i, B_j\,\, (i=0,1,\cdots; j=0,1,\cdots)$ be the filtration defined in $R$ by the pseudovaluations $v$ and $w$ repectively. Since the topologies defined by $v$ and $W$ are equivalent we can find $j>k$ such that 
$A_j\subseteq B_{l+1}$.  Let $\bar{R}=R/A_j$. We have in $\bar{R}$ pseudovaluations  $\bar{v}$ and $\bar{w}$; since the $v$-values and $w$-values of elements from $ A_j$ are greater then $l$ we obtain that 

\begin{equation} \bar{v}(\bar{r})=k, \bar{w}(\bar{r})=l\end{equation}

The ring $\bar{R}$ has now a pseudovaluation $\bar{v}$ defined by the filtration $\bar{A}_i=A_i/A_j\,\, (i=1,2,\cdots,j$, the graded ring 
$gr_{\bar{v}}(\bar{R})$ is an isomorphic image of the polynomial ring 
$(R/A)[\tilde{T}_v]$,  it is generated over the subring $R/A$ by all the standard monomials whose  $v$-value less than or equal to $j-1$ and  the standard monomials with $v$-value $i$ form a basis of the homogeneous component $\bar{A}_i/\bar{A}_{i+1}\,\, (i=1,2,\cdots, j-1)$.  A routine argument shows that every element $\bar{x}\in \bar{R}$ has a unique representation

\begin{equation} \bar{x}=\sum_{\alpha=0}^m\mu_{\alpha}\pi_{\alpha}\,\,\end{equation}

where $ \mu_{\alpha}\in \mathcal{X}\,\, (\alpha=1,2,\cdots,m)$ and $v(\pi_{\alpha})\leq j-1\,\, (\alpha=1,2,\cdots,m)$. 

Since $\bar{v}(\bar{r})=k$  we conclude that $\bar{r}$ has a unique representation 

\begin{equation} \bar{r}=\sum_{\beta}\mu_{\beta}\pi_{\beta}\end{equation}

 where all the standard monomials $\pi_{\beta}$ have $v$-value greater than or equal to $k$.

We consider now the following element $r_1\in R$ 

\begin{equation} r_1=\sum_{\beta}\mu_{\beta}\pi_{\beta}\end{equation}

The image of $r_1$ in $\bar{R}$ is $\bar{r}$; on the other hand  representation  $(2.8)$ implies that $v(r_1)=w(r_1)=k$. Since all the elements of $A_j$ have $v$-values and $w$-values greater than $k$ we obtain that 
$v(\bar{r})=w(\bar{r})=k$. 
We obtained a  contradiction with  relation $(2.5.)$ and the asssertion follows.\bigskip

The proof of  following  analog of Lemma $2.2.$   for rings is obtained by the same straightforward argument.\bigskip

{\bf Lemma 2.3.} {\it Let $R$ be a ring with a pseudovaluation $\rho$ and associated graded ring $gr(R)$, $r_1,r_2,\cdots, r_l$ be elements of $R$ with the same weight $q$. }

i) {\it If the weight of the element $r_1+r_2+\cdots +r_l$
 is $q$ then its homogeneous component in $gr(R)$ is equal to $\tilde{r}_1+\tilde{r}_2+\cdots \tilde{r}_l$;  if the weight of element $r_1+r_2+\cdots +r_l$
is greater than $q$ then} $\tilde{r}_1+\tilde{r}_2+\cdots +\tilde{r}_l=0$.

ii) {\it If the weight of the element $r_1 r_2 \cdots r_l$ is $lq$ then its homogeneous component is $\tilde{r}_1\tilde{r}_2\cdots \tilde{r}_l$; if the weight of this element is greater than $lq$ then $\tilde{r}_1\tilde{r}_2\cdots \tilde{r}_l=0$.} \bigskip

{\bf 2.3.} {\bf Lemma 2.4.} {\it Assume that the algebra $L_p(H,H_i)$ is generated by the first $l$ homogeneous components. Let $h\in H$ be an element of weight $r$. Then the homogeneous component $\tilde{h}$  is a sum of Lie monomials
 \begin{equation}[\tilde{h}_{\alpha_1},\tilde{h}_{\alpha_2},\cdots,\tilde{h}_{\alpha_s}]^{[p]^{n_{\alpha}}}\end{equation} 
where the homogeneous elements  $\tilde{h}_{\alpha_1},
\tilde{h}_{\alpha_2},\cdots, \tilde{h}_{\alpha_s} $ are taken from the first $l$ factors $H_i/H_{i+1}$, the weight of every monomial $ (2.9)$ is  $r$, and} 

\begin{equation} 
w([\tilde{h}_{\alpha_i},\tilde{h}_{\alpha_2},\cdots,\tilde{h}_{{\alpha}_s}])=    w(\tilde{h}_{\alpha_1})+ w(\tilde{h}_{\alpha_2})+\cdots+w(\tilde{h}_{\alpha_s})=r_1  \end{equation}

{\it where $r=p^{n_{\alpha}}r_1$}. 

{\it Further, if $h_{\alpha_i}$ is the coset representative of $\tilde{h}_{\alpha_i}\,\, (i=1,2,\cdots,s)$ then the element } 

\begin{equation} [h_{\alpha_1},h_{\alpha_2}\cdots, h_{\alpha_s}]^{p^{n_{\alpha}}}\end{equation} 

{\it is the representative in $H$ of the homogeneous component $(2.9)$}.\bigskip

{\bf Proof.} Since the algebra $L_p(H,H_i)$ is graded the element $\tilde{h}$ must be a sum of homogeneous Lie monomials $(2.9)$ of degree $r$.  Lemma $2.2.$ implies that if a  Lie monomial $(2.9)$ is non-zero then condition $(2.10)$ holds, and that the element $(2.11)$ is a coset representative of the element $(2.9)$  in $H_r/H_{r+1}$, and the proof is complete.\bigskip

We will make an essential  use of the following fact which is proven in 
Lichtman [10]. (See [10], Lemmas $2.6.$ and  $3.1.$)\bigskip

{\bf Proposition 2.5. } {\it Let $H$ be a polycyclic group with Hirsch number $r$. Assume that there exists a $p$-series $(1.1)$    such that the Lie algebra $L_p(H,H_i)$ is abelian. Then   the rank of the algebra $L_p(H,H_i)$ does not exceed $r$; if this rank is equal   $r$ then the subgroup $\bigcap_{i=1}^{\infty}H_i$ is finite. }\bigskip

We will need also the following result. \bigskip
 
{\bf Proposition 2.6.} {\it Let $H$ be a  polycyclic group with Hirsch number $r$, $F$ be an   normal subgroup   with a Hirsch number $r_1$. Assume that there exists a $p$-series with associated graded algebra $L_p(H,H_i)$ abelian of rank $r$. 

Then the  Lie algebra $L_p(F,F_i)$ associated to the $p$-series $F_i=F\bigcap H_i\,\, (i=1,2,\cdots)$ is abelian of rank $r_1$. The intersection of the $p$-series $\bar{H}_i=(H_iF)/H_i\,\, (i=1,2,\cdots)$ is a finite subgroup and the Lie algebra $L_p(\bar{H},\bar{H}_i)$ of the group $\bar{H}=H/F$ associated to the $p$-series $\bar{H}_i\,\, (i=1,2,\cdots)$  is abelian of rank $r-r_1$.} \bigskip

{\bf Proof.}  We can assume that $H$ is infinite. Let $U$ be a torsion free normal subgroup of finite index, $V=U\bigcap F$.    The  subagebra $L_p(U,U_i)$ associated to the $p$-series $U_i=H_i\bigcap U\,\, (i=1,2,\cdots)$ has a finite index in  $L_p(H,H_i)$ so its rank is equal to $r$; similarly the rank of the subalgebra $L_p(V,V_i)$ is equal to the rank of $L_p(F,F_i)$. We see that we can assume that the group $U$ is torsion free. The first statement follows now from Proposition $3.5.$ in Lichtman [10]; the second statement follows from Proposition $3.2.$ in [10].\bigskip

{\bf 2.4.}  The study of $p$-series in groups and  Lie algebras associated to them  is connected to filtrations and valuations in the group rings of these groups. We have already observed that in the study of the properties of the algebra $L_p(H,H_i)$ we can assume that $\bigcap_{i=1}^{\infty} H_i=1$ and   $p$-series $(1.1)$  defines
in a natural way a weight function in the group $H$ and  this weight function defines a filtration in the group ring $KH$ (see Passman, [13]):
\begin{equation}A_0 = KH\supseteq A_1\supseteq A_2\supseteq\cdots\end{equation}
where $ A_n\,\,(n\geq 1)$ is the $K$-linear span of the set of all the products
\begin{equation} \{(x_{\alpha_1} - 1) (x_{\alpha_2} -1) \cdots (x_{\alpha_s} - 1)\mid\sum_{i
=1}^s w(x_{\alpha_i}) \ge n\} \end{equation}
 We recall that if $H_i=M_i(KH)\,\,(i=1,2,\cdots)$ is the Lazard-Zassenhaus-Jennings $p$-series (see Passman [13], section $11$) then  the filtration defined by it is $\omega^n(KH)\,\, (n=1,2,\cdots)$ where $\omega(KH)$ is the augmentation ideal of $KH$.  Series $(2.1)$ defines the $p$-topology in the group $H$ and  filtration $(2.12)$ defines a topology in the group ring $KH$ which   we will call   $p$-topology in $KH$.  The following fact follows easily from the classical theorems on group rings, we give a sketch of the proof.\bigskip

{\bf Lemma 2.6.} {\it If the topology defined by $p$-series $(2.1)$ in the group $H$ is equivalent to the $p$-topology then the topology defined by the filtration $(2.12)$ in the group ring $KH$ is equivalent to the $p$-topology of $KH$.}\bigskip

{\bf Proof.} Let $\omega^k(KH)$ be a given power of the augmentation ideal. We have to find $i=i(k)$ such that $ A_i\subseteq \omega^n(KH)$. Since $(h-1)\in \omega^n(KH)\,\, (h\in M_n(H))$  the proof is easily reduced to the case  of the group $\bar{H}=H/M_n(H)$ and the $p$-series $\bar{H}_i=(H_iM_n(H))/M_n(H)\,\, (i=1,2,\cdots)$;    the series $\bar{H}_i\,\, (i=1,2,\cdots)$ has a finite length because there  exists a number $j$ such that $H_j\subseteq M_n(H)$. We see that we can assume that the original series $H_i\,\, (i=1,2,\cdots)$ has 
a finite length.

We obtain now that the weights of the factors $x_{\alpha}-1$ in the left side of $(2.13)$ are bounded; this implies   that if  the number $n$ in $(2.13)$ runs to infinity  then  $s\longrightarrow\infty$.  We obtain from this that if $n$ is sufficiently large then $s\geq k$ and $A_i\subseteq \omega^k(KH)$ if $i\geq s$.\bigskip 

 If  $\bigcap_{n=1}^{\infty}A_n=0$  the  filtration $(2.12)$ defines in a natural way a pseudovaluation $\rho$ in $KH$: if $x\not=0$ then $\rho(x)$ is equal to the maximal $n$ such that $x\in A_n$, and $\rho(0)=\infty$. We will say that this pseudovaluation is defined by $p$-series $(2.1)$. We recall that this  pseudovaluation is a valuation if the graded ring $gr(KH)$ is a domain; in the case of the   group ring  $KH$ this condition can be replaced by the equivalent condition that $U_p((L_p(H, H_i))$ is a domain (see Proposition $2.7.$ below). Further it is known that  if filtration $(2.12)$ defines a valuation in $KH$   the topological completion of $KH$ must be a domain.

We will need the following  fact (see Lichtman [10], Lemma $3.2.$)\bigskip 

{\bf Lemma 2.7.} {\it Let $H$ be a group which contains a $p$-series $(2.1)$ with unit intersection, $N$ be a normal subgroup of $H$. Assume that the topology defined by the  series $N_i=N\bigcap H_i\,\, (i=1,2,\cdots)$ in $N$ is equivalent to the  $p$-topology in $N$, and the topology defined by the series $\bar{H}_i=(H_iN)/N\,\, (i=1,2,\cdots)$  in the quotient group   $\bar{H}=H/N$ is equivalent to the  $p$-topology in $\bar{H}$. Then the topology defined in $H$ by series $(2.1)$ is equivalent to the $p$-topology in $H$.} \bigskip

{\bf Lemma 2.8.} {\it Let $H$ be a free abelian   group  of rank $r\geq 1$ which contains a $p$-series $(2.1.)$ with unit intersection and the associated graded  algebra 
$L_p(H,H_i)$  finitely generated. }

i) {\it If the rank of $H$ is $1$  then so is the rank of $L_p(H,H_i)$}.

ii) {\it If the rank of  $L_p(H,H_i)$ is equal to $r$ then the topology defined by series $(2.1)$ is equivalent to the $p$-topology.}\bigskip

{\bf Proof.}
We will prove  first both statements for the case when $H$ has rank $1$.
 Since $L_p(H,H_i)$  is finitely generated we can find a number $i_0$ such that the subalgebra 

\begin{equation}\sum_{i\geq i_0}H_i/H_{i+1}\end{equation}

contains no nilpotent elements. Hence this subalgebra is free abelian of rank $1$. 

On the other hand the subalgebra $(2.14)$ is the restricted Lie algebra of the subgroup $V=H_{i_0}$ which corresponds to the $p$-series 

\begin{equation} V=H_{i_0}\supseteq H_{i_0+1}\supseteq\cdots \end{equation}

Let  $v$ be  the generator of $V$, $\tilde{v}$ be the homogeneous component of $v$. The weight of $v$ is $i_0$, since $\tilde{v}$ is not nilpotent we obtain from Lemma $2.2.$ that the homogeneous component of the element $v^{p^n}$ is equal to $\tilde{v}^{[p]^n}$ and the weight of $v^{p^n}$ is $p^n i_0$. We see that the topology defined in $V$ by series $(2.14)$ is equivalent to the $p$-topology. Since the quotient group $H/H_{i_0}$ is a  finite cyclic $p$-group we conclude  from Lemma $2.7.$ that the topology defined in $H$ by series $(2.1)$  is equal  to the $p$-topology. This completes the proof for the special case when the rank of $H$ is $1$.

We consider now the general case.  Let $u=h_1,h_2,\cdots,h_r$ be a free system of generators for $H$, $U$ be the subgroup generated  by the element $u$. Proposition  $2.6.$ implies that the algebra $L_p(U,U_i)$ associated to the $p$-series $U_i=H \bigcap U_i\,\, (i=1,2,\cdots)$ is finitely generated  abelian of rank $1$, that the $p$-series $\bar{H}_i=(H_iU)/U\,\, (i=1,2,\cdots)$ in the group $\bar{H}=H/U$ has a unit intesection and the algebra  $L_p(\bar{H},\bar{H}_i)\,\, (i=1,2,\cdots)$ is finitely generated abelian of rank $r-1$. We have 
already proven that the topology defined in $U$  by the series $U_i\,\, (i=1,2,\cdots)$ is equivalent to the $p$-topology, and the assumption of the induction implies that the topology defined in $\bar{H}$ by the series $\bar{H}_i\,\, (i=1,2,\cdots)$ is equivalent to the $p$-topology. The assertion now follows from Lemma $2.7.$\bigskip 

{\bf 2.5.} Let $H$ be a group, $N$ be its normal  $p$-subgroup and
$V\supseteq N^{\prime}N^p$ be an $H$-invariant subgroup of $H$, $G=H/N$.  The
conjugation in $H$ defines in a natural way  a structure of a $Z_pH$ module in the quotient group $\bar{N}=N/V$; this module is in fact a $Z_pG$-module.\bigskip

{\bf Lemma 2.9.} {\it Assume that  the module $\bar{N}$ has  dimension  less than or equal  $r$. Then the following conditions are equivalent}

i) {\it For every $h\in H$ there exists a number $n(h)$ such that} 
$[h^{p^{n(h)}},N]\subseteq V$

ii)  $[h^{p^r},N]\subseteq V$

iii) $\bigcap_{i=1}^{\infty}\omega^i(Z_pG)^i\bullet \bar{N}=0$

iv) $\omega^r(Z_pG)\bullet \bar{N}=0$ \bigskip 

{\bf Proof.} See [10], Lemma $2.15$. \bigskip

{\bf Corollary 2.5.}  {\it Assume that  $\bar{H}=H/N^{\prime}N^p$   is a residually \{finite $p$-group\} and the rank of the vector space 
$\bar{N}=N/N^{\prime}N^p$ does not exceed $r$. Then}

 \begin{equation}\omega^r(Z_pG)\bullet \bar{N}=0\end{equation}

{\it and}

 \begin{equation} [h^{p^r},N]\subseteq N^{\prime}N^p\end{equation}

{\bf Proof.} If  $\bar{H}$ is  a residually \{finite $p$-group\} then condition iii) of Lemma $2.9.$ holds and the assertion now follows from Lemma $2.8.$\bigskip

The folowing lemma  is a known facts,  its   proof is   obtained by a routine  argument.\bigskip

{\bf Lemma  2.10.} {\it Let $H$ be a group which contains a finite normal $p$-subgroup $U$ such that that quotient group  $H/U$  is an
extension of a finitely generated torsion free nilpotent group by a finite $p$-group. Assume also that there exists a number $l$ such that every $h^{p^l}\,\, (h\in H)$ centralizes the factor $U/U^{\prime}U^p$. Then $H$ contains a finitely generated torsion free nilpotent normal subgroup $V$ whose index is a power of $p$.}\bigskip

{\bf 2.6.} We will need  the following fact which follows immediately from Theorem $2.4.$ in  Lichtman [12]; this fact  is a generalization of Quillen's theorem [15].\bigskip

{\bf  Proposition  2.7.}  {\it Let  H be a group,   $(2.1)$ be an arbitrary $p$-series in $H$, $K$ be a commutative domain of characteristic $p$ and  $(2.13)$ be the filtration defined in $KH$ by this $p$-series. Let $gr(KH)$ be the graded ring associated to this filtration. Let $h\in H$ be an element of weight $i$. Then the correspondence 

\begin{equation}hH_{i+1} \longrightarrow (h-1)+A_{i+1}\end{equation} 

defines an isomorphism between
the  graded algebras $gr(KH)$  and \newline $K\otimes U_p(L_p(H, H_i))$;
in particular the restricted Lie subalgebra generated in $gr(KH)$ by all the elements $h-1+ A_{i+1}(KH)$ is isomorphic to the restricted Lie algebra $L_p(H,H_i)$. If $U$ is a normal subgroup of $H$ then the elements $u-1+A_{i+1}\,\, (u\in U)$  generate a restricted Lie subalgebra isomorphic to the subalgebra $L_p(U,U_i)$. Further, the subring $gr(KU)$ is isomorphic to the universal $p$-envelope $U_p(L_p(U,U_i)$ and the algebra $gr(KH)$ is isomorphic to  a suitable smashed product of $ gr(KU)$ with the restricted Lie algebra  of the group $H/U$ associated to the $p$-series $(H_iU)/U\,\, (i=1,2,\cdots)$.}\bigskip

Let $R$ be an algebra over a field $K$, $v$ be a pseudovaluation in $R$, $R_i\,\, (i\in Z)$ be the filtration defined by $v$, i.e. $R_i=\{r\in R| v(r)\geq i\}$, $gr(R)$ be the associated graded ring.  We extend $v$ in a natural way to the Laurent polynomial ring $R[t,t^{-1}]$ assuming that $v(t)=1$. Let $V$ be the valuation ring of  $R[t,t^{-1}]$, i.e. $V=\{x\in R[t,t^{-1}]| v(x)\geq 0$.  The following fact is Lemma $4.3.$ and Corollary $4.1.$ in Lichtman [9].\bigskip

{\bf Proposition 2.8.} {\it There exists an isomorphism $\psi$ between the rings $gr(R)$ and $V/(t)$.  If $e_j\,\, (j\in J)$ is a system of elements of $R_i$ which gives a basis of $R_i/R_{i+1}$ then the images of the elements $e_jt^{-i}\,\, (j\in J)$ in $V/(t)$ form a basis of the subspace $\psi(R_i/R_{i+1})$ of $V/(t)$.}  \bigskip

We consider once again now a group $H$ with a $p$-series $(2.1)$, let $v$ be the pseudovaluation defined in $KH$ by this series, we have also the corresponding filtration $(2.12)$. Let $h_j$ be an arbitrary element of $H$ with weight $n_j$ so  $v(h_j-1)=n_j$, let $\tilde{h}_j$ be the homogeneous component of $h_j$ in the algebra $L_p(H,H_i)$. All the elements of $KH$ have non-negative values and we see that   the quotient ring $V/(t)$ is generated over $K$ by the images $\overline{(h-1)t^{-n_j}}\,\, (j\in J)$ of the elements 
$(h-1)t^{-n_j}\,\, (j\in J)$. We apply now Propositions $2.7.$ and $2.8.$ and obtain  the following representation for the algebra $L_p(H,H_i)$.\bigskip

{\bf Proposition 2.9.} {\it  There exists an isomorphism   $\theta\colon V/(t)\cong U_p(L_p(H,H_i))$ defined by the map $\overline{(h_j-1)t^{-n_j}}\longrightarrow \tilde{h}_j\,\, (j\in J)$.  The elements $\overline{(h_j-1)t^{-n_j}}\,\, (j\in J)$ in $V/(t)$ generate with respect to the Lie operations in $V/(t)$ a subalgebra isomorphic to $L_p(H,H_i)$.  } \bigskip

Now let  $v$ be a discrete pseudovaluation in the group ring $KH$ where $K$ is a field of characteristic $p$ such that $v(h-1)\geq 1$ for every $h\in H$. This pseudovaluation defines a filtration in
$KH$
\begin{equation}A_n(KH)=\{x\in KH|v(x)\geq n\}\,\, (n=0,1,\cdots)\end{equation} and  also  a $p$-series in $H$
\begin{equation} H_i=\{h\in H|v(h-1)\geq i\}\,\, (i=1,2,\cdots)\end{equation}
Let $\widetilde{(h-1)}$ be the homogeneous component of the element $h-1$ in the graded ring associated to the filtration $(2.19)$.
The following fact is a special case of Theorem $6.6.$ in Lazard [7].\bigskip

{\bf Proposition 2.10.} {\it The elements $\widetilde{(h-1)}\,\, (h\in H)$ generate in $gr_v(KH)$ a graded Lie subalgebra
isomorphic to the Lie algebra $L_p(H,H_i)$.}\bigskip

{\bf Corollary 2.6.} {\it If the pseudovaluation $v$ in Proposition $2.10.$ is a valuation then the algebra $L_p(H,H_i)$ contains no nilpotent elements.}\bigskip

{\bf Proof.} Follows from the fact that an element $\widetilde{(h-1)}\in gr(KH)\,\, (h\in H)$ is not nilpotent.\bigskip

Now pick  a natural number $n$ and consider a discrete pseudovaluation $v_1$ in $KH$ which is defined as $v_1(x)=nv(x)\,\, (x\in KH)$; clearly, $v_1$ is equivalent to $v$,  the filtrations defined by $v$ and $v_1$ coincide and we obtain the same graded ring $gr(KH)$ for both pseudovaluations.  We obtain from this the following fact which will be used in the proof of Theorem XII.\bigskip

{\bf Corollary 2.7.} {\it The pseudovaluations $v$ and $v_1$ in $KH$ define  in $H$ the same $p$-series $H_i\,\, (i=1,2,\cdots)$ and define the same algebra $L_p(H,H_i)$.} \bigskip

{\bf 2.7.}  We will need  a few  facts about free abelian graded Lie algebras. We assume that the algebras are graded by a set $I$, which is either  the set of natural numbers or   the set of integers modulo some number $m$. \bigskip 

{\bf Lemma 2.11.} {\it Let $F$ be a  restricted graded abelian Lie algebra without nilpotent elements over a field $K$. Assume that }

\begin{equation}F=\bigoplus_{i\in I} F_i\end{equation}

{\it is a  grading in $F$. Then the algebra  $F$ is free abelian. Further, there exists a subset $I_1\subseteq I$ such that} 

\begin{equation}F^{[p]}=\bigoplus_{i\in {I_1}} F_i\end{equation}

{\it and $F$ has a direct sum representation} 

\begin{equation} F=M+F^{[p]}\end{equation}

{\it where $M$ is  graded subspace of $F$ and}
 
\begin{equation}M=\bigoplus_{i\in I_2} F_i\end{equation}

{\it where $I_2$ is the complement of $I_1$ in $I$. Let $E_i$ be  system of homogeneous elements of $F$  which forms a basis of the vector subspace 
$F_i\, (i\in I_2)$ and let $E=\bigcup_{i\in I_2 }E_i$. Then $E$  is a free system of generators for $F$,   i.e. the subset }

\begin{equation} \bigcup_{n=1}^{\infty}E^{[p]^n}\end{equation}

{\it is a basis of $F$ over $K$.  The universal $p$-envelope of $F$ is isomorphic to the symmetric algebra $K[M]$ which is the polynomial algebra over $K$ in the system of variables $E$.}  \bigskip

{\bf Proof.}  The map $x\longrightarrow x^{[p]}\,\, (x\in F)$ defines an epimorphism of  graded algebras $F\longrightarrow F^{[p]}$.  Since $F$ contains no nilpotent elements we obtain immediately that the images of distinct homogeneous components $F_{i_1}$ and $F_{i_2}$ are distinct homogeneous components $F_{i_1}^{[p]}$ and $F_{i_2}^{[p]}$ in $F^{[p]}$ and relations $(2.23)-(2.25) $ follow easily.

Since $E$ is a basis for $M$ we obtain that it is a free sysytem of generators for $F$. This completes the proof.\bigskip

  The  subspace $M\cong \bar{F}$  is in fact a graded subspace of $F$  and it  is natural to say that the free  restricted algebra $F$ is freely generated  by  the vector subspace $M$, and the universal $p$-envelope of $F$ is isomorphic to the symmetric algebra $K[M]$. Every subspace $F_i\,\, (i\in I_2)$ generate an ideal $K[F_i]$.

We will need a refinement  of Lemma  $2.11.$  for the case when $F$ is a restricted Lie algebra of a group $H$.\bigskip

{\bf Corollary 2.8.} {\it Let $H$ be a group which contains a $p$-series $(2.1)$ with unit intersection  such that  the algebra $L_p(H,H_i)$ free abelian. Let $E$ be a system of elements obtained in Lemma $2.11.$ Then every element of $\in E$ is a homogeneous component of an element $h\in H$ that is $e=\tilde{h}$ for a suitable $h\in H$. Moreover, every non-zero element  of the subspace 
$F_i\,\, (i\in I_2) $ is a homogeneous component in the algebra  
$L_p(H,H_i)$.}\bigskip

{\bf Proof.} The first statement follows from the definition of $E$. 

We prove the second statement. Let $i\in I_2$ and $0\not=h\in F_i$. Then 
\begin{equation} x=\sum_{j=1}^n\lambda_je_j\end{equation}

where $e_j\in E_i\,\, (j=1,2,\cdots,n)$ are homogeneous elements of weight $i$,
 $e_j=\tilde{h}_j\,\, (j=1,2,\cdots,n)$
and $\lambda_j \in Z_p\,\, (j=1,2,\cdots,n)$. Let $1\leq n_j\leq p-1$ be a natural number which is a coset representative of $\lambda_j\,\, (j=1,2,\cdots,n)$.  We obtain now from Lemma $2.3.$ ii)  that the weight of the element $h=\prod_{j=1}^n h_j^{n_j}$ is  $i$ and it is a coset representative for the the element $x$, that is $x=\tilde{h}$, and the assertion follows. \bigskip

{\bf 2.8.} Let $R$ be a ring, $\phi$ be an  automorphism of $R$. A  pseudovaluation $\rho$ is $\phi$-invariant if for every $r\in R$

\begin{equation} \rho(\phi(r))=\rho(r)\end{equation}

If the pseudovaluation $\rho$ is $\phi$-invariant the automorphism $\phi$ defines in a natural way an automorphism $\tilde{\phi}$ of the ring $gr(R)$. If $\tilde{r}$ is a homogeneous element of degree $k$ in $gr(R)$ and $r$ is an arbitrary   element of $R$ with $\rho$-value $k$ then $\tilde{\phi}(\tilde{r})=0$ if $\rho(\phi(r))>k$ and

\begin{equation} \tilde{\phi}(\tilde{r})=\widetilde{\phi(r)}\end{equation}

if $\rho(\phi(r))=k$. We have also 

 \begin{equation} \tilde{\phi}(\sum_{i=1}^n\tilde{r}_i)=\sum_{i=1}^n\tilde{\phi}(\tilde{r}_i)\end{equation}

for homogeneous elements $\tilde{r}_i\,\, (i=1,2,\cdots,n)$.

We will say that an automorphism $\phi$ of $R$ centralizes the ring $gr(R)$ if
 it  acts trivially on $gr(R)$. This is equivalent to the condition

\begin{equation}\tilde{\phi}(\tilde{r})=\tilde{r}\end{equation}

for all the homogeneous elements $\tilde{r}\in gr(R)$  and it  means 
 that for every $r\not=0$ there exists $u\in R$ such that
 
\begin{equation}\phi(r)=r+u\,\, \mbox{where}\,\, \rho(u)>\rho(r)\end{equation}  

or that 

\begin{equation} \rho(\phi(r)-r)>\rho(r)\,\, (r\not=0)\end{equation} 

We will need the version of these relations for group rings.\bigskip

 {\bf Lemma 2.12.} {\it Let $H$ be a group,  $\phi$ be an automorphism of $H$  and assume that a filtration and a pseudovaluation $\rho$ of  $KH$ are defined by  a  $\phi$-invariant $p$-series $(2.1)$. The automorphism $\phi$ centralizes the graded ring $gr_{\rho}(KH)$ iff it centralizes all the factors $H_i/H_{i+1}\,\, (i=1,2,\cdots)$ of the series i.e. if $h\in H_i\backslash H_{i+1}$ then there exists $u\in H_{i+1}$ such that $\phi(h)=hu$.}\bigskip

{\bf Proof.} If $h$ and $u$ are as in the condition of the lemma then relation

\begin{equation} hu-1=(h-1)+(u-1)+(h-1)(u-1)\end{equation}

implies that the homogeneous components of the elements $h-1$ and $hu-1$ coincide,  hence the homogeneous components of the elements $(h-1)$ and $\phi(h-1)$ are equal iff the condition $\phi(h)=hu$ holds, and the assertion follows.\bigskip

{\bf Corollary 2.9.} {\it The automorphism $\phi$ centralizes the graded ring $gr_{\rho}(KH)$ iff it centralizes the algebra $L_p(H,H_i)$.}\bigskip

\setcounter{section}{3}
\setcounter{equation}{0}

\section*{\center \S 3.  The Lifting of the   Valuations.}

{\bf 3.1.}  Let $R$ be a ring, $A$ be a residually nilpotent  ideal  
\begin{equation} \bigcap_{i=1}^{\infty}A^i=0\end{equation} 

and $gr(R)$ be the graded ring associated to the filtration  
\begin{equation} R\supseteq A\supseteq A^2\supseteq\cdots\end{equation} 

   Let  $\rho$ be the pseudovaluation defined by filtration $(3.2)$,  $\tilde{R}$ be  the completion of $R$ in the topology defined by  $\rho$  and $\mathcal{X}$ be a system of coset representatives for the quotient ring $R/A$. Now assume that $A$ is generated by a central element $t$  then   an arbitrary  element  $r\in \tilde{R}$ a representation 

\begin{equation} r=\sum_{j=0}^{\infty}\lambda_jt^j\,\, (\lambda_j\in \mathcal{X}; \,\,(j=0,1,\cdots) \end{equation}

   In this notation we have the following simple  fact.\bigskip

{\bf Lemma 3.1.}{ \it  Let $A$ be a residually nilpotent ideal of $R$ generated by a central element $t$.   The following three  conditions are equivalent:

$1)$  The graded ring $gr(R)$ associated to the pseudovaluation $\rho$  is isomorphic to the polynomial ring $(R/A)[t]$.

$2)$ Representation $(3.3)$ is unique}. 

$3)$ {\it The element $t$ is regular.}\bigskip

{\bf Proof.} If condition $1)$ hold then every element $x$ in the factor $A^n/A^{n+1}$ has a unique representation $x=\lambda t^n\,\, 
(\lambda \in \mathcal{X})$.  We conclude from this that if we take an   arbitrary   natural number $n$ and consider the   quotient ring $R/A^n$ then  every element of $R/A^n$  
has a unique representation

 \begin{equation}x=\sum_{j=0}^m\lambda_j t^j\end{equation}
where  $\lambda_j\in \mathcal{X}\,\, (j=1,2,\dots,m)$. Since the ring $\tilde{R}$ is the  inverse limit of the system of rings  $R/A^n\,\, (i=1,2,\cdots)$ we obtain that every element of $\tilde{R}$ has a unique  representation $(3.3)$.
 
Assume that condition $2)$ holds. We derive from this that every element of $R$ modulo $A^n$ has a unique representation as a polynomial  of degree less than or equal $n$  with coefficients from $\mathcal{X}$ and hence every element of $A^n/A^{n+1}$ has a unique representation $x=\lambda t^n\,\, (\lambda\in \mathcal{X})$.

We pick now elements $t^k$ and $t^l$. The element $t^k t^l$ has a unique representation as   a  monomial $t^{k+l}$. Hence the associated graded ring $gr(\tilde{R})$ is isomorphic to the polynomial ring $(R/A)[t]$ and the same is true for the ring $gr(R)\cong gr(\tilde{R})$. This proves that $2)\longrightarrow 1)$.

We prove now the equivalence of $2)$ and $3)$. If $2)$ holds then the element $t$ is regular in $R$ because its homogeneous component is regular in $gr(R)$. Conversely, assume that $t$ is regular and let $x\in A^n$. There exists $\lambda\in R$ such that $x=\lambda t^n$. If now $x\not \in A^{n+1}$ then $\lambda\not \in A$; further we can assume that $\lambda\in \mathcal{X} $. We see that for every element  $A^n/A^{n+1}$ there exists a representation $x=\lambda t^n\,\, (\lambda\in \mathcal{X})$. We will now  verify that if  $\lambda_1\not=\lambda_2$ then $\lambda_1 t^n\not=\lambda_2 t^n\,\, (mod \,A^{n+1}) $. This will prove that the factor $A^n/A^{n+1}$ is isomorphic to the vector space $(R/A)t^n$ and hence $gr(R)\cong (R/A)[t]$.

In fact if 
$(\lambda_1-\lambda_2) t^n\in A^{n+1}$ then there exists $y\in R$ such that 
$(\lambda_1-\lambda_2) t^n=t^{n+1}y$ which yields 
$t^n(ty-\lambda_1+\lambda_2)=0$.  This  contradicts the assumption that $t$ is regular and the proof is complete.\bigskip

{\bf 3.2.} We consider in this subsection an ideal $A$ of $R$ which is generated by a polycentral independent  system of elements $<t_1,t_2>$. We recall that this means that  $t_1$ is a central  element in $R$ such that $\bigcap_{i=1}^{\infty} (t_1)^i=0$ and the graded ring $gr(R)$ associated to the filtration defined by the powers of the ideal $(t_1)$ is isomorphic to the polynomial ring $R_1[t]$, where $R_1=R/(t_1)$;   the element  $t_2$ is central  modulo the ideal $(t)$,    the ideal $A_1=A/(t_1)$ generated in the  ring  $R_1$  by the element  $t_2$ is residually nilpotent and the graded ring associated to the filtration defined by the powers of $A_1$ is isomorphic to the polynomial ring $(R_1/A_1)[t_2]\cong (R/A)[t_2]$. 

Pick an arbitrary $m$  and let  $\bar{X}$ be 
the image of a subset $X$ under the natural homomorphism 
$R\longrightarrow R/(t_1)^m$. Clearly we have a natural homomorphism $\bar{R}\longrightarrow \bar{R}/(\bar{A})\cong R/A$ and the system of elements $\mathcal{X}$  can be considered also as a system of coset representatives for quotient ring $\bar{R}/\bar{A}\cong R/A$. We order the  system of elements $<t_1,t_2>$  assuming that $t_1<t_2$; since  there is one to one correspondence between the systems $<t_1,t_2>$ and $<\bar{t}_1,\bar{t}_1>$  the system $<\bar{t}_1,\bar{t}_2>$ is also well ordered. The standard monomials on $<t_1,t_2>$ or on $<\bar{t}_1,\bar{t}_2>$ are defined in the usual way. \bigskip

{\bf Lemma 3.2.} {\it Let 
\begin{equation}x=r_1\bar{t}_2r_2\bar{t}_2\cdots r_n\bar{t}_2r_{n+1}\end{equation} where $r_{\alpha}\in \bar{R}\,\, (\alpha=1,2,\cdots,n+1)$. Then there exists a representation 

\begin{equation} x=r \bar{t}_2^n+\bar{t_1} x_1\end{equation}

where $s\in \bar{R}$ and $x_1\in \bar{A}^{n-1}$ }.\bigskip

{\bf Proof.} Since the element  $\bar{t}_1$ is central modulo the ideal $(t_1)$  we obtain that for an arbitrary $\bar{r}\in \bar{R}$ there exists $\bar{s}\in \bar{R}$ such that 

 \begin{equation}\bar{t}_2\bar{r}=\bar{r}\bar{t}_2+\bar{t_1}\bar{s}\,\, \end{equation}

 We apply this identity to the factors  $\bar{t}_2$ and $\bar{r}_{n+1}$ in 
$(3.5)$ and obtain 

\begin{equation} x=\bar{r}_1\bar{t}_2\bar{r}_2\bar{t}_2\cdots\bar{r}_n
\bar{r}_{n+1}\bar{t}_2+\bar{t_1}\bar{y}\end{equation}

where $\bar{y}\in \bar{A}^{n-1}$. We repeat  this procedure, and after $n$ steps obtain $(3.6)$
where $r=\bar{r}_1\bar{r}_2\cdots \bar{r}_{n+1}$.\bigskip

{\bf Corollary 3.1.} {\it For an arbitrary natural $n$ }

 \begin{equation}\bar{A}^n\bigcap ( \bar{t_1})\subseteq \bar{t}_1\bar{A}^{n-1}\end{equation}

{\bf Proof.} Let $x\in \bar{A}^n\bigcap (\bar{t}_1)$.  Since $x\in \bar{t}_1$ we obtain  from  representation $(3.6)$ that $ \bar{r}\bar{t}_2^n\in (\bar{t}_1)$. We obtain from the last inclusion that $\bar{r}\in (\bar{t}_1)$ because the   element $t_2$ is regular in the quotient ring $\bar{R}/(\bar{t}_1)\cong R/(t_1)=R_1$.
We conclude  from this and Lemma $3.2.$ that $x\in \bar{t_1}\bar{A}^{n-1}$ and the proof is complete.\bigskip

{\bf Lemma  3.3.} \begin{equation} \bigcap_{n=1}^{\infty} \bar{t}_1^{m-1}\bar{A}^n=0\end{equation}

{\bf Proof.} We consider the regular representation $\rho$ of $\bar{R}$ in the ideal $(\bar{t}_1)^{m-1}$ generated by the central element $\bar{t}_1$. We prove  first of all that the kernel of this representation is the ideal $(\bar{t}_1)$. Clearly, the kernel contains the ideal $(\bar{t}_1)$.  On the other hand  assume that  $\bar{r}\not\in (\bar{t}_1)$,  $\bar{r}\bar{t}_1^{m-1}=0$ and let $r$ is an element of $R$ which is mapped in $\bar{r}$ under the homomorphism $R\longrightarrow R/(t_1)^m=\bar{R}$;  then $rt_1^{m-1}=t_1^ms$ for some element $s\in R$. Hence $t_1^{m-1}(r-t_1s)=0$ which is impossible because $t_1$ is a regular element.  This  proves that the the kernel of $\rho$ is the ideal $(\bar{t}_1)$.

We obtain now that the ideal $(\bar{t}_1)^{m-1}=\bar{t}_1^{m-1}\bar{R}$ is a one dimensional free module with generator $\bar{t}_1^{m-1}$  over the ring $\bar{R}/\bar{t}\cong R/(t_1)\cong R_1$ and the $R_1$-module $\bar{t}_1^{m-1}A_1^n$ is isomorphic to $A_1^n$.  We have now 
\begin{equation}\bar{t}_1^{m-1}\bar{R}=\bar{t}_1^{m-1}R_1\end{equation}

and we obtain from this that   $\bar{t}_1^{m-1}\bar{A}^n=\bar{t}_1^{m-1}A^n_1$
for every natural $n$
. Relation $(3.10)$ now follows from the relation $\bigcap_{n=1}^{\infty}A_1^n=0$. This completes the proof. \bigskip

{\bf Proposition 3.1.} {\it Let $\bar{X}$ denote the image of a subset $X\subseteq R$ under the natural homomorphism  the $R\longrightarrow R/(t)_1^m$ and assume that  $n\geq (m-1)$. Then} 
\begin{equation}\bar{A}^n\bigcap(\bar{t}_1)^{m-1}\subseteq (\bar{t}_1)^{m-1}\bar{A}^{n-(m-1)}\end{equation} 

{\bf Proof.} We apply induction by the number $m$; the initial step of the induction when $m=1$  is obvious and we can assume that the assertion has already been proven for all the quotient rings $R/(t_1)^k\,\, (k=1,2,\dots,m-1)$. This assumption implies that

 \begin{equation}\bar{A}^n\bigcap (\bar{t}_1)^{m-2}\subseteq \bar{t}_1^{m-2}\bar{A}^{n-(m-2)}\end{equation}

 and hence

 \begin{equation}\bar{A}^n\bigcap (\bar{t}_1^{m-1})\subseteq \bar{t}_1^{m-2}\bar{A}^{n-(m-2)}\end{equation}

 and it is enough to prove that if $n\geq (m-1)$ then 

\begin{equation}\bar{t}_1^{m-2}\bar{A}^{n-(m-2)}\bigcap(\bar{t}_1)^{m-1}\subseteq \bar{t}_1^{m-1}\bar{A}^{n-(m-1)}\end{equation}
if $n\geq (m-1)$.

Let $\bar{y}$ be an element from the left side of $(3.15)$. We obtain from Lemma  $3.2.$  that   
$\bar{y}$ is a sum of elements which either  have type
\begin{equation}\bar{t}_1^{m-2}r\bar{t}_2^k  \,\,    (k\geq n-m+2)\end{equation}

 or they have type

\begin{equation} \bar{t}_1^{m-2}\bar{t}_1\bar{y}_1=\bar{t}_1^{m-1}\bar{y}_1\end{equation}

with $ \bar{y}_1\in \bar{A}^{n-(m-2)-1}$.
 
Since   element $(3.17)$  belongs  to $\bar{t}_1^{m-1}\bar{A}^{n-(m-1)}$
  we can consider only the case when    $\bar{y}$ is a sum of elements of type 
$(3.16)$, that is 
\begin{equation}\bar{y}=\bar{t}_1^{m-2}\sum\bar{r}_{\alpha}\bar{t}_2^{k_{\alpha}}\end{equation}

 where $k_{\alpha}\geq n-m+1$.

  Since we assumed that $\bar{y}\in (\bar{t}_1)^{m-1}$ we conclude from $(3.18)$ that 

\begin{equation} (\bar{t}_1^{m-2}\sum \bar{r}_{\alpha}\bar{t}_2^{k_{\alpha}})\in (\bar{t})^{m-1}\end{equation}   

Let $u$ be an element of $R$ which is mapped in the element 
$\sum \bar{r}_{\alpha} \bar{t}_2^{k_{\alpha}} $ under the homomorphism 
$R\longrightarrow \bar{R}=R/(t)_1^m$. Relation $(3.19)$ implies that 
there exist  $a,b\in R$  such that $ut_1^{m-2}=at_1^{m-1}+bt_1^m $.   Since 
$t_1$ is a regular element we can cancel the last relation by $t_1^{m-2}$ and obtain that  
$u=t_1a+t_1^2b$.  Hence 
$\sum_{\alpha=1} (\bar{r}_{\alpha} \bar{t}_2^{k_{\alpha}})\in (\bar{t}_1)$ and the relation   $k_{\alpha}\geq (n-m+2)$  together with  Corollary  $3.1.$ imply that  

\begin{equation}(\sum_{\alpha=1} \bar{r}_{\alpha} \bar{t}_2^{k_{\alpha}})\in \bar{t}_1\bar{A}^{n-m+1}\end{equation}

 We obtain from this and  $(3.18)$ that $\bar{y}\in (\bar{t}_1)^{m-1}\bar{A}^{(n-m+1)}$  and the proof is complete.\bigskip

{\bf Corollary 3.2.} {\it Assume that the conditions of Proposition $3.1.$ hold. Then the ideal $\bar{A}$ of $\bar{R}$ is residually nilpotent}.\bigskip

{\bf Proof.} We recall also that the definition of the independent polycentral system  imply that the  ideal $A$ is residually nilpotent modulo $(t)_1$ and  we
can assume that $A$ is residually nilpotent modulo $(\bar{t}_1)^{m-1}$, i. e.

 \begin{equation}\bigcap_{n=1}^{\infty}A^n\subseteq (\bar{t}_1)^{m-1}\end{equation} 
We have now 

\begin{eqnarray}\bigcap_{n=1}^{\infty}\bar{A}^n=(\bigcap_{n=1}^{\infty} \bar{A}^n)\bigcap
(\bar{t}_1)^{m-1}\subseteq  (\bigcap_{n=1}^{\infty} \bar{A}^n\bigcap 
(\bar{t}_1)^{m-1})\subseteq\nonumber\\
\subseteq  \bigcap_{n=2(m-1)}^{\infty}(\bar{A}^n\bigcap (\bar{t}_1)^{m-1})\subseteq
 (\bigcap_{n=2(m-1)}^{\infty} (\bar{t}_1)^{m-1}\bar{A}^{(n-m+1)})
\end{eqnarray}
Lemma  $3.3.$  implies that the the last term in $(3.22)$ is zero. This completes the proof.\bigskip

{\bf Theorem 3.1.} {\it Let $R$ be a ring, $t_1,t_2>$ be a polycentral independent system in $R$. Let $A=<t_1, t_2>$ be the ideal generated by the elements 
$t_1, t_2$.  Then the ideal $A$ is residually nilpotent.}\bigskip

{\bf Proof.} We  pick an arbitrary natural
$m$ and consider the quotient ring $\bar{R}=R/(t_1)^m$; let $\bar{X}$ denote 
the image of a subset $X\subseteq R$ under the natural homomorphism $R\longrightarrow R/(t_1)^m$. Since the ideal $\bar{A}\subseteq
\bar{R}$ is residually nilpotent by Corollary $3.2.$  the assertion follows now from the condition $\bigcap_{m=1}^{\infty}(t_1)^m=0$. \bigskip

{\bf 3.3.} We assume throughout this subsection that the conditions of Theorem  
$3.1.$ hold and hence $\bigcap_{n=1}^{\infty}A^n=0$. We see that the powers of the ideal $A$ define a topology in $R$; let $\tilde{R}$ be the completion of $R$. Similarly we pick  an arbitrary  $m$ and consider the ring $R_m=R/(t_1)^m$; we denote by $\bar{X}$ the image of a subset $X\subseteq R$ under the homomorphism $R\longrightarrow R_m$. Corollary $3.2.$ yields   that the ideal $\bar{A}$ is residually
nilpotent in $R_m$. We denote by $\widetilde{R_m}$ the  completion of of $\bar{R}$ in the  topology defined by the powers of $\bar{A}$.

Let $\mathcal{X}$ be a system of coset representatives for the quotient ring 
$R_m/(\bar{t})\cong  R/(t)$.  and define  the standard  monomials on the set $\bar{T}$ in the usual way.  
\bigskip

{\bf Proposition 3.2.} {\it Every element  $r\in \widetilde{R_m}$ has a unique representation 

\begin{equation}x=\lambda_0+\lambda_1\bar{t}_1+\cdots+\lambda_{n-1}\bar{t}_1^{m-1}\end{equation}
where $\lambda_{\alpha}\,\, (\alpha=0,1,\cdots,m-1)$ is a power series 
\begin{equation}\sum_{i=1}^{\infty}u_it_2^i\end{equation}
with $u_i\in \mathcal{X} \,\, (i=1,2,\cdots)$.}\bigskip

{\bf Proof.}   We consider   the completion of the ideal $(\bar{t}_1)^{m-1}\subseteq R_m$. 
Proposition $3.1.$  implies that the topology induced in the subring $\bar{t}_1^{m-1}R_m$ is equivalent to the topology defined by the system of ideals 
$ (\bar{t}_1)^{m-1}\bar{A}^n$; hence the completion of  the ideal  
$(\bar{t}_1)^{m-1}$ is isomorphic to the ideal  $\bar{t}_1^{(m-1)}\widetilde{R_m}$ of 
$\widetilde{R_m}$ and hence   every element of this completion has a unique representation 

\begin{equation} u=\bar{t}_1^{m-1}\lambda\end{equation}

where $\lambda$ is  a power series of type $(3.24)$. Further the     ideal $\bar{t}_1^{(m-1)}\widetilde{R_m}$    is the kernel of the natural homomorphism $\widetilde{R_m}\longrightarrow \widetilde{R_{m-1}}$. We can assume that every element   $y\in \widetilde{R_{m-1}}$ has a unique representation 

\begin{equation} y=\mu_0+\mu_1\bar{t}_1+\mu_2\bar{t}_1^2+
\cdots\mu_{m-2}\bar{t}_1^{m-2}\end{equation}

where the coefficients $\mu_0,\mu_1,\cdots,\mu_{m-2}$ are power series of type $(3.24)$ and this representation is unique.

 The vector space $\tilde{R}_m$ is a direct sum $\tilde{R}_m=\tilde{R}_{m-1}+\bar{t}_1^{m-1}\tilde{R}_m$. We obtain from this  and from the representations 
$(3.25)$ and $(3.26)$ that $x$ has repesentation $(3.23)$; the uniqueness of this representation follows easily from the uniqueness of the representations of $(3.25)$ and $(3.26)$. This completes the proof.\bigskip

{\bf Proposition 3.3.} {\it Let $\tilde{R}$ be the completion of $R$ in the topology defined by the powers of the ideal $A$. Then every element of $r\in \tilde{R}$ has a unique representation 
\begin{equation} r=\sum_{j=1}^{\infty}\lambda_jt_1^j\end{equation}
where $\lambda_j=\sum_{i=1}^{\infty}u_{ji}t_2^i\,\, (j=1,2,\cdots)$
with $u_{ji}\in \mathcal{X}\,\, (i=1,2,\cdots;j=1,2,\cdots)$}.\bigskip

{\bf Proof.} The ring $\tilde{R}$ is an inverse limit of the system of rings $R/(t)^m\,\, (m=1,2,\cdots)$ and the existence and the uniqueness of representation $(3.27)$ follows from Proposition $3.2.$\bigskip

{\bf Corollary 3.3.} {\it Every element $r\in \tilde{R}$ has a unique representation 
\begin{equation} r=\sum_{i=1}^{\infty}u_i\tau_i\end{equation}
where $u_i\in \mathcal{X}\,\, (i=1,2,\cdots)$, $\tau_i \,\, (i=1,2,\cdots)$ are standard monomials on the set $<t_1,t_2>$ with $\lim_{i\to\infty}l(\tau_i)=\infty$}.\bigskip

 {\bf 3.4.} We set up now the notation for Theorem III and Lemmas $3.4.-3.10.$  Let $R$ be a ring,  $T=<t_1,t_2,\cdots, t_n>$ be an independent polycentral system which is composed from the central independent systems $T_1,T_2,\cdots,T_k$, $A$ be the ideal generated by the system $T$,  $f$ be a function on $T$ whose values form a bounded set of natural numbers and $f(t_1)>2f(t_2)$ for an arbitrary pair of elements $t_1\in T_i, t_2\in T_{i+1}\,\, (1\leq i\leq k-1 )$.  The definition of the polycentral independent system implies that  for an arbitrary $(1\leq m\leq k)$ the subsystem $t_1,t_2,\cdots, t_m$ is an independent polycentral system in $R$; let
  $A_m$ be the ideal generated $<t_1,t_2,\cdots, t_m>$.  Further,  
the system $t_{m+1},t_{m+2},\cdots, t_n$ is an independent polycentral system in the quotient ring $R/A_m$.

 Every element  $r\in A$ is sum  
\begin{equation}\sum \mu_1t_{\beta_1}\mu_2t_{\beta_2}\cdots\mu_k t_{\beta_k}\mu_{k+1}\end{equation}

where $t_{\beta_1},t_{\beta_2 }\cdots,t_{\beta_k} $ are elements from $T$, $\mu_1,\mu_2,\cdots,\mu_k \in R$. Let $M$ be the maximum of the values of $f(t)$ on $T$.  We  define now    a system of subsets  
$B_j\,\, (j=0,1,\cdots)$ in $R$ as follows.  Define  $B_0=R$ and then for $j\geq 1$ define that an   element $x$ belongs to $B_j$ if there exists for it a representation  $(3.29)$  such that every summand satisfies the condition 

\begin{equation} \sum_{\alpha=1}^k f(t_{\beta_{\alpha}})\geq j\end{equation}
 
It is clear that $B_1=A$ and that $0\in B_j\,\, (j=0,1,\cdots)$. Let $\mathcal{X}$ be a system of coset representatives for the ideal $A$.\bigskip
 
{\bf Theorem III.} i) {\it The ideal $A$ is residually nilpotent and $\bigcap_{j=0}^{\infty}B_j=0$. The system of subsetes $B_j\,\, (j=0,1,\cdots)$ is a filtration in $R$.}

ii) {\it  This   filtration  defines a pseudovaluation $v$ such that $v(t)=f(t)$ for an arbitrary $t\in T$;  the homogeneous components $\tilde{t}_i\,\, (i=1,2,\cdots,n)$ are central in $gr_v(R)$  and the ring $gr_v(R)$ is isomorphic to the polynomial ring $(R/A)[\tilde{t}_1,\tilde{t}_2,\cdots,\tilde{t}_n]$.

The pseudovaluation $v$ is defined uniquely by the conditions $v(t)=f(t)\,\, (t\in T)$,  and $gr_v(R)\cong (R/A)[\tilde{t}_1,\tilde{t}_2,\cdots,\tilde{t}_n]$; the topology defined by this pseudovaluation is equivalent to the topology defined by the powers of $A$. }\bigskip

{\bf 3.5.} We will prove in this subsection  a few auxilary facts.

 {\it We will  use throughout this subsection the notation of subsection $3.4.$  and we will assume that  the ideal $A$ generated by the system $T$ is residually nilpotent.  } We will prove under this  condition Lemmas $3.4.-3.10.$      We  point out that  Theorem $3.1.$ shows that the ideal $A$ is residually nilpotent if $T$ is the polycentral independent system $<t_1,t_2>$ . This fact  will make possible  to apply  Lemmas $3.4.-3.10.$  in the proofs  of Theorems $3.2.$ and $3.3.$\bigskip

{\bf Lemma 3.4.} {\it Every $B_j$  is an ideal in $R$}, 

\begin{equation}R=B_0\supseteq B_1\supseteq B_1\supseteq \cdots\end{equation}  and 

\begin{equation}A^j\subseteq B_j\,\, (j=1,2,\cdots)\end{equation}

{\it Further}  

 \begin{equation}B_{jM}\subseteq A^j\,\, (j=1,2,\cdots)\end{equation}

{\bf Proof.} The proof of the first statement is straightforward  and   relation $(3.29)$ follows immediately  as well as   $(3.32)$     

  Now  assume that $r\in B_{jM}$. Then there exists for $r$  representation  where all the summands  $(3.31)$  satisfy relation 

\begin{equation} \sum_{\alpha=1}^k f(t_{\beta_{\alpha}})\geq jM\end{equation}

 We pick an arbitrary of these summands.  Since  $f(t_{\beta_{\alpha}})\leq M\,\, (\alpha =1,2,\cdots,k)$ we must have $k\geq j$ and hence this summand belongs to $A^j$.   This proves relation $(3.33)$. The proof of the lemma is complete.\bigskip

{\bf Corollary 3.4.} {\it The system of ideals $B_j\,\, (j=0,1,\cdots)$ defines a pseudovaluation $v$ in $R$ as follows 
\begin{equation}v(0)=\infty; \, \mbox{if}\,\,  x\neq 0\,\, \mbox{then}\,\,  v(x)=max\{j|x\in B_j\}\end{equation}

 The topology defined by this pseudovaluation is equivalent to the topology defined by the powers of the ideal $A$.  

An element $r\in R$ has $v$-value greater than zero iff it belongs to the ideal $A$.

If $t\in T$ then $v(t)\geq f(t)$.}\bigskip

{\bf Proof.}   We obtain from $(3.32)$  and the assumption  $\bigcap_{j=1}^{\infty} A^j=0$   that 

\begin{equation}\bigcap_{j=1}^{\infty}B_j=0\end{equation}

Further the definition of the system of ideals $B_j\,\, (j=0,1,\cdots)$ shows that $B_{j_1}B{j_2}\subseteq B_{j_1+j_2}$. We see that the system of ideals $B_j\,\, (j=0,1,\cdots)$ forms a filtration in $R$ so $v$ is a pseudovaluation in $R$. 

 The definition of the pseudovaluation $v$ shows that  $v(x)>0$ iff  $x$ has a  representation $(3.29)$  with   every summand 
containing   elements of $T$ which means that $x\in A$.

The definition of $B_j$ implies that if $f(t)=j$ then $t\in B_j$, so $v(t)\geq j$ which proves the last statement.  \bigskip

 {\bf Definition 3.1.} {\it The pseudovaluation $v$ and the filtration $B_j\,\, (j=0,1,\cdots)$ will be called the  pseudovaluation and the filtration defined in $R$ by the polycentral system $T$ and the function $f$.}\bigskip

{\bf Lemma 3.5.} {\it The graded ring associated to the pseudovaluation $v$ is generated by the zero component $R/A$ and the set of homogeneous components 
$\tilde{T}=\{\tilde{t}_i\,\,|(t_i\in T,i\in I\}$.}\bigskip

{\bf Proof.} Let $r$ be an element of $R$ with $v(r)=j$. Then $r\in B_j$ and we will  consider the image of this element  in the factor $B_j/B_{j+1}$. Since we consider the image of $r$ modulo $B_{j+1}$ we  can assume that  it has  representation   $(3.29)$   where every summand $r_1=\mu_1t_{\beta_1}\mu_2 t_{\beta_2}\cdots\mu_k t_{\beta_k}\mu_{k+1}$ 
satisfies condition

\begin{equation} \sum_{\alpha=1}^k f(t_{\beta_{\alpha}})=j\end{equation}

and noone of these summands  belongs to $B_{j+1}$. We conclude  that  $v(r_1)=j$ otherwise we would  have $r_1\in B_
{j+1}$.   The relation   $v(r_1)=j$ implies  that  $v(\mu_{\alpha})=0\,\, (\alpha=1,2,\cdots,k+1)$. In fact, if we assume   $v(\mu_{\alpha})>0)$ for some $\mu_{\alpha}$ then we obtain from this assumption together with condition $(3.37)$ that $v(r_1)>j$.  We conclude from this and from Lemma $2.3.$ that the homogeneous component $\tilde{r}_1$ of the summand  $r_1$  has representation 

\begin{equation}\tilde{r_1}=\tilde{\mu}_1 \tilde{t}_{\beta_1}\tilde{\mu}_2\tilde{t}_{\beta_2}\cdots\tilde{\mu}_k\tilde{t}_{\beta_k}\tilde{\mu}_{k+1}\end{equation}

and that

 \begin{equation} \tilde{r}=\sum\tilde{\mu}_1 \tilde{t}_{\beta_1}\tilde{\mu}_2\tilde{t}_{\beta_2}\cdots\tilde{\mu}_k\tilde{t}_{\beta_k}\tilde{\mu}_{k+1}\end{equation}

and  the assertion  follows.\bigskip

The following fact was established in the proof of Lemma $3.5.$\bigskip

{\bf Corollary 3.5.} {\it Let $r$ be an element of $R$ such that $v(r)=j$. Then the homogeneous component $\tilde{r}$ is a sum of monomials $(3.38)$ where $v(\tilde{u}_{\alpha})=0\,\, (\alpha=1,2,\cdots,k+1)$ and 
$\sum_{\alpha=1}^k v(\tilde{t}_{\beta_{\alpha}})=j$.}\bigskip

{\bf Lemma 3.6.} {\it Assume that  the homogeneous components $\tilde{t}_i\,\, (t_i\in T)$ are central in $ gr(R)$. Let $x\in B_j\backslash B_{j+1}$. Then }

i) {\it there exists a representation }
\begin{equation}x=x_1+y\end{equation}
{\it  where }

\begin{equation}x_1=\sum_{i=1}^n\lambda_i\pi_i\end{equation}

 $\lambda_i\in \mathcal{X}\,\, (i=1,2,\cdots,n)$;
{\it  $\pi_i\,\, (i=1,2,\cdots,n)$ are standard monomials on $T$ with $v$-value  equal  $j$ and $y\in B_{j+1}$.}  

ii) {\it Let  $\tilde{R_v}$  be the completion of $R$ in the topology defined by $v$. Then the element $y$ in the right side of $(3.40)$  has a 
representation }

\begin{equation} y=\sum_{i=n+1}^{\infty} \lambda_i\pi_i\end{equation}

{\it where $\lambda_i\in \mathcal{X}$,  $v(\pi_i)\geq j+1\,\, (i=1,2,\cdots)$,
  and } $lim_{i\to\infty}v(\pi_i)=\infty$. 

iii) {\it Representations $(3.40)-(3.42)$  yield a power series representation} 

\begin{equation} x=\sum_{i=1}^{\infty} \lambda_i\pi_i\end{equation}

{\bf Proof.} i)  We have for $x$  a representation $(3.29)$. We consider now the homogeneous component $(3.38)$ of one of the summands in this representation. Since the homogeneous component $\tilde{t}_i\,\, (i\in I)$ are central in $gr(R)$ we obtain that the element $(3.38)$ is equal to the element 
$\tilde{\mu}_1\tilde{\mu}_2\cdots\tilde{\mu_{k+1}}\pi(\tilde{t}_{\beta_1},\tilde{t}_{\beta_2},\cdots, \tilde{t}_{\beta_k})$    where $v(\tilde{\mu}_1\tilde{\mu}_2\cdots\tilde{\mu_{k+1}})=0$  and  $\pi$ is a suitable  standard monomial of value $j$ on the set of elemenrs $\tilde{t}_{\beta_1},\tilde{t}_{\beta_2},\cdots,\tilde{t}_{\beta_k}$. 
 
We obtain from this

\begin{equation}\mu_1t_{\beta_1} \mu_2t_{\beta_2}\cdots \mu_kt_{\beta_k}\mu_{k+1}=\mu_1\mu_2\cdots\mu_k\mu_{k+1}\pi+x_0\end{equation}

 where $\pi $ is a standard monomial with value $j$ on the of elements 
$t_{\beta_1},t_{\beta_2},\cdots, t_{\beta_k} $ obtained from $\pi(\tilde{t}_{\beta_1},\tilde{t}_{\beta_2},\cdots, \tilde{t}_{\beta_k})$ by substitution $\tilde{t}_{\beta_{\alpha}}\longrightarrow t_{\beta_{\alpha}}\,\, (\alpha=1,2,\cdots,k)$,   $x_0\in B_{j+1}$ and $\mu_1\mu_2\cdots\mu_k\mu_{k+1}=r$ is an element of $R$.  We obtain from this that $x$ has a representation 
\begin{equation}x=\sum_{i=1}^m r_i\pi_i+r_0\end{equation}

 where  $\pi_i\,\, (i=1,2,\cdots,m)$  are standard monomials  on $T$ which belong to $B_j$,   
$r_i\in R\,\, (i=1,2,\cdots,m)$, $r_0\in B_{j+1}$. By adding, if necessary, the coefficients of the same  monomial, we can assume that   that all the monomials $\pi_i\,\, (i=1,2,\cdots,m)$ are lexicographically distinct; further if a coefficient $r_i$ belongs to $A$ then the product $r_i\pi_i$ belongs to $B_{j+1}$ so we can assume that $r_i\in R\backslash A\,\, $ for $i=1,2,\cdots,  m$. Hence, 
  $r_i=\lambda_i+u_i,0\not= \lambda_i\in \mathcal{X}, u_i\in A=B_1\subseteq B_1\,\, (i=1,2,\cdots,m)$. We substitute these expressions in $(3.45)$ and obtain

\begin{equation} x=\sum_{i=1}^n \lambda_i\pi_i+\sum_{i=1}^mu_i\pi_i+r_0\end{equation}

The summands $u_i\pi_i\,\, (i=1,2,\cdots,m)$ belong to $B_{j+1}$, we denote now $y=\sum_{i=1}^mu_i\pi_i+r_0$ and  obtain relations $(3.40)$- $(3.42)$.

 This proves statement i).  \bigskip

ii)  Let $x\in R$ with $v(x)=j$. We obtain from statement i) representation $(3.40)$.  It is important that we can assume that all the  coefficients $\lambda_i$ in the element $x_1=\sum_{i=1}^n\lambda_i\pi_i $ are non-zero otherwise we would get $x\in B_{j+1}$. 
We  consider now the element $y$ which has $v$-value $j_1>j$. Once again, statement i)   yields
       that $y=x_2+z$, where $x_2$ is a linear combination with coefficients from $\mathcal{X}$ of  standard monomials  with $v$-value $j_1$ and $z\in B_{j_1+1}$. 
  
We obtain by this argument that there exists  a representation
 \begin{equation}x=x_1+x_2+x_3+\cdots\end{equation} 

where $x_1$ is a linear combination of monomials with value $j$, and $x_2,x_3,\cdots$ are linear combinations of monomials with values  $j_1,j_2,\cdots$ greater than $j$.  This series converge in the $\rho$-topology and in the topology defined by the pseudovaluation $v$. This completes the proof of statement ii) .

Statement iii) follows from i) and ii).\bigskip 

{\bf Corollary 3.6.} {\it Assume that the conditions of Lemmas $3.6.$ and $3.7.$ hold.  Then the  homogeneous component $B_j/B_{j+1}$ is a left module  over ring $R/A$, the standard monomials $\pi_i$ with value $j$  on the set $T $  form a system of generators for it.}\bigskip

{\bf Proof.} Since  $AB_j\subseteq B_{j+1}$  we obtain that $B_j/B_{j+1}$ is a left module over $R/B_1\cong R/A$ and the assertion now follows from Lemma 
$3.6.$ \bigskip

   {\bf Lemma 3.7.} {\it Assume that conditions of Lemmas $3.6.$  hold. Then the  following  statements are equivalent:

$1)$  Representation $(3.42)$ is unique for every element $x\in B_j\,\, (j=0,1,\cdots)$.

$2)$ The associated graded ring $gr_v(R)$ is isomorphic to the polynomial ring $(R/A)[\tilde{t}_1,\tilde{t}_2,\cdots, \tilde{t}_n]$.}\bigskip

{\bf Proof.}  We prove first   that  $1)\longrightarrow 2)$. Assume that representation $(3.42)$ is unique.  We recall that  the zero component of $gr(R)$ is $R/A$ and Corollary $3.6.$ implies   that $B_j/B_{j+1}$ is a module  over the quotient ring $R/A$.   

We show now  that standard monomials with value $j$ must be linearly independent in $B_j/B_{j+1}$ over $R/A$.  In fact if this is not true then there exists an element $y\in B_{j+1}$ such that 

\begin{equation} y=\sum_{i=1}^m\lambda_i\pi_i\end{equation} 

where 

\begin{equation}v(\pi_i)=j,   0\not=\lambda_i\in \mathcal{X}\,\, (i=1,2,\cdots,n)\end{equation}

On the other hand, since $y\in B_{j+1}$ we obtain from Lemma $3.6.$  a power series representation for $y$ where the minimal value of the monomials is greater than or equal $j+1$. We obtained  two representations for the element $y$. This contradiction shows that the standard monomials of weight $j$ must be linearly independent. Further, if $\tilde{B}_j$ is the completion of $B_j$ in the $v$-topology then $\tilde{B}_j/\tilde{B}_{j+1}\cong B_{j}/B_{j+1}$. But $\tilde{B}_j/\tilde{B}_{j+1}$ is generated over $R/A$ by the standard monomials with value $j$, we obtain from  that the standard monomials on $T$ with value $j$ form a basis of 
$B_j/B_{j+1}$,  the homogeneous components of the elements of $T$  are central in $gr(R)$  so we conclude that the ring $gr(R)$ is isomorphic to $(R/A)[\tilde{t}_1,\tilde{t}_2,\cdots,\tilde{t}_n]$.  This  proves that $1\longrightarrow 2)$.

 Now prove that $2) \longrightarrow 1)$. Assume that $gr(R)\cong (R/A)[T]$. Let $x\in B_j\backslash B_{j+1}$. We have for $x$ representation $(3.40)$ where $x_1$ is given by $(3.41)$ with  the coefficients $\lambda_1, \lambda_2,\cdots,\lambda_n$ uniquely defined. We repeat this procedure with the element $y$ as in the proof of statement ii) of Lemma $3.6.$ and obtain a power series representation for $x$ with uniquely defined coefficients $\lambda_i\in \mathcal{X}\,\, (i=0,1,\cdots)$. This completes the proof. \bigskip

{\bf 3.6.} We use throughout this section the same notation as in subsections $3.4.-3.5.$ In particular, $T=<T_1,T_2,\cdots, T_k>$ is an independent polycentral system, $A$ is the ideal generated by this system and we assume that this ideal is residually nilpotent.   We consider now the central independent  subsystem  $T_1=<t_1,t_2,\cdots, t_{n_1}>$ and the ideal $(t_1)$ of $R$. 
Let $\bar{X}$ denote the image of a subseteq $X\in R$ under the natural homomorphism $R\longrightarrow \bar{R}=R/(t_1)$.  If $n_1>1$ then the definition of the independent polycentral system implies that the system $\bar{T}=<t_2,t_3,\cdots,t_{n_1-1};T_2,T_3,\cdots, T_k>$ is an independent polycentral system in $\bar{R}$; if $n_1=1$ we obtain in $\bar{R}$ an independent polycentral system $\bar{T}=<T_2,T_3,\cdots,T_k>$. In both cases the system $\bar{T}\subseteq \bar{R}$ has smaller length then the system $T$,  we keep for the restriction of the function $f$  on $\bar{T}$ the same notation $f$.   We will assume in the following Lemmas $3.8.$ and $3.9$ that  the conclusions of Theorem III hold in the ring $\bar{R}=R/(t_1)$ for the polycentral independent system $\bar{T}$ and  the function $f$. This means that:  

$1)$ the ideal $\bar{A}$ is residually nilpotent; 

$2)$ there exists a pseudovaluation $v_1$ such that $v_1(\bar{t})=f(\bar{t})$;  

$3)$ the  homogeneous components of the elements $\bar{t}\in \bar{T}$ are central in the ring $gr_{v_1} \bar{R}$; 

$4)$ these homogeneous components are algebraically independent in   the graded ring $gr_{v_1}(\bar{R})$ over the zero degree homogeneous component   $(\bar{R}/\bar{A})\cong (R/A)$ and that they  generate the ring $gr_{v_1}(\bar{R})$ over $\bar{R}/\bar{A}\cong R/A$.

Lemmas $3.8.$ and $3.9.$ will be used in the induction proof of Theorem III in subsection $3.9.$ The initial step of the induction is the case when the system $T$ consists from one central  element $t$ such that $\bigcap_{i=1}^{\infty}(t)^i=0$ and the graded ring $gr(R)$ asociated to the pseudovaluation $\rho$ defined by the powers of $t$  is isomorphic to the polynomial ring in the variable $t$ over $R/(t)$. The weight function $f(t)=M$ defines a pseudovaluation $v$ which is equivalent to $\rho$ and we see that all the conclusions of Theorem III hold hold for the initial step of the induction.

 We consider now the homomorphism  $R\longrightarrow \bar{R}$.  Proposition $2.3.$ implies that this homomorphism defines the  filtration $\bar{B}_j\,\, (j=0,1,\cdots)$ and a pseudovaluation $\bar{v}$ defined by this filtration. We need the following fact.\bigskip

{\bf Lemma 3.8.}  {\it Assume that the conclusions of Theorem III hold in the ring $\bar{R}$. Then the  pseudovaluation $\bar{v}$ coincides with $v_1$}:

\begin{equation}\bar{v}(\bar{r})=v_1(\bar{r})\end{equation}

{\it for  every $\bar{r}\in \bar{R}$.} \bigskip

{\bf Proof.}  We consider the filtration  $B^{\ast}_j\,\, (j=0,1,\cdots)$  defined by the pseudovaluation $v_1$ and will  prove that $\bar{B}_j=B^{\ast}_j\,\, (j=0,1,\cdots)$. 

Assume that $0\not=\bar{x} \in \bar{B}$ and let $x\in R$ be an element whose image in $\bar{R}$ is $\bar{x}$. Then there exists for $x$ a representation as a sum of monomials $(3.29)$  which satisfy condition $(3.30)$.
If a summand $(3.29)$ in the representation of $x $ contains the factors from $t_1$ its image in $\bar{R}$ is zero. Otherwise its image is the element 

 \begin{equation}\bar{\mu}_1t_{\beta_1}\bar{\mu}_2t_{\beta_2}\cdots\bar{\mu}_kt_{\beta_k}\bar{\mu}_{k+1}\end{equation}

 We obtain from this  that $\bar{x}\in B^{\ast}_j$ which implies that $\bar{B}_j\subseteq B^{\ast}_j\,\, (j=0,1,\cdots)$. The proof of the reverse  inclusion is immediate, and the asertion follows. \bigskip

{\bf Lemma 3.9.}    {\it    Assume that the  conclusions of Theorem III hold in the quotient ring $\bar{R}=R/(t_1)$.
Then $v(t)=f(t)\,\, \mbox{if}\,\, t\not=t_1$ and all the  homogeneous components $\tilde{t}\,\, (t\in T)$ are central in the graded ring $gr_v(R)$}. \bigskip

 {\bf Proof.} The definition of $v$ shows that $v(t)\geq f(t)$ if  $t\in T$. Since $\bar{v}(\bar{r})\geq v(r)$ for an arbitrary $r\in R$ and 
$v_1(t))=\bar{v}(t)=f(t)$ for $t\not=t_1$ we obtain the first statement of the lemma.

We prove now the second statement. Lemma $3.5.$ implies that the graded ring $gr_v(R)$ is generated by 
the zero homogeneous component $R/A$ and the homogeneous components of the elements $t\in T$ so we have to prove only that

\begin{equation}v([t_{i_1},t_{i_2}])>v (t_{i_1})+v(t_{i_2})\end{equation}

 and for every $t\in T$ and $r\not \in A$  

\begin{equation} v([r,t])> v(r)+v(t)\end{equation}

 If any  of the elements $t_{i_1}$ or $t_{i_2}$ belongs   to $T_1$  then it  is central and relations $(3.52)$ and $(3.53)$ are   obvious, so  we can assume that $t_{i_1},t_{i_2}\in T_2\bigcap T_3\bigcap\cdots\bigcap T_k$.  We have already proven that   $v(t_{i_1})=f(t_{i_1}), v(t_{i_2})=f(t_{i_2})$ and  the condition of the assertion imply that  the homogeneous components of $t_{i_1}$ and $t_{i_2}$ commute in  the graded ring $gr_{v_1}(\bar{R})$. We have therefore in $\bar{R}$

 \begin{equation} [t_{i_1},t_{i_2}]=\bar{a}\end{equation}

where 

\begin{equation}v_1(\bar{a})=\bar{v}(\bar{a})>
(\bar{v}(t_{i_1})+\bar{v}(t_{i_2}))\end{equation}

We find now an element $a\in R$ whose image in $\bar{R}$ is $\bar{a}$ and $v(a)=\bar{v}(\bar{a})$ and obtain from $(3.54)$

\begin{equation} [t_{i_1},t_{i_2}]=a+b\end{equation}

 where $b\in (t_1)$ and $a$ is an element  with value 
$v(a)>f(t_{i_1})+f(t_{i_2}) $. Since $b$ is a multiple of $t_1$  is a multiple of $t_1$  the $v$-value of $b$  is  greater than 
$  f(t_{i_1})+f(t_{i_2})$ because   $t_1\in T_1$,   $t_{i_1},t_{i_2}\in (T_2\bigcap T_3\bigcap\cdots\bigcap T_k)$ and $f(t_1)>2f(t)\,\, \mbox{for}\,\, t\not\in T_1$.   We conclude therefore that $v(a+b)>f(t_{i_1})+f(t_{i_2}) $, and this proves $(3.52)$

We will now  prove $(3. 53)$.   The image $\bar{r}$ of $r$  does not belong to $A$, so $v_1(\bar{r})=\bar{v}(\bar{r})=0$ and 
we have now  in $\bar{R}$

$$\bar{v}([\bar{r},t])>\bar{v}(\bar{r})+\bar{v}(t)=\bar{v}(t)$$
and the proof of $(3.53)$ can be completed in the same way as of $(3.52)$.\bigskip

 We can now   prove Theorem $3.2.$   which is an important step in the proof of Theorems I-III. \bigskip

{\bf Theorem 3.2.} {\it Let $R$ be a ring, $<t_1,t_2> $ be an   independent polycentral system, $A$ be the ideal generated by these elements, $M_1,M_2$ be two natural number such that $M_1>2M_2$.  Then there exists a pseudovaluation $v$ such that $v(t_i)=M_i\,\, (i=1,2,\cdots)$,   
the  graded ring associated to the pseudovaluation $v$  is isomorphic to the polynomial ring $(R/A)[\tilde{t}_1,\tilde{t}_2]$,  the topology defined by this valuation is equivalent to the $\rho$-topology defined by the powers of the ideal $A$. The pseudovaluation $v$ is defined uniquely by these properties.}\bigskip

{\bf Proof.}  Theorem  $3.1.$   states that the ideal  $A$ formed by the system $T=<t_1,t_2>$  is residually nilpotent. We  consider the filtration $B_j\,\, (j=0,1,\cdots)$ and the pseudovaluation $v$ defined by the system $T$ and the function $f(t_i)=M_i\,\, (i=1,2)$ (see Definition $3.1.$)
 Since the element $t_1$ is central in $R$ its homogeneous component is central in $gr(R)$.  Lemma $3.9.$  implies  that the homogeneous component $\tilde{t}_2$ is  central in $gr(R)$. Since the topology defined by the pseudovaluation $v$ is equivalent to the $\rho$-topology defined by the powers of the ideal $A$ we obtain from Corollary $3.3.$ that every element of the completion $\tilde{R}$ of $R$ has a unique representation $(3.43)$ .
We obtain now from Lemma $3.7.$ that $gr_v(R)\cong (R/A)[\tilde{t}_1,\tilde{t}_2]$.

This completes the proof.\bigskip

{\bf 3.7.  }  We will prove in this subsection  the following Theorem $3.3.$ which is the main step in the proof of Theorems I-III.\bigskip

{\bf  Theorem 3.3.} {\it  Let $R$ be a ring, $T=<t_1,t_2,\cdots,t_n>$  be an independent  polycentral system in $R$,  $k\geq n$  be a  natural
 numbers.  Then the weight function    defined by the function $f_k$ on $T$ as  

\begin{equation} f_k(t_i)=3^{k-i+1}  (i=1,2,\cdots,n)\end{equation}
 
 extends to $t$-adic pseudovaluation $v_0$ of $R$ with associated graded ring  isomorphic to the  polynomial ring $(R/A)[\tilde{t}_1,\tilde{t_2},\cdots,\tilde{t}_n]$.  The topologies defined by the pseudovaluation $v_0$  and by the powers of the ideal $A$  are equivalent}.\bigskip

 Let $R$ be a ring, $t_i\,\, (i\in I)$  be  a  system of elements in $R$, $A$ be the ideal generated by them. Let 
 $v$ be a pseudovaluation in $R$ such that 
$v(t_i)=mm_i\,\, (i\in I)$  where $m,m_i\,\, (i\in I)$ is a system of natural numbers and let $\tilde{t}_i$ be the homogeneous component of $t_i$ in the associated graded ring $gr(R)$.  Assume that the   associated   graded ring $gr(R)$ is  isomorphic to the polynomial ring over $(R/A)$ in the system of variables $\tilde{t}_i\,\, (i\in I)$.    Let $R[t,t^{-1}]$ be the Laurent polynomial ring over $R$, we extend the pseudovaluation $v$ to this ring by defining $v(t)=m$.  Let  $V=\{x\in R[t,t^{-1}]| v(x)\geq 0\}$ be the ring of integers of $R[t,t^{-1}]$ and $\bar{X}$ be  the image of a subseteq $X\subseteq V(R[t,t^{-1}])$ under the homomorphism $V\longrightarrow V/(t)$.\bigskip

{\bf Lemma 3.10.} {\it\it The  quotient ring $\bar{V}=V/(t)$ is isomorphic to the  polynomial ring over $R/A$ in the system of variables  $u_i=\overline{t_it^{-m_i}}\,\, (i\in I)$ }.\bigskip

{\bf Proof.} The pseudovaluation $v$ is equivalent to the pseudovaluation $v_1$ defined as $v_1(t_i)=m_i\,\, (i\in I); v_1(t)=m$ and the assertion now follows from Proposition $2.8.$\bigskip

We have the immediate corollary of Lemma  $3.10.$\bigskip

{\bf Corollary 3.8.} {\it The map $\tilde{t}_i\longrightarrow 
\overline{tt_i^{-m_i}}$ extends to an isomorphism between the rings $gr(R)$ and $V/(t)$}.\bigskip

{\bf Proof of Theorem 3.3.} We  recall that the definition of the independent polycentral system $<t_1,t_2,\cdots, t_n>$ implies that  for every $m<n$ the subsystem $<t_1,t_2,\cdots,t_m>$ is an independent polycentral system in $R$ and the   system $<t_{m+1},t_{m+2},\cdots, t_{n}>$  is a polycentral 
independent  system in the quotient ring $R_m=R/A_{m}$ where $A_m$ is the ideal generated by $t_1,t_2,\cdots,t_m$. We make an assumption of an induction that the assertion is proven for an arbitrary ring $S$ with a polycentral independent system $x_1,x_2,\cdots, x_{n_1}$ of length $n_1<n$ and a weight function $f_{k_1}\,\, (k_1\geq n_1)$ defined on this system by the function 

\begin{equation} f_{k_1}(x_1)=3^{k_1}, f_{k_1}(x_2)=3^{k_1-1},\cdots, f_{k_1}(x_{n_1-1})=3^{k_1-n_1+1}\end{equation} 
We have already pointed out in subsection $3.6.$  that the initial step of the induction when $T$ contains only one element is true. Further,    assumption of the  induction implies in particular that for every $k\geq n-1$ there exists a pseudovaluation 
 $v_{k-1}$ of $R$ such that  $v_{k-1}(t_i)=3^{k-i+1}\,\, ( i=1,2,\cdots,n-1) $ 
and the associated graded ring $gr_{v_{k-1}}(R_1)$  is isomorphic to the polynomial ring $R_{n-1}[\tilde{t}_1,\tilde{t}_2,\cdots,\tilde{t}_{n-1}]$; the definition of the independent polycentral system implies that the element  $t_n$ is central  in $R_{n-1}$, the ideal $(t_n)R_{n-1}$ is residually nilpotent and the graded ring associated to the filtration $(t_n)^i R_n\,\, (i=1,2,\cdots)$ is isomorphic to the polynomial ring  generated over $R_{n-1}/(t_nR_{n-1})\cong R/A$ by the homogeneous component of $t_n$.  

We  consider now the  group ring  
$R[t,t^{-1}]$ and extend $v_{k-1}$ to this group ring defining $v_{k-1}(t)=3^{k-n+2}$; let $V=\{x\in R[t,t^{-1}] | v_{k-1}(x)\geq 0\}$ be the subring of $v_{k-1}$-integers of  $V$. We obtain from Lemma  $3.10.$   that  the quotient ring $V/(t)$ is isomorphic to the polynomial ring $R_{n-1}[\tilde{u}_1,\tilde{u}_2,\cdots,\tilde{u}_{n-1}]$  where $u_1=t_1t^{-3^{n-2}}, u_2=t_2t^{-3^{n-3}},\cdots, u_{n-1}=t_{n-1}t^{-1}$ and $v_{k-1}(u_i)=0\,\, (i=1,2,\cdots,n-1)$. We have  in $R_{n-1}$ a central  element  $t_n$ which is also  central  in the polynomial ring $R_{n-1}[\tilde{u}_1,\tilde{u}_2,\cdots, \tilde{u}_{n-1}]$ and we obtain now in $V$ an independent  polycentral system $<t, t_n>$. We obtain now  from Theorem    $3.2.$  that there exists   a pseudovaluation $v_0$ in $V$ such that 
$v_0(t)=3^{k-n+2}, v_0(t_n)=3^{k-n+1}$ with the associated graded ring a polynomial ring in the system of variables $<\tilde{t},\tilde{t}_n>$ over  the ring 

$$(R_{n-1})[\tilde{u}_1,\tilde{u}_2,\cdots,\tilde{u}_{n-1}])/(t_n)\cong (R/A)[\tilde{u}_1,\tilde{u}_2,\cdots,\tilde{u}_{n-1}]$$

The restriction of this pseudovaluation to $R$ defines a pseudovaluation in 
$R$;  we will use for it the same notation $v_0$.  We will now verify that this pseudovaluation satisfies all the conclusions of the theorem.

We have 
\begin{equation} t_1=u_1t^{3^{n-2}}, t_2=u_2 t^{3^{n-3}},  t_{n-1}=u_{n-1}t\end{equation}  
which yields

 \begin{eqnarray} v_0(t_i)= v_0(u_i 3^{n-i-1}t)=v_0(u_i)+3^{n-i-1}v_0(t)= \nonumber\\
=3^{n-i-1}v_0(t)=3^{k-i+1}\,\, (1\leq i\leq n-1);\,\,    v(t_n)=3^{k-n+1}\end {eqnarray}  

 We obtain from $(3.60)$ and $(3.57)$ that   $v_0(t_i)=f_k(t_i)\,\, (i=1,2,\cdots,n)$; this proves the first conclusion of Theorem $3.3.$

We have a natural imbedding $gr_{v_0}(R)\subseteq gr_{v_0} (V)$. We will prove first that the zero component of $gr_{v_0}(R)$ is $R/A$. This is equivalent to the fact that if  $x\in R$ then $v_0(x)=0$ iff $x\not \in A$ which we will now verify. Clearly, if $x\in A$ then $v_0(x)>0$ because $A$ is generated by the system of elements $<t_1,t_2,\cdots, t_n>$ and the elements of this system have values greater than zero. 

Assume now that $x\not \in A$. Then $x\not \in A_{n-1}$ and the asssumption of the induction implies that $v_1(x)=0$ and the homogeneous component of $x$ 
 in the ring $gr_{v_1}(R)$ belongs to the subring $R/A_{n-1}$ of $gr_{v_1}(R_{n-1})$ and it coincides with  the image $\bar{x}$ of $x$ under the homomorphism $R\longrightarrow R_{n-1}=R/A_{n-1}$. Since $x\not \in A$ we obtain that it image $\bar{x}$ is not contained in the ideal $(A/A_{n-1})=(t_n)R_{n-1}$. We recall now that $V/(t)\cong gr_{v_1}(R_{n-1})$ and obtain that the element $x$ does not belong to the ideal $(t,t_n)$. We obtain from Theorem  $3.2.$ that $v_0(x)=0$ and our claim is proven.

 The  system of elements    $\tilde{u}_1,\tilde{u}_2,\cdots,\tilde{u}_{n-1},\tilde{ t}_n,\tilde{t}$ is algebraically independent in the ring $gr_{v_0}(V)$ over the subring $R/A$. We obtain now from $(3.60)$ that $\tilde{t}_i=\tilde{u}_i+3^{n-i-1}\tilde{t}\,\, (i=1,2,\cdots, n-1)$ and conclude from this that  the system of homogeneous components $\tilde{t}, \tilde{t}_i\,\, (i=1,2,\cdots,n)$  is   also algebraically independent and generates $gr(V)$.   
The homogeneous components $ \tilde{t}_i\,\, (i=1,2,\cdots,n)$ belong to the subring $gr_{v_0}(R)\subseteq gr_{v_0}(V)$ as well as the subring $R/A$,  we conclude from this that these homogoneous components 
  are algebraically independent in the ring  $gr_{v_0}(R)$ over  subring  $R/A$.

It remains to prove that $gr_{v_0}(R)$ is generated by the elements $\tilde{t}_1,\tilde{t}_2,\cdots, \tilde{t}_n$. We consider now the natural homomorphism $\phi\colon R\longrightarrow R_1=R/(t_1)$. We have in $R_1$ the polycentral independent system $t_2,t_3,\cdots, t_n$ and the weight function obtained by restriction of $f_k$ on the subsystem $t_2,t_3,\cdots,t_n$

\begin{equation} f_k(t_2)=3^{k-1}, f_k(t_3)=3^{k-3},\cdots, f_k(t_n)=3^{k-n+1}\end{equation}

The assumption of the induction yields  that this weight function extends to a $t$-adic pseudovaluation $v_1$  in the ring $R_1=R/(t_1)$ with associated graded ring isomorphic to the polynomial ring in the system of homogeneous components of the elements  $t_2,t_3,\cdots,t_n$ over the  ring $R_1/(t_2,t_3,\cdots,t_n)\cong R/A$; for every $2\leq i\leq n$ we  denote by $t_i^{\ast}$ the  homogeneous component  of  the element $t_i$ in the ring $gr_{v_1}(R_1)$.

Lemma $3.8.$ shows that the pseudovaluation $v_1$ coincides with the pseudovaluation $\bar{v}$ defined by the natural homomorphism $\phi\colon  R\longrightarrow R/(t_1)=R_1$ via Proposition $2.3.$ We obtain also  from Proposition $2.3.$  that the homomorphism $\phi$ defines the homomorphism  of graded rings $\tilde{\phi}\colon gr_{v_0}(R)\longrightarrow gr_{v_1}(R_1)$; the kernel of this homomorphism is the ideal $gr_{v_0}(A_1)$ where $A_1=(t_1)$. The homogeneous component of degree zero in $gr(R)$  and in $gr(R_1)$ is $R/A$;   since $ker(\phi)\subseteq A$ we obtain that the restriction of  $\tilde{\phi}$ on $R/A$ defines an isomorphism between the homogeneous components of degree zero of $gr_{v_0}(R)$ and $gr_{v_1}(R_1)$.

For every $2\leq i\leq n$ we have $\phi(t_i)=t_i$ and $v_0(t_i)=v_1(t_i)=3^{k-i+1}\,\, (i=2,3,\cdots,n)$. Since  the weights of the elements $t_i\,\, (i=2,3,\cdots,n)$ in $R$ and $R_1$ coincide we  obtain that for every $2\leq i\leq n$ $\tilde{\phi}$  maps the homogeneous component  $\tilde{t}_i\in gr_{v_0}(R)$ on the  homogeneous component $t_i^{\ast}\in gr_{v_1}(R_1)$. We see that  the  homomorphism $\tilde{\phi}$ is an epimorphism of the ring $gr(R)$ on the polynomial ring $(R/A)[t^{\ast}_2,t^{\ast}_3,\cdots, t^{\ast}_n]$ which is isomorphic to the polynomial subring $(R/A)[\tilde{t}_2,\tilde{t}_3,\cdots,\tilde{t}_n]\subseteq gr(R)$.  We obtain therefore that 
$gr(R)$ is a direct sum of the subring $(R/A)[\tilde{t}_2,\tilde{t}_3,\cdots,\tilde{t}_n]$ and the ideal $gr(A_1)$:

\begin{equation} gr(R)\cong (R/A)[\tilde{t}_2,\tilde{t}_3,\cdots,\tilde{t}_n]+gr(A_1)\end{equation}

 We will now show that the  ideal $gr(A_1)$ of $gr(R)$ is generated by the homogeneous component $\tilde{t}_1$.

Let $x=rt_1$ be an element of $A_1$. If $v_0(r)+v_0(t_1)>v_0(x)$ then we would have in $gr(R)$  $\tilde{r}\tilde{t}_1=0$ for the homogeneous components of $r$ and $t_1$. This is impossible because $gr(V)\cong (R/A)[\tilde{t}_1,\tilde{t}_2,\cdots, \tilde{t}_n,t]$ so the homogeneous component of $t_1$ is regular in $gr(V)$  and because of this it is regular in the subring $gr(R)\subseteq gr(V)$.
We see that $v_0(x)=v_0(r)+v_0(t)$ and we obtain $\tilde{x}=\tilde{r}\tilde{t}_1$. 
This proves that 

\begin{equation} gr(A_1)\cong (\tilde{t}_1)gr(R)\end{equation}

Let $y$ be a homogeneous element of degree $l$ in $gr(R)$. We obtain from $(3.62)$ and $(3.63)$ that there exist a homogeneous element  
$y_1\in (R/A)[\tilde{t}_2,\tilde{t}_3,\cdots,\tilde{t}_n]$ of degree $l$  such that $y=y_1+\tilde{t}_1y_2$ where $y_2\in gr(R)$ is a homogeneous element  of degree  $l_ 1<l$  because the degree of $t_1y_2$ is $l$.  We can assume that $y_2$ has a representation as a polynomial in $ \tilde{t}_1,\tilde{t}_2,\cdots,\tilde{t}_n$ over $R/A$ and we obtain from this that $y$ is also a polynomial.

This completes the proof of Theorem $3.3.$\bigskip

 {\bf 3.8.} {\bf Theorem I}. {\it Let $R$ be a ring, $T=<t_1,t_2,\cdots, t_n>$  be an independent  polycentral system in $R$.}

 {\it The ideal $A$ generated by the system $t$ is residually nilpotent. If $\tilde{R}$ is the completion of $R$ in the topology defined by this ideal and $\mathcal{X}$ is a system of coset representatives for the elements of the quotient ring $R/A$ then every element $x\in \tilde{R}$ has a unique representation }

\begin{equation} x=\sum_{n=0}^{\infty} \lambda_n\pi_n\end{equation}

{\it where $\lambda_n\in \mathcal{X}\,\, (n=0,1,\cdots) $ 
  $\pi_n$ are standard monomials  on $T$, and $lim_{n\to\infty}v(\pi_n)=\infty$.}\bigskip
 
{\bf Proof of Theorem I.}   Let    $x\in \bigcap_{i=1}^{\infty} A^i$; we    obtain from Theorem $3.3.$  that $v_0(x)=\infty$,  hence $x=0$. This proves the residual nilpotence of $A$.

 We  will now prove the second statement. Theorem $3.3.$ yields that  $gr_{v_0}\cong (R/A)[\tilde{t}_1,\tilde{t}_2,\cdots,\tilde{t}_n]$. We obtain from this and Lemma $3.7.$ that every element of the completion of $R$ in the $v_0$-topology has a unique representation $(3.43)$. Since the $v_0$ and $\rho$ topology are equivalent we obtain that every element of $\tilde{R}$ has a unique representation $(3.64)$. 

The proof of Theorem I is complete. \bigskip

 {\bf 3.9. Proof of Theorem III.} We assume that the assertion is proven for the ring $\bar{R}=R/(t_1)$ and the  polycentral system $\bar{T}=<t_2,t_3,\cdots,t_n>$; we have already pointed out in subsection $3.6.$ that the initial step of the induction is true.   Lemma $3.9.$ now implies that the homogeneous components $\tilde{t}\,\, (t\in T)$ are central in $gr_v(R)$. Lemma $3.6.$  implies that for every element $r$ there exists in $\tilde{R}_v$ a representation $(3.43)$. Since the $v$-topology is equivalent to the $\rho$-topology we obtain now from Theorem $3.3.$  that it is equivalent also to the $v_0$-topology, and we conclude from this that representation $(3.43)$ is unique in $\tilde{R}_v$. Finally, we  obtain now from Lemma $3.7.$  that $gr_v(R)\cong (R/A)[\tilde{t}_1,\tilde{t}_2,\cdots,\tilde{t}_n]$. 

The uniqueness of $v$ follows from Proposition $2.4.$

This completes the proof of Theorem III.\bigskip

{\bf Corollary 3.8.}  {\it Let $v_m$ be the pseudovaluation  defined in the quotient ring $R_m=R/A_m$ by the independent polycentral system $t_{m+1},t_{m+2},\cdots, t_n$ and the restriction of the function $f$ on this system.   Then the epimorphism $\psi\colon R\longrightarrow R_m$ defines in a natural way an epimorphism $\tilde{\psi}\colon gr_v(R)\longrightarrow gr_{v_m}(R_m)$.}\bigskip

{\bf Proof.}  Lemma $3.6.$ implies that this statement is true for the homomorphism $\psi_1\colon R\longrightarrow R_1=R/(t_1)$. We can assume that it is true for the epimorphism $\tau\colon R_1\longrightarrow R_m$, and obtain that it is true for the product of homomorphisms  $\psi_1\tau\colon R\longrightarrow R_m$. 
This completes the proof.\bigskip

{\bf Theorem II.} {\it Let $R$ be a ring, $T=<t_1,t_2,\cdots,t_n>$ be a polycentral independent system in $R$ which is composed from the central systems $T_1,T_2,\cdots,T_k$, $A$ be the ideal generated by the system $T$. Let $f$ be a function on  on $T$ whose values are  natural numbers and  $f(t_1)>2 f(t_2)\,\, \mbox{for}\,\, t_1\in T_i,t_2\in T_{i+1}\,\,  (i=1,2,\cdots,k-1)$.

 Then there exists a pseudovaluation $v$ of $R$ such that }

\begin{equation} v(t)=f(t)\,\ \mbox{if}\,\, (t\in T) \end{equation}

{\it and the graded ring $gr_v(R)$ is isomorphic to the polynomial ring \\
$(R/A)[\tilde{t}_1,\tilde{t}_2,\cdots,\tilde{t}_n]$, the topology defined in $R$ by this pseudovaluation is equivalent to the topology defined by the powers of the ideal $A$.  Furthermore,   $v$ is the  unique pseudovaluation such that $v(t)=f(t)\,\, (t\in T)$  and the graded ring associated to it is isomorphic 
to\\ $(R/A)[\tilde{t}_1,\tilde{t}_2,\cdots,\tilde{t}_n]$.}\bigskip

{\bf Proof. } Follows immediately from Theorem III.\bigskip

We finish this section with the following  corollary of the results of this section.\bigskip

{\bf Corollary 3.9.} {\it  Let $R$ be a ring with a discrete pseudovaluation $\rho$, $gr_{\rho}(R)$  be the associated graded ring. Assume that there exists in $gr(R)$ an independent polycentral sustem  $T$, let $A$ be the ideal generated in $gr(R)$ by this system. Then  there exists in $R$ a discrete pseudovaluation $v$ such that the graded ring $gr_v(R)$ is isomorphic to a subring of the Laurent polynomial ring in the system of variables $T,t,t^{-1}$ over the subring $(gr(R))/A$}.  \bigskip

{\bf Proof.} We extend the pseudovaluation $\rho$ to  the ring $R[t,t^{-1}]$ and apply Proposition $2.8.$ We obtain that $V/(t)\cong gr(R)$, so the system of elements $t,T$ is an independent polycentral system in $V$. Theorem II implies that there exists in $V$ a discrete pseudovaluation $v$ such that $gr_v(V)\cong (gr(R)/A)[T,t]$. 

Since the  ring $R[t,t^{-1}]$ is isomorphic to the ring of fractions of $V$ with repect to the subsemigroup $S$ generated by the element $t$ the pseudovaluation $v$ extends in a natural way to the ring to $R[t,t^{-1}]$ and the graded ring of $R[t,t^{-1}]$ associated to $v$ is isomorphic to the ring of fractions 
 $(gr(R)/A)[T,t,t^{-1}]$ of $(gr(R)/A)[T,t]$. Since $R\subseteq R[t,t^{-1}]$ we obtain that $gr_v(R)$ is the subring of this ring of fractions, and the proof is complete.

\setcounter{section} {4}
\setcounter{equation} {0}
\section*{\center \S 4. }

 {\bf 4.1.} We will prove  in this section Theorem IV  which is an   application of Theorems I-III   to the group ring of torsion free nilpotent group.  
Its proof is based on the fact that if $H$ is a finitely generated torsion free nilpotent group then every central series with torsion free factors provides in a natural way an independent polycentral system (see Proposition  $4.1.$,  statement ii). If $H$ is not finitely generated  we will reduce the proof to the finitely generated case by direct limit arguments.\bigskip

{\bf  Lemma 4.1.} i) {\it Let $H$ be a torsion free nilpotent group which contains no elements of infinite $p$-height.  Then }

i) {\it all the factors of the upper central series   contain no elements of infinite $p$-height.}

 ii) {\it  Let $H$ be a torsion free nilpotent group which has a central  series} 
\begin{equation} H=U_1\supseteq U_2\supseteq\cdots\supseteq U_{k-1}\supseteq U_k=1\end{equation}

{\it  with unit intersection and all the factors $U_i/U_{i+1}$ are torsion free abelian group without elements of infinite $p$-height. Then $H$ contains no elements of infinite $p$-height.}  \bigskip

 {\bf Proof.} i) Let $H$ be a torsion free nilpotent group without  elements of infinite $p$-height, $Z$ be its center. We will prove that  the quotient group $H/Z$ contains no elements of infinite $p$-height. Statement i) will follow from this by an induction on the nilpotency classs of $H$. 

 Assume that there exists a non-central  element $h\in H$ which has an infinite $p$-height in $H/Z$. This means that for every given $n$  there exist $u\in H$ and $z_n\in Z$ such that $h=u^{p^n}z_n$. We pick an arbitrary element $g\in H$ which does not commute with $h$ and obtain $1\not= [g,h]=[g,u^{p^n}]$. Now assume that the nilpotency class of $H$ is $c$, and pick $n=mc$ where $k$ is an arbitrary number. Since the element  $h=[g,u^{p^n}]$
is a products of two elements $g^{-1}u^{-p^{mc}}g$ and $u^{p^{mc}}$ we obtain from Malcev's Lemma (see, for instance, Hartley  [2], Lemma $2.4.2.$) that $h=v^{p^m}$ for a suitable $v\in H$. Since $m$ was an arbitrary integer we conclude that the element $h$ has an infinite $p$-height. We obtained a contradiction and the proof of statement i) is complete.

ii) We can assume by an  induction argument  that the quotient group $H/U_{k-1}$ contains no elements of infinite $p$-height. Assume that  $H$ has an element $h$ of infinite $p$-height;  if the  image $\bar{h}$ of $h$ in  $H/U_{k-1}$ is non-unit it must have  infinite $p$-height in $\bar{H}$. We conclude therefore that $h\in U_{k-1}$. Once again, as in the proof of statement i), if there exist elements 
$h_1,h_2,\cdots,h_n\in H$ such that $h=h_1^{p^{mc}}h_2^{p^{mc}}\cdots h_n^{p^{mc}}$
then there exists $u\in H$ such that $h=u^{p^m}$. Since the quotient group $H/U_{k-1}$ is torsion free we conclude that $u\in U_{k-1}$. Since $h$ has infinite $p$-height this  contradicts the asssumption that all the factors contain no elements of infinite $p$-height and the proof is complete.\bigskip

{\bf 4.2.} We need a few   elementary facts about Lie algebras of groups without elements of infinite $p$-height. Let $H$ be a torsion free abelian group without elements of infinite $p$-height, $E$ be a system of elements which forms a basis for the vector space $H/H^p$. The system of elements $E$ forms a free system of generators for every quotient group $H/H^{p^m}\,\, (m=1,2,\cdots)$ and for the inverse limit of this system of groups. This inverse limit is isomorphic to the vector space 
$C_p\otimes H$ where $C_p$ is the subring of $p$-integers of the field of rationals, and the system $E$ is also the basis of this  vector space.
If $H_0$ denote the subgroup generated by $E$ then the quotient group $H/H_0$ is a periodic group without elements of order $p$.    Let $E_j\,\, (j\in J)$ be the direct  system of all the finite subsystems of $E$, and  for every $j\in J$ let $H_{0j}$ be the subgroup generated by $E_j$. Then $H$ is isomorphic to the direct limit of the system of group $H_{0j}\,\, (j\in J)$. \bigskip

{\bf Lemma 4.2. } {\it Let  $H$ be a torsion free abelian group without elements of infinite $p$-height, $E$ be a system of elements which forms a basis of the vector space 
$H/H^p$, $K$ be a field of characteristic $p$. Then }

i) { \it The restricted Lie algebra $L_p(H,H_i)$ associated to the $p$-series 

\begin{equation} H\supseteq H^p\supseteq\cdots\end{equation} 

is free abelian Lie algebra with system of generators $\tilde{E}$ formed by homogeneous components of elements $e\in E$.}

ii) {\it The graded ring of $KH$ associated to the filtration 
$\omega^i(KH)\,\, (i=1,2,\cdots)$ is isomorphic to the symmetric  algebra $K[V]$ of the vector 
space $V=H/H^p$.}

iii) {\it  Algebra $L_p(H)$ is the direct  limit of the system of free abelian algebras  $L_p(\tilde{E}_i) \,\, (i\in I)$ and the graded ring $gr(KH)$ is the direct limit of the symmetric  algebras $K[V_i]\,\, (i\in I)$.}\bigskip

{\bf Proof.} i) The group $H$ is an inverse limit of the system of groups $H^{p^n}\,\, (n=1,2,\cdots)$. For every given $n$ the Lie algebra of the group  $H/H^{p^n}$ is an abelian algebra of exponent $p^n$ and it is freely generated by the first factor $H/H^p$. Since the system $E$ forms a basis of the first factor we obtain from this statement i).

ii) The group ring $KH$ is an inverse limit of the system of rings $KH/\omega^i(KH)\,\, (i=1,2,\cdots)$ where  every of these rings is generated by the system $E$, and is isomorphic to the quotient ring of the polynomial ring $K[E]$ by the ideal $(E)^i$ where $(E)$ is the ideal generated by the system of elements $E$ and the proof can be completed easily.\bigskip

 Now let  $f$ be an arbitrary function defined on the basis $E$ of the  vector space $V=H/H^p$ with values in the set of the  natural numbers, $K$ be a field of characteristic $p$.  This function defines in a natural  way  a valuation of the ring $K[V]$ and of its completion $\widetilde{K[V]}$. Since the group ring $KH$  imbeds into $\widetilde{K[V]}$ we obtain a valuation $v$ of $KH$ whose values on the elements of $V$ coincide with the values of the function $f$. We see  that the valuation $v$  can be defined by its  values on an arbitrary basis of the vector space $H/H^p$.  Since $H/H^p$ is the homogeneous component of degree $1$ in the ring $K[V]$ we obtain also that the weight function $f$ defines a weight function on the vector space $H/H^p$.
\bigskip

{\bf 4.3.}    Let  $H$ be  torsion free nilpotent group without elements of infinite $p$-height,  $(4. 1)$  be a central series in $H$ whose factors $U_i/U_{i+1}\,\, (i=1,2,\cdots,k-1)$ contain no elements of infinite $p$-height. For every $1\leq i\leq k-1$ let $E_i $ be a system of elements in $U_i$ which is a basis of the vector space $U_i/U_i^pU_i^{\prime}$ and let $E_i-1$ denote the system of elements $e-1\,\, (e\in E_i)$, and 

\begin{equation}E-1=< E_1-1,E_2-1,\cdots,E_{k-1}-1>\end{equation}

{\bf Proposition  4.1.} {\it  Let $H$ be a finitely generated torsion free nilpotent group,    $R$ be a ring. Then }

i) {\it The system of elements $E_1$ is an independent polycentral system in  the group ring $RH$.} 

ii) {\it Let $\pi$ be a central regular element in $R$ such that $\bigcap_{n=1}^{\infty}(\pi)^n=0$ and  the quotient ring  $K=R/(t)$ is a field of finite characteristic $p$. 
Then the  system of elements $\pi, E-1$  is an independent polycentral system in  the group ring $RH$.} \bigskip

{\bf Proof.}  We will prove statement ii); the proof of i) is obtained by an obvious simplification of the argument. 

The system of elements $E_{k-1}-1$ is a central  independent system in $KU_{k-1}$ by statement ii) of  Lemma $4.2.$   and it is central and independent in $KH$. The ideal    generated by this system coincides with the ideal $\omega(KU_{k-1})KH$ and the quotient ring $KH/(\omega(KU_{k-1})KH) $ is isomorphic to the group ring of the group $H/U_{k-1}$.   The induction argument now implies that the system $E-1$ is  polycentral and independent in $KH$. We obtain from this that the system $<\pi,E-1>$ is polycentral and independent in $RH$.   \bigskip

{\bf Theorem IV.} {\it  Let $H$ be a torsion free nilpotent group without elements of infinite $p$-height, $R$ be a ring which contains a central regular element $\pi$  such  that $\bigcap_{n=1}^{\infty}(\pi)^n=0$ and  the quotient ring  $K=R/(t)$ is a field of finite characteristic $p$.

Let $f$ be a weight function  on the system $<\pi,E-1>$ whose    values are natural numbers  and  }

\begin{eqnarray} f(\pi)&>&2f(e-1)\,\, ( e\in E_{k-1});\qquad \nonumber\\  
f(e_1-1)&>&2 f(e_2-1)\,\, (e_1\in E_{i+1}, e_2\in E_i)\,\, (1\leq i\leq k-2)\end{eqnarray}

{\it This weight function has a unique extension  to a $t$-adic pseudovaluation of $RH$ with  associated graded ring isomorphic to the polynomial ring 
$K[\tilde{\pi},\widetilde{(E-1)}]$ generated by the homogeneous components $\tilde{\pi}, \widetilde{(e-1)}\,\, (e\in E)$. The epimorphism  $\colon R\longrightarrow R/(\pi)$ defines an epimorphism of graded rings $\tilde{\phi}\colon K[\tilde{\pi},\widetilde{(E-1)}]\longrightarrow K[{\widetilde{(E-1)}}]$.} \bigskip

{\bf Proof.} The system of elements $E_i$ is a basis of the vector space $U_i/U_i^pU_{i+1}$, we have seen in the end of subsection $4.2.$ that   the weight function on the system $E_i$ extends in a natural way to a weight function  on the vector space $U_i/U_i^pU_{i+1}$; conversely a weight function on the vector space $U_i/U_i^pU_{i+1}$  defines a weight function on $E_i$.  Because of this it is possible to assume that for every $i$ the weight function $f$ is defined on $U_i/U_i^pU_{i+1}$. 
We pick now an arbitrary finitely generated subgroup $V\subseteq H$ such that $V=\sqrt{V}$ i.e. if $h^n\in V$ for an element $h\in H$ then  $h\in V$. We denote now $V_i=V\bigcap U_i\,\,(i=1,2,\cdots, k)$ and  obtain now in $ V$ a central series

\begin{equation} V=V_1\supseteq V_2\supseteq\cdots\supseteq V_{k-1}\supseteq V_k=1\end{equation} 

The condition $V=\sqrt{V}$ implies that for every $i$ the torsion free abelian subgroup  $V_i/V_{i+1}$ is pure in the group $U_i/U_{i+1}$ and that the vector space $M_i^{\prime}=V_i/V_i^pV_{i+1}$ naturally imbeds in the vector space $M_i=U_i/U_i^pU_{i+1}$. We obtain therefore a restriction of the weight function $f$ on the subspaces $V_i/V_i^pV_{i+1}\,\, (i=1,2,\cdots,k-1)$; we denote this restriction by $f^{\prime}$.  We consider the group ring $KV$ and obtain  from Proposition $4.1. $  
and Theorem II that there exists a valuation $v$  in $KV$ which extends the weight function  $f^{\prime}$
  such that the graded ring $gr_v(KV)$ is isomorphic to the  symmetric algebra   $K[M^{\prime}]$ where $M^{\prime}=\bigcup_{i=1}^{k-1}M_i^{\prime}$. Further, Theorem II yields that there exists a valuation $\rho$  on the ring $RV$ such that $\rho(\pi)=f(\pi)$, $\rho(V_i)=f(V_i)\,\, (i=1,2,\cdots,k-1)$, the graded ring $gr_{\rho}(RV)$ is isomorphic to the polynomial ring  $K[\tilde{\pi}, T]$ where 
$T$ is an arbitrary basis of $V$ over $K$, and the epimorphism  $RH\longrightarrow KH$ defines the epimorphism  of graded rings $K[\tilde{\pi}, T]\longrightarrow KT$. 

We recall now that $V$ was an arbitrary finitely generated subgroup of $H$ such that $V=\sqrt{V}$. We construct for every $V_j\,\, (j\in J)$ a valuation $v_j$ and the proof is completed now by a standard direct limit argument. \bigskip

We have the following immediate corollary of Theorem IV.\bigskip

{\bf Corollary 4.1.}  {\it Assume that the ring $R$ in Theorem IV is either the ring of integers $C$ or the ring of $p$-adic integers $\Omega$. Then the graded ring associated to the valuation $v$  is isomorphic  to the polynomial ring 
$Z_p[t,T]$   generated by the homogeneous components $\tilde{p}=t$, and all the homogeneous components $\widetilde{e-1}\,\, (e\in E)$.}

 {\it If $H$ is finitely generated  the system of elements $p, E_1-1,E_2-1,\cdots, E_{k-1}-1$ is an independent polycentral system in the group ring $CH$ or  $\Omega H$.}\bigskip

{\bf 4.4.} Let $K$ be a field of characteristic zero, $H$ be a torsion free nilpotent group.  The same type of an argument as in Theorem IV, with some simplifications, proves  the following result. \bigskip

{\bf Theorem IV$^{\prime}$.} {\it Let $H$ be a torsion free nilpotent group, $K$ be a field of characteristic zero. Let $(4.1)$ be a central series with torsion free factors, and let $E_i$ be a system elements in $U_i$  which is a maximal linearly independent system modulo $U_{i+1}$ and let $E=\bigcup E_i\,\, (i=1,2,\cdots,k-1)$.  Let $f$ be a weight function on on the system $E_1$ such that}

\begin{equation}  f(e_1-1)>2 f(e_2-1)\,\, (e_1\in E_{i+1}, e_2\in E_i)\,\, 
(1\leq i\leq k-2)\end{equation}

{\it This  weight function  defines filtration $A_j\,\ (j=0,1,\cdots)$  in $KH$   with zero intersection and  a valuation  with graded ring isomorphic to $K[E]$. 

If $H$ is finitely generated the system  $E$ is a polycentral independent system in $H$. } 

We point out that Theorem IV$^{\prime}$ gives a new proof   Hall-Hartley Theorem about the residual nilpotence of the augmentation ideal $\omega(KH)$ (see Passman $[13]$.  In fact the Hall-Hartley Theorem follows from Theorem
 IV$^{\prime}$ because 
$\omega^j(KH)\subseteq A_j(KH)\,\, (j=1,2,\cdots)$. 

The original proof of  Hall-Hartley Theorem was  based on extensive use of commutator calculus.

\setcounter{section}{5}
\setcounter{equation}{0}

\section*{ \center \S  5 }
 {\bf Lemma 5.1.} {\it Let $L$ be a finitely generated restricted Lie algebra. Assume that every element of $ L$ is nilpotent, i.e. for every $x\in L$ there exists a number $p^n$ such that $x^{{[p]}^n}=0$  and that the Lie algebra $L$ is nilpotent, i.e. $\gamma_c(L)=0$ for some number $c$. Then the algebra $L$ is finite,  and its order is a power of $p$.}\bigskip

{\bf Proof.} We will use induction by the nilpotency class of $L$.  The assertion is obvious if $L$ is abelian, and we assume that it is true for the quotient algebra of $L$ by its center $Z$. Since $L/Z$ is finite and $L$ is finitely generated we obtain (see [1], $2.7.5.$) that the subalgebra $Z$ is finitely generated. We obtain therefore that $L$ is an extension of a finite algebra $Z$ by a finite algebra $L/Z$ hence it is finite. The same argument shows that the order of $L$ is a power of $p$. \bigskip
 
{\bf Lemma 5.2.} {\it Let $L$ be a finitely generated restricted Lie algebra which contains an ideal $V$ which  is nilpotent as a Lie algebra.  Then there exists $n$ such that the ideal $V^{[p]^n}$ generated by all the elements $v^{[p]^n}\,\, (v\in V)$  is central  in $L$. If the quotient algebra $L/V$ is finite  then there exists $m\geq n$ such that  $V^{[p]^m}$ is either zero or a central free abelian subalgebra of finite rank and the quotient algebra $L/V^{[p]^n}$ is a restricted finite nilpotent algebra.}\bigskip

{\bf Proof.} Let $v$ be an arbitrary element of $V$. Then $[u,v]\in V$ for an arbitrary element $u\in U$. Since $V$ is nilpotent we can find a number $p^n$ such that $[u,\underbrace{v,\dots,v]
}_{p^n}=0$ for every $u\in L$, and hence $[u,v^{[p]^n}]=0$. The last equation means that every element $v^{[p]^n}$ is central.\par
If  $L/V$ is finite   then $V$ must be finitely generated. So $V^{[p]^n}$ is a finitely generated central subalgebra of $L$.  Since all the nilpotent elements in $V^{[p]^n}$ form a finite subalgebra we can take   number $m$  which is a suitable multiple of $n$  and to obtain from Lemma $2. 4.$ that $V^{[p]^m}$ is either zero or free abelian,   and the quotient algebra $V/V^{[p]^m}$ must be finite  by Lemma $5.1.$;   hence $L/V^{[p]^m}$ is finite    and the proof is complete. \bigskip

We need the following fact.\bigskip

{\bf Lemma 5.3.} {\it Let $H$ be a polycyclic-by-finite group. Then $H$ contains a normal subgroup $V$ of finite index such that the nilpotent radical $N$ of $V$ is torsion free and  the quotient group $V/N$ is free abelian. Moreover $V$ has
 an additional property that     every element of  it  centralizes $N$ modulo the dimension subgroup $N^{\prime}N^p$.}\bigskip

{\bf Proof.} The first statement follows from the Kolchin-Malcev Theorem (see Kargapolov and Merzlyakov [6], Theorem VII. $3.3.$, or Segal [17]), the second one is obtained by a routine argument.\bigskip

 {\bf Theorem  V.} {\it Let $H$ be an infinite    polycyclic-by-finite group with Hirsh number $r$.  Assume that there exists    a $p$-series $(1.1)$ with unit intersection such that the corresponding restricted Lie algebra $L_p(H,H_i)$ is finitely generated.    Then there exists a 
torsion free  normal  subgroup $F$   with index a power of $p$ such that the ideal $L_p(F,F_i)$ associated to the $p$-series $F_i=F\bigcap H_i\,\, (i=1,2,\cdots)$ is a restricted free abelian  subalgebra of the center of  $L_p(H,H_i)$
   with  index  a power of $p$ and rank $1\leq r_1\leq r$.

 Hence the   center $Z$ of $L_p(H,H_i)$  has a finite index which is a power of $p$ and $L_p(H,H_i)$ is a nilpotent Lie algebra.  }\bigskip

{\bf Proof of Theorem V. }  Let $V$ be a normal subgroup obtained in Lemma $5.3.$, $N$ be its nilpotent radical.   Since the quotient group $H/V$ is finite the ideal $L_p(V,V_i)$ has a finite index in $L_p(H,H_i)$ and we obtain  that the algebra $L_p(V,V_i)$ is finitely generated. Further the nilpotence of $N$ implies  via Corollary  $2.4.$ that the restricted Lie algebra $L_p(N,N_i)$ is  nilpotent as a Lie algebra.
Lemma $5. 2.$ implies that there exists a number $p^n$ such that the subalgebra $L_p(N,N_i)^{[p]^n}$ is  a  central subalgebra of $L_p(H,H_i)$.  Let $c$ be the nilpotency class of $N$.  Malcev's Lemma (see  Hartley [4]) implies that every element of the subgroup $N^{p^{nc}}$ has a form $x=y^{p^n}$ for a suitable $y\in N$. We consider now the normal subgroup   $W=N^{p^{nc}}$ generated by all the elements $h^{p^{nc}}\,\, (h\in N)$; every element of $L_p(W,W_i)$ is a homogeneous component of some  element $y^{p^n}$ so either its homogeneous component is $\tilde{y}^{[p]^n} $ or $\tilde{y}^{[p]^n} =0$. We    obtain from this that the subalgebra $L_p(W,W_i)$ is central in $L_p(V,V_i)$.  The quotient group $\bar{V}=V/W$ is an extension of the finite $p$-group  $N/N^{p^{nc}}$ by the free abelian group $\bar{V}/\bar{N}\cong V/N$, and $V$ acts trivially on the factor $\bar{N}/\bar{N}^{\prime}\bar{N}^p\cong N/N^{\prime}N^p$,  hence the group $\bar{V}$ is nilpotent. Corollary  $2.4.$  implies that the restricted Lie algebra $L_p(\bar{V},\bar{V}_i)$ is a nilpotent Lie algebra. We see that the restricted Lie algebra $L_p(V,V_i)$ is an extension of a central ideal $L_p(W,W_i)$ by the algebra  $L_p(\bar{V},\bar{V}_i)$  which is a nilpotent Lie algebra,    hence $L_p(V,V_i)$  is a nilpotent Lie algebra. Lemma $5.2.$ now implies that there exists a number $l$ such that the subalgebra 
$L_p(V,V_i)^{[p]^l}$ is central in $L_p(H,H_i)$. The index of $L_p(V,V_i)$ is finite because the subgroup $V$  has a finite index in $H$; since $L_p(V,V_i)^{[p]^l}$ has finite index in $L_p(V,V_i)$ we conclude that the index of $L_p(V,V_i)^{[p]^l}$ in $L_p(H,H_i)$ is finite.

We take now the subgroup $Q=V^{p^{lc}}$ and obtain by the same argument as above that the subalgebra $L_p(Q,Q_i)$ is central in $L_p(H,H_i)$ and has a finite index. Hence $L_p(Q,Q_i)$ is finitely generated once again  by [1], $2.7.5.$     Since $L_p(Q,Q_i)$ contains a finite number of nilpotent elements we can find a number 
$i_0$ such that the central subalgebra 

\begin{equation} \sum_{i\geq i_0}Q_i/Q_{i+1}\end{equation} 

 is either the zero subalgebra or is free abelian.  Subalgebra $(5.1)$ is in fact the restricted Lie algebra 
$L_p(R,R_i)$ of the subgroup $R=Q_{i_0}$ associated to the $p$-series $R_i=Q_{i_0}\bigcap H_i$ ; it   has a finite index in $L_p(H,H_i)$ because $Q$ has a finite index in $H$. If $L_p(R,R_i)=0$ then $\bigcap_{i=1}^{\infty}R_i=R$ which together with the relation $\bigcap_{i=1}^{\infty}R_i=1$ implies $R=1$. Since $H$ is infinite and $R$ has a finite index in it we obtain a contradiction. This means that  $L_p(R,R_i)$ is a central free abelian subalgebra  of rank $r_1\geq 1$. Proposition $2.5.$ implies that $r_1\leq r$. 

We apply now Corollary $2.1.$  and obtain a normal subgroup $F\supseteq R$ with index a power of $p$ such that $L_p(F,F_i)=L_p(R,R_i)$. The subgroup $F$ is 
torsion free because the subalgebra $L_p(F,F_i)$ is free abelian. We see that the subgroup $F$ satisfies all the conclusion of the theorem.

We prove now the last statement of the theorem.   Since $Z\supseteq L_p(F,F_i)$ we obtain that the index of $Z$ is a power of $p$. Further, since the the quotient algebra $L_p(H,H_i)/L_p(F,F_i)$ is a nilpotent Lie algebra and $L_p(F,F_i)$ is a central subalgebra we obtain that the Lie algebra $L_p(H,H_i)$ is nilpotent. 

The proof of Theorem V is complete.

\setcounter{section}{6}
\setcounter{equation}{0}

\section*{\center \S 6. Proof of Theorem VI.}

{\bf 6.1.}  {\bf Proposition 6.1.} {\it Let $H$ be a group, $U$ be its normal subgroup which satisfies a law $\omega(x_1,x_2,\cdots,x_n)=1$. Assume that $H\in \mbox{res}\,\mathcal{N}_p$ and let $V\supseteq U$ be the normal subgroup formed by all the elements which have an infinite $p$-height in the quotient group $H/U$. Then $V$ satisfies law $\omega$.}\bigskip

{\bf Proof.} Let $\phi_i\colon H\longrightarrow H_i$ be a homomorphism of $H$ on a finite  $p$-group $H_i$, $X_i$ be the image of a subset $X\subseteq H$ under this homomorphism. The group $V_i$ is an extension of the normal subgroup $U_i$ by the group $V_i/U_i$; since the elements of $V$ have an infinite $p$-height in the quotient group $V/U$,  all the non-unit elements of $V_i/U_i$ have an infinite $p$-height.  But the group  $V_i/U_i$ is a subgroup of $H_i$ which is a finite  $p$-group . Hence  
$V_i/U_i=1$, i.e. $V_i=U_i$. We obtain therefore that 
$\omega(V_i)=1$, which implies that  $\omega(V)=1$ and the assertion follows.\bigskip 

{\bf Corollary 6.1.} {\it Let $H$ be a finitely generated  nilpotent-by-finite group, $U$ be the unique maximal nilpotent normal subgroup of $H$. If $H$ is a residually-\{finite $p$-group\} then the index of $U$ is a power 
 of $p$}.\bigskip

{\bf Proof.} Let $V$ be as in Proposition $6.1.$ Then $V$ is nilpotent, so we obtain that $V=U$.  The definition of $V$ implies  that the order of an arbitrary  element of $H/V=H/U$  must be a power of $p$, so $H/U$ is a $p$-group.  \bigskip

 {\bf Proposition 6.2.} {\it Let $H$ be a finitely generated abelian-by-finite group with Hirsh number $r$ and without finite normal subgroups.  Assume that there exists a $p$-series 
\begin{equation} H=H_1\supseteq H_2\supseteq\cdots\end{equation}
with unit intersection 
such that the Lie algebra $L_p(H,H_i)$ is abelian of rank $r$. Then the topology defined by  $p$-series $(6.1.)$  is equivalent to the $p$-topology, 
  there exists a number $i_0$ such that $H_{i_0}\subseteq U$ where $U$ is the unique maximal free abelian  characteristic subgroup of $H$, the index of $U$ is a power of $p$ }.\bigskip

{\bf Proof.} The nilpotent radical $\rho(H)=U$ is a torsion free nilpotent group. Since it is abelian-by-finite group it must be abelian. If $C(U)$ is the centralizer of $U$ then the commutator subgroup of $C(U)$ is finite by Schur's Theorem, so $C(U)$  must be abelian because $H$ contains no finite normal subgroups.  Since $U=\rho(H)$ we obtained that $C(U)=U$.

Let $\bar{H}=H/U$, $\bar{H}_i\,\, (i=1,2,\cdots)$ be the image of the $p$-series $H_i$ in $\bar{H}$. We will prove now  that \begin{equation}\bigcap_{i=1}^{\infty}\bar{H}_i=1\end{equation}  In fact if 
$1\not =\bar{h}\in \bigcap_{i=1}^{\infty}\bar{H}_i$ then there exists $h\in H$ such that $h\not\in U$ but for every $i$ there exists $u_i\in U$ such that $hu_i\in H_i$. This implies that $[hu_i,u]=[h,u]\in H_i\,\, (i=1,2,\cdots)$ for every $u\in U$ and hence $[h,u]\in \bigcap_{i=1}^{\infty} H_i=1$. But the centralizer of $U$ coincides with $U$, hence $h\in U$ and we obtained a contradiction, and relation $(6.2)$ is proven.

Since the  group $\bar{H}$ is finite we obtain from $(6.2)$  that there exists $i_0$ such that $\bar{H}_{i_0}=1$ that is $H_{i_0}\subseteq U$; this proves also that the index of $U$ is a power of $p$ because the index of $H_{i_0}$ is a power of $p$.

It remains to prove that the topology defined by series $(6.1)$ is equivalent to the $p$-topology. The group $U$ is free abelian and the Lie algebra $L_p(U,U_i)$ associated to the $p$-series $U_i=U\bigcap H_i$ is abelian of rank $r$. Lemma $2.8.$ 
implies that the topology defined by this series in $U$ is equivalent to the topology defined by the lower $p$-series $U^{p^j}\,\, (j=1,2,\cdots)$. We can find  therefore for every $j$ a number $i=i(j)$ such that 
$H_i\subseteq U^{p^j}\subseteq M_j(H)$  which proves that the series $(6.1)$ defines the same topology as the  series $M_j(H)\,\, (j=1,2,\cdots)$. This completes the proof.\bigskip

{\bf Proposition 6.3.} {\it Let $H$ be a finitely generated abelian-by-finite group with Hirsch number $r$ which contains a $p$-series $(6.1)$  with unit intersection such that the algebra $L_p(H,H_i)$ is free abelian of rank $r$.  Then the topology defined by this $p$-series is equivalent to the $p$-topology and  there exists a number $i_0$ such that the subgroup $H_{i_0}$ is torsion free abelian.}

{\bf Proof.} Let $tor(H)$ be the unique maximal normal finite subgroup of $H$. Since the intersection of the terms of  series $(6.1)$ is trivial we can finite  a subgroup  $F=H_{i_0}$ such that $F\bigcap tor(H)=1$. The normal subgroup $F$ contains no finite normal subgroups and contains a series $F_i=F\bigcap H_i$ with unit intersection and  if $i\geq i_0$ we have $F_i=H_i$. The algebra $L_p(F,F_i)$ is naturally imbedded in $L_p(H,H_i)$ so it must be abelian. The quotient algebra $L_p(H,H_i)/L_p(U,U_i)$ is finite because the index of $F$ in $H$ is finite. We conclude from this that the rank of $L_p(F,F_i)$ is equal $r$, and it coincides with the Hirsch number of $F$. We apply now Proposition $6.2.$ to the group $F$ and series $F_i\,\, (i=1,2,\cdots)$ and obtain that there exists a number $i_1\geq i_0$ such that the subgroup $F_i=H_i$  is free abelian if $i\geq i_1$.

Proposition $6.2.$ implies that the topology defined in $F$ by the series $F_i\,\,(i=1,2,\cdots)$ is equivalent to the $p$-topology of this group. We obtain from this that for every subgroup $M_n(F)$ there exists $i(n)$ such that $F_i\subseteq M_n(F)\subseteq M_n(H)$ if $i\geq i(n)$. This together with the fact that $F_i=H_i\,\, (i\geq i_1)$ implies that $H_i\subseteq M_n(H)$ if 
$i\geq \mbox{max} \{i(n),i_1\}$ and the proof is complete.\bigskip

{\bf Proposition 6.4.} {\it Let $H$ be a finitely generated nilpotent-by-finite group  with Hirsch rank $r$, $U$ be the nilpotent radical of $H$.  Assume that there exists a 
$p$-series $H_i$ with unit intersection such 
that the Lie algebra $L_p(H,H_i)$ is abelian of rank $r$. Then the topology defined by this series is equivalent to the $p$-topology. Further,  there exists an index $i_0$ such that the subgroup $H_{i_0}$ is torsion free nilpotent and    the index of the nilpotent radical $U$ is a power of $p$.}\bigskip 

{\bf Proof.} {\it  We  consider  first the special case when $H$ contains no finite normal subgroups}. Let in this case  $V$ be a maximal  abelian $H$-invariant subgroup of $U$ and $W$ be the inverse image in  $H$ of  the maximal finite normal subgroup of $H/V$.  Clearly $V\subseteq W$ and  $W$ is an abelian-by-finite normal subgroup of $H$ which  contains no finite normal subgroups. Let $\rho(W)$ be the nilpotent radical of $W$. Then $\rho(W)\supseteq V$ and   $\rho(W)$ is the unique maximal abelian normal subgroup of $W$; we obtain from this that $\rho(W)\subseteq U$. Since $V$ is a maximal abelian  $H$-invariant  subgroup of $U$ and $\rho(W)\supseteq V$ we conclude that $\rho(W)=V$.  The group $W$ contains a $p$-series $W_i=H\bigcap W\,\, (i=1,2,\cdots)$,  the  Lie algebra $L_p(W,W_i)$ is abelian of rank $r_1$ equal to the Hirsh number of $W$ by Proposition $2.6.$  Proposition $6.2$ now implies that the index $(W:V)$ is a power of $p$ and the topology defined by the series $W_i$ in $W$ is equivalent to the $p$-topology in $W$. 

Let $\bar{H}=H/W$.  The definition of $W$ implies that  $\bar{H}$ contains no finite normal subgroups; we obtain from Proposition $2.6.$ that   the image of the $p$-series $H_i\,\, (i=1,2,\cdots)$  is  a $p$-series 
$\bar{H}_i$ with the associated Lie algebra $L_p(\bar{H},\bar{H}_i)$ free abelian of rank $r-r_1$ which is equal to the Hirsch number of $\bar{H}$. We can assume by induction on the Hirsch number that the topology defined by the $p$-series $\bar{H}_i$ in $\bar{H}$ is equivalent to the $p$-topology. This, together with the fact that the topology defined by $W_i$ in $W$ is equivalent to the $p$-topology in $W$ implies via Lemma  $2.7.$ that the topology defined by the $p$-series $H_i\,\, (i=1,2,\cdots)$ is equivalent to the $p$-topology. 

Since   $H$ is a residually finite $p$-group and $U$ is the nilpotent radical of it we can apply  now   Corollary $6.1$ and to conclude that $H/U$ is a finite $p$-group. Hence  there exists $n$ such that $M_n(H)\subseteq U$. Since the series $H_i\,\, (i=1,2,\cdots)$ and $M_i(H)$ define the same topology we can find $i_0$ such that $H_{i_0}\subseteq M_n(H)\subseteq U$. 

This completes the proof   for the special case. 

Now consider  the general case.  Let    $tor(H)$ be the unique maximal finite  normal subgroup of $H$, $F=H_{i_1}$ be a term of series $(6.1)$ that have a unit intersection with $tor(H)$. Then $F$ contains no finite normal subgroups and the assertion now follows from the proven special case by the same argument that was used in Proposition $6.3.$\bigskip

{\bf Proposition 6.5.}  {\it Let $H$ be a finitely generated nilpotent-by-finite group with Hirsch number $r$. Assume that there exists a 
$p$-series $H_i$ such that the Lie algebra $L_p(H,H_i)$ is abelian of rank $r$.
Let $F$ be a subgroup of $H$ with Hirsch rank $k$, $F_i=H\bigcap H_i\,\, (i=1,2,\cdots)$. Then the subalgebra $L_p(F,F_i)$ of $L_p(H,H_i)$ has rank $k$. }\bigskip

{\bf Proof.} Let $H_{i_0}=Q$ be the torsion free nilpotent normal subgroup of $H$,  obtained in Proposition $6.4.$, $N=F\bigcap Q$,   $Q_i=Q\bigcap H_i, N_i=N\bigcap H_i\,\, (i=1,2,\cdots)$. The algebra $L_p(Q,Q_i)$ is abelian of rank $r$ by Proposition $2.6.$; further,   $N_i=Q_i\bigcap N\,\, (i=1,2,\cdots)$. Since $N$ has a finite index in $F$ the subalgebra $L_p(N,N_i)$ has a finite index in $L_p(F,F_i)$, so it is enough to prove that $L_p(N,N_i)$ has rank $k$. We see that we can assume from the very beginning that the group $H$ is torsion free nilpotent, and $F$ is a subgroup of it.

Let  $V$ be an infinite cyclic central subgroup of $H$ such that the quotient group $H/V$ is torsion free, $T=V\bigcap F,   T_i=T\bigcap H_i\,\, (i=1,2,\cdots)$. 
We consider now the    natural homomorphism $H\longrightarrow H/V=\bar{H}$ which defines the homomorphism $L_p(H,H_i)\longrightarrow L_p(\bar{H},\bar{H}_i)$ where $\bar{H}_i=(H_iV)/V\,\, (i=1,2,\cdots)$. The restriction of this homomorphism on $F$ defines a homomorphism $F\longrightarrow F/T$ and  a homomorphism of Lie algebras $L_p(F,F_i)\longrightarrow L_p(\bar{F}/\bar{F}_i)$ where $\bar{F}=F/T, \bar{F}_i=(F_iT)/T\,\, (i=1,2,\cdots)$; the kernel of the last homomorphism is a Lie subalgebra $L_p(T,T_i)$ of $L_p(F,F_i)$. 

We have two possible cases: $T$ is the unit subgroup, or $T$ is an infinite cyclic subgroup  of finite index in $V$.  In the first case we obtain that $F\cong \bar{F}$ and it is naturally imbedded in the quotient group $\bar{H}=H/V$, and $L_p(\bar{F},\bar{F}_i)$ is the Lie algebra of the subgroup $\bar{F}=(FV)/V$   associated to the $p$-series $\bar{F}_i=\bar{F}\bigcap \bar{H}_i\,\, (i=1,2,\cdots)$. The Hirsch number of $\bar{H}$ is $r-1$ and we can assume by induction that the assertion holds for this case,  that is the rank of the algebra $L_p(\bar{F},\bar{F}_i)\cong L_p(F,F_i)$ is equal to its Hirsh number, so it is equal to $k$.

We consider now the second case. In this case the algebra $L_p(F,F_i)$ is an extension of the algebra $L_p(T,T_i)$ by the Lie algebra $L_p(\bar{F}, \bar{F}_i)$. The algebra $L_p(T,T_i)$ is an abelian subalgebra of $L_p(V,V_i)$, its index is finite because   $T$ has a finite index in $V$.  Hence $L_p(T,T_i)$ is an abelian algebra of rank $1$. The rank of $L_p(\bar{F},\bar{F}_i)$ is $k-1$ by the induction argument. We obtain from this that the rank of $L_p(F,F_i)$ is $k$ and the proof is complete.\bigskip

{\bf 6.2.} {\bf Theorem VI.} {\it Let $H$ be a polycyclic group with Hirsch number $r$. Assume that there exists a $p$-series $H_i\,\, (i=1,2,\cdots)$ with unit intersection such that $L_p(H,H_i)$ is  abelian of rank $r$. }

i) {\it Let $U$ be an arbitrary subgroup of $H$ with Hirsch number $k$, $U_i=U\bigcap H_i\,\, (i=1,2,\cdots)$. Then the algebra $L_p(U,U_i)$ is abelian  of rank $k$.} 

ii)  {\it  Let $U$ be a normal subgroup of $H$ with Hirsch number $k$,   $\bar{H}_i$ be the image of the subgroup $H_i$ in $H/U$. The subgroup $\bigcap_{i=1}^{\infty} \bar{H}_i$ is finite and the algebra $L_p(\bar{H},\bar{H}_i)$ is abelian of rank $r-k$. In particular, if
 $\bar{H}=H/U$ contains no finite normal subgroups then
$\bigcap_{i=1}^{\infty}\bar{H}_i=1$ and $\bar{H}$  is a  residually
 \{finite $p$-group\}. }

iii) {\it  Let $U$ be a normal subgroup of $H$. If  $\bar{H}=H/U$ is a residually \{finite $p$-group\} then $\bigcap_{i=1}^{\infty}\bar{H}_i=1$}.

iv) {\it Let $W$ be the unique maximal normal nilpotent-by-finite subgroup of $H$. Then $W$ is an extension of a torsion free nilpotent group by a finite $p$-group.  The quotient group $H/W$ is an extention of a free abelian group by a finite $p$-group.}

v) {\it The topology defined in $H$ by the $p$-series $H_i\,\, (i=1,2,\cdots)$  is equivalent to the $p$-topology. 
The topology defined in an arbitrary subgroup $U$   by the series  $U_i=H_i\bigcap U\,\, (i=1,2,\cdots)$ and    $M_n(H)\bigcap U\,\, (n=1,2,\cdots) $ are equivalent to the $p$-topology in $U$.}

vi)  {\it There exists an index $i_0$ such that if $i\geq i_0$ then the   subgroup  $Q=H_i$, contains  a torsion free nilpotent subgroup $N$  which is invariant in $H$,     $Q/N$ is free abelian, }  

\begin{equation}H/N^{\prime}N^p\in res\, \mathcal{N}_p\end{equation}
{\it Clearly, $H/Q$ is a finite $p$-group.}

vii) {\it Let $F\supseteq S$ be two normal subgroups in $H$ such that $H/F$ and $F/S$ are residually \{finite $p$-groups\}. Then $H/S$ is a residually \{finite $p$-group\}}. 

viii)  {\it  Let 

 \begin{equation} H=H^{\ast}_1\supseteq H^{\ast}_2\supseteq\cdots\end{equation}   

be a series in $H$  with finitely generated associated graded Lie algebra $L_p(H,H^{\ast}_i)$. If the  topology defined by series $(6.4)$ is equivalent to the $p$-topology then the center of  $L_p(H,H^{\ast}_i)$ has rank $r$. Moreover, there exists a number $k$ such that if  $U=H_i \,\, (i\geq k)$ then the subalgebra $L_p(U,U_i^{\ast})\cong \sum_{i\geq k} H_i/H_{i+1}$ associated to the $p$-series $U_i^{\ast}=U\bigcap H_i^{\ast}\,\,  (i=1,2,\cdots)$ is a  central free abelian subalgebra of rank $r$.}       \bigskip

{\bf Proof of Theorem VI.}    Statement ii)  follows from Propositions $2.6.$ and $2.5.$ 

We  prove now statement i). Let $r_1$ be the Hirsch number of the unique maximal nilpotent-by-finite subgroup $W$.  Proposition $2.6.$ implies that the rank of the algebra $L_p(W,W_i)$ is $r_1$. We  consider now  the  subgroup of $F=U\bigcap W$ and  obtain from Proposition $6.5.$ that the rank of the algebra $L_p(F,F_i)$ is equal to the Hirsch number of $F$, say $l$.  Now consider the natural homomorphism $H\longrightarrow H/W$ and its restriction $U\longrightarrow \bar{U}=U/F$.  Statement ii) implies once again that the rank of the algebra $L_p(\bar{H},\bar{H}_i)$ is equal to the Hirsch number $r-l$ of $\bar{H}$; since the group $\bar{H}$ is abelian-by-finite we obtain from    Proposition $6.5.$  that  the rank of $L_p(\bar{U},\bar{U}_i)$ is equal to the Hirsch number $k-l$ of $\bar{U}$.  We see that the restricted abelian Lie algebra $L_p(U,U_i)$ is an extension of an algebra $L_p(F,F_i)$ with rank $l$ by the algebra $L_p(U,U_i)$ with rank $k-l$. This implies that the rank of $L_p(U,U_i)$ is $k$ which proves ii).\bigskip

{\bf The proof of statements iv and v).}  The quotient group $\bar{H}=H/W$ is abelian-by-finite and it does not contain   finite normal subgroups. Statement ii) implies that $L_p(\bar{H},\bar{H}_i)$ is abelian with rank equal to the Hirsch number of $\bar{H}$ which is  rank $r-r_1$,  where $r_1$ is the Hirsch number of $W$ and $\bigcap_{i=1}^{\infty}\bar{H}_i=1$; we obtain from Proposition $6.2.$  that  the group $\bar{H}$ is an extension of a  free abelian group by a finite $p$-group. This completes the proof of statement  iv).\bigskip

 Corollary $6.1.$  implies that  $W$ is an extension of a torsion free nilpotent group by a finite $p$-group, and  
the topology defined by the series $W_i=H_i\bigcap W\,\, (i=1,2,\cdots)$ is equivalent to the $p$-topology of $W$.  Proposition $6.3.$ implies that the topology defined in $\bar{H}$ by series $\bar{H}_i\,\, (i=1,2,\cdots)$ is equivalent to the
 $p$-topology in $\bar{H}$.   We obtain  now  from Lemma  $2.7.$  that the topology defined in $H$ by series $H_i\,\, (i=1,2,\cdots)$ is equivalent to the $p$-topology. We can apply this result to the subgroup $U$ because $U$ contains a $p$-series $U_i=U\bigcap H_i\,\, (i=1,2,\cdots)$ with unit intersection and with associated graded  algebra $L_p(U,U_i)$   abelian of  rank equal to the Hirsch number of $U$; the last fact  follows from statement i).  We obtain  that the topology defined in $U$ by the series $U_i\,\, (i=1,2,\cdots)$  is equivalent to the $p$-topology. 

We consider now the series $M_n(H)\bigcap U\,\, (n=1,2,\cdots)$. We have already proven that there exists 
a number $i(n)$ such that $U_{i(n)}\subseteq M_n(U)$. Hence 

\begin{equation} (M_{i(n)}\bigcap U)\subseteq (H_{i(n)}\bigcap U)\subseteq
 U_{i(n)}\subseteq M_n(U)\end{equation}

which proves that the topology defined by the series $M_n(H)\bigcap U\,\, (n=1,2,\cdots)$ is equivalent to the $p$-topology in $U$.

This completes the proofs of statement  v).

To  prove statements vi) and vii) we need the following lemma.\bigskip

{\bf Lemma 6.1.} {\it Let $H$ be a group with Hirsch number $r$. Assume that there exists  a $p$-series $(6.1)$ with unit intersection and the restricted Lie algebra $L_p(H,H_i)$  abelian of rank $r$. Let $R$ be a normal subgroup such that the quotient group $H/R$ is an extension of a free abelian group by a finite $p$-group.  Then $H/M_n(R)\in res\, \mathcal{N}_p$ for every $n$.}\bigskip

{\bf Proof.} Let $H_0\supseteq R$ be a normal subgroup of index $p^k$ such that the quotient group $H_0/R$ is free abelian. The group $H_0$ has Hirsch number $r$ and it contains a $p$-series $H_{0i}= H_0\bigcap H_i\,\, (i=1,2,\cdots)$ with unit intersection such that the  algebra $L_p(H_0,H_{0i})$    is abelian of rank $r$.  Since $H/H_0$ is a finite $p$-group  it is enough to prove that $H_0/M_n(R)\in res\,\mathcal{N}_p$. We see that we can assume from the very beginning that $H=H_0$, i.e. $H/R$ is free abelian.

  Statement i) of Theorem VI implies that the group $R$ has a $p$-series $R_i=H_i\bigcap R\,\, (i=1,2,\cdots)$ such that the algebra $L_p(R,R_i)$ has rank equal to the Hirsch number of $R$;  statement v)  implies that we can find $m$ such $H_m\bigcap R=R_m\subseteq M_n(R)$. Let  $\bar{H}=H/H_m$, $\bar{R}=R/R_m$.  Since $H^{\prime}\subseteq R$ we obtain 

\begin{equation}[R,\underbrace{H,H,\cdots,H}_{m-1}]\subseteq (R\bigcap H_m)=R_m\end{equation}

and then 

\begin{equation} [\bar{R},\underbrace{\bar{H},\bar{H},\cdots,\bar{H}}_{m-1}]=1\end{equation}

Since $\bar{H}^{\prime}\subseteq  \bar{R}$ we obtain from $(6.7)$ that the group $\bar{H}$ is nilpotent. The group $H/M_n(R)$ is a homomorphic image of $\bar{H}$, hence it is also nilpotent. On the other hand, the group $H/M_n(R)$ is an extension of a finite $p$-group $R/M_n(R)$ by a free abelian group $ H/R$; this together with the nilpotency of this group implies that 
$H/M_n(R)\in \mbox{res}\, \mathcal{N}_p$. This completes the proof of the lemma.\bigskip

{\bf Proof  of statement vi) of Theorem VI. }  Let once again $W$ be the nilpotent-by-finite radical of $H$. Since the rank of the Lie algebra $L_p(W,W_i)$ is equal to the Hirsch number of $W$ by statement i)  Proposition $6.4.$ implies that there exists in  $H$ a normal subgroup $W_{i_0}=H_{i_0}\bigcap W$ such that $W_{i_0}$ is torsion free nilpotent and hence the group $W$ is an extension of a torsion free nilpotent group  by a finite $p$-group. Hence we  can find $l$ such that $M_l(W)$ is a torsion free nilpotent group.  On the other hand we  obtain from   Lemma $6.1.$ that the group $\bar{H}=H/M_l(W)\in \mbox{res}\,\mathcal{N}_p$  so $\bar{H}$    is an extension of a finite $p$-group $\bar{W}=W/M_l(W)$ by the group $\bar{H}/\bar{W}\cong H/W$, which contains no finite normal subgroups,  and it is also an extension of a free abelian group by a finite $p$-group;  hence $\bar{W}$ is  the unique maximal finite normal subgroup of $\bar{H}$. Since $\bar{H}\in \mbox{res}\, \mathcal{N}_p$  we  can now find    a normal subgroup $\bar{S}$ with  index $p^n$ in $\bar{H}$ which does not intersect $\bar{W}$. Hence $\bar{S}$ is an extension of a free abelian group by a finite $p$-group; let $\bar{V}$ be a free abelian normal subgroup of index $p^m$ in $\bar{S}$.  Its inverse image $V$ has index $p^{n+m}$ in $H$, and $V$ is an extension of the torsion free nilpotent group $M_l(W)$ by the free abelian group $\bar{V}$. Since the index of $V$ is a power of $p$ we can find $k\leq n+m$ such that $M_k(H)\subseteq V$   and   statement v) implies that there exists an index $i_0$ such that  $H_i\subseteq M_k(H)\subseteq V$  if $i\geq i_0$. We pick now an arbitrary $i\geq i_0$ and a  subgroup $Q=H_i\subseteq V$ and $N=H_i\bigcap M_l(W)$. The quotient group $Q/N$ is free abelian because it is a subgroup of the free abelian group $V/M_l(W)$. Lemma $6.1.$ implies that $(Q/N^{\prime} N^p)  \in \mbox{res}\,\mathcal{N}_p$. Since $H/Q$ is a finite $p$-group we obtain that $(H/N^{\prime}N^p)\in \mbox{res}\,\mathcal{N}_p$.  

 This completes the proof of statement vi).   \bigskip

{\bf Proof of statement vii).} Since the group $H/F$ is polycyclic it has  a  unique maximal finite normal subgroup $F_0/F$.  Since $H/F\in res\,\mathcal{N}_p$ we obtain that $F_0/F$ is a $p$-group, and hence $F_0/S$ is an extension of a residually\{finite $p$-group\} $F/S$ by a finite $p$-group $F_0/F$. Hence 
$(F_0/S)\in res\,\mathcal{N}_p$. Since $H/F_0$ contains no non-unit finite subgroups we obtain from statement ii)  that  $H/F_0\in res\,\mathcal{N}_p$. We see that it is enough to prove the statement for  the normal subgroups   $F_0$ and $S$; we can assume  from the very beginning that $F=F_0$ and hence $H/F$ contains no non-trivial finite normal subgroups.

Once again, we obtain from statement ii)  that  the quotient group $H/F$ contains a $p$-series with a unit intersection with  associated restricted Lie algebra of rank equal to the Hirsch rank of $H/F$. Statement vi) now implies  that $H/F$ contains a poly-\{infinite cyclic\} normal  subgroup $H_0/F$ with index $(H:H_0)$  a power of $p$.  Since the index of $H_0$ is a power of $p$ it is enough to prove that $H_0\in res\, \mathcal{N}_p$. We can  assume therefore that $H=H_0$ i.e. the quotient group $H/F$ is poly-\{infinite cyclic\}.

We will now use induction on the Hirsch number of $H/F$. Let $R\supseteq F$ be a normal subgroup of $H$ such that $H/R$ is infinite cyclic group generated by an element $h\in H$. Statement i) implies that the algebra $L_p(R,R_i)$ is abelian with rank $r-1$ which is the Hirsch number of $R$. Since the Hirsch number of $R$ is $r-1$ we can assume that the 
assertion is proven for the group $R$ and  its  normal subgroups $F\supseteq S$, so $R/S\in \mbox{res}\,\mathcal{N}_p$.
 
Let $\bar{X}$ be the image of a subset $X\subseteq H$ under the natural homomorphism $H\longrightarrow H/S$ and  $\bar{h}\not=1$ be an element of $\bar{H}$. Since  $R/S\in res\,\mathcal{N}_p$ we can find a dimension subgroup  $M_n(\bar{R})$ which does not contain $\bar{h}$ and hence the image of $\bar{h}$ in the quotient group $\bar{H}/M_n(\bar{R})$ is an element $h^{\ast}\not=1$.  But  $\bar{H}/M_n(\bar{R})\in res\,\mathcal{N}_p$ by Lemma $6.1.$, hence there exists a homomorphism of this group on a finite $p$-group $G$  which maps $h^{\ast}$ in a non-unit element. We see that the composition of homomorphisms $\bar{H}\longrightarrow \bar{H}/M_n(\bar{R})$ and $ \bar{H}/M_n(\bar{R}\longrightarrow G$ maps $\bar{h}$ on a non-unit element in $G$. This proves that $\bar{H}\in res\,\mathcal{N}_p$ and  completes the proof of statement vii). \bigskip

{\bf Proof of statement iii).} Assume that there exists an element 

\begin{equation} 1\not=\bar{x}\in \bigcap_{i=1}^{\infty}\bar{H_i}\end{equation}

Since $\bar{H}\in res\,\mathcal{N}_p$ we can find  $M_n(\bar{H})$ such that $\bar{x}\not\in M_n(\bar{H}$. If $x$ is an element which is mapped on $\bar{x}$ under the homomorphism $H\longrightarrow \bar{H}$ then $x\not\in M_n(H)U$. Since the topology defined by the series $H_i\,\, (i=1,2,\cdots)$ is equivalent to the $p$-topology we can find $i(n)$ such that $H_{i(n)}\subseteq M_n(H)$ and hence 
$x\not\in H_{i(n)}U$. This implies that $\bar{x}\not\in \bar{H}_{i(n)}$ which  contradicts $(6.8)$\bigskip

{\bf Proof of statement viii)}.  Theorem V implies that there exists a torsion free normal subgroup $F$ with index $(H\colon F)=p^n$ such that the subalgebra $L_p(F,F_i)$ corresponding to the $p$-series $F_i=F\bigcap H_i^{\ast}\,\, (i=1,2,\cdots)$ is finitely generated free abelian  and central in $L_p(H,H^{\ast}_i)$ and its rank is $r_1\leq r$.  Since the topology defined by the series $H_i\,\, (i=1,2,\cdots)$ is equivalent to the $p$-topology by statement v) we  can find a number $i_0$ such that  $H^{\ast}_i\subseteq F$ if $i\geq i_0$. 

On the other hand, the topology defined in $H$ by  series $H_i\,\, (i=1,2,\cdots)$ is equivalent to the $p$-topology by statement v) so there exists $i_1$ such that $H_i\subseteq F\,\, \mbox{if}\,\ i\geq i_1$. Since the algebra $L_p(H,H_i)$ is abelian of finite rank there exists $i_2$ such that the subalgebra $\sum_{i\geq i_2}H_i/H_{i+1}$ is free abelian of rank $r$. We pick now an arbitrary $i$ greater than or equal $k=max\{ i_0,i_1,i_2\}$ and obtain a normal subgroup $U=H_i$   in $H$ with the following properties:

$1)$ $U\subseteq F$; $2)$ The index of $U=H_i$ is a power of $p$; $3)$ The algebra $L_p(U,U_i)$ associated to the $p$-series $U_i=U\bigcap H_i\,\, (i=1,2,\cdots)$ coincides with the subalgebra $\sum_{i\geq i_2}H_i/H_{i+1}$ so it is 
 free abelian of rank $r$; $4)$ The algebra $L_p(U,U_i^{\ast})$ associated to the $p$-series $U_i^{\ast}=U\bigcap F_i^{\ast}=U\bigcap H_i^{\ast}   \,\, (i=1,2,\cdots)$   has a finite index in $L_p(H,H_i)$ because $U$ has a finite index in $H$, so   $L_p(U,U^{\ast}_i)$ is finitely generated; it is free abelian with rank $r_1\leq r$ because it is a subalgebra of finite index in $L_p(F,F^{\ast}_i)$.

 Further since  the index of $U$ in $H$ is a power of $p$ a  straightforward argument show also that the topologies defined in $U$ by the series $U_i\,\, (i=1,2,\cdots)$ and $U_i^{\ast}\,\, (i=1,2,\cdots)$ are equivalent to the $p$-topology of $U$. 

We see that the relation $r_1=r$ will follow from the following lemma.\bigskip

{\bf Lemma 6.2.}   {\it Let $U$ be an arbitrary  group,

\begin{equation} U=U_1\supseteq U_2\cdots\end{equation}

\begin{equation} U=U^{\ast}_1\supseteq U^{\ast}_2\supseteq\cdots\end{equation}

be two $p$-series with unit intersection. Assume that the topologies defined by these series are equivalent and that the algebras $L_p(U,U_i)$ and $L_p(U,U_i^{\ast}$ are free abelian;  then they are isomorphic.} \bigskip

{\bf Proof.} The algebra $gr_{\rho}({Z_p(U,U_i)})$ is isomorphic to the polynomial ring $Z_p T$ where $T$ is a free system of generators for $L_p(U,U_i)$.  We consider now the completion $\widetilde{Z_pU}$ in the $\rho$-topology. A routine argument (see, for instance, the proof of statement ii) of Lemma $3.6.$ or the end of the proof  of Lemma $3.7.$)  shows that every elenment   $x\in \widetilde{Z_pU}$ has a unique representation 

$$ x=\sum_{i=0}^{\infty}\lambda_i\pi_i$$

where $\lambda_i\in Z_p\,\, (i=0,1,\cdots)$, $\pi_i\,\, (i=0,1,\cdots)$ are standard monomials on the set of variables $T$.  The algebra $\widetilde{Z_pU}$ has an augmentation ideal $\omega(\widetilde{Z_pU})$ which is topologically generated by  the of elements  $ T$. We consider the filtration defined by the powers of this ideal and obtain that the graded ring of $\widetilde{Z_pU}$ associated to this filtration is isomorphic to the polynomial ring $Z_p[T]$.

Since the completion of $Z_pU$ with respect to the topology defined by the valuation $\rho^{\ast}$ is isomorphic to $\widetilde{Z_pU}$ (see Lemma $2.6.$)  we obtain that
$gr_{\rho^{\ast}}(Z_pU)\cong Z_p[T]$. On the other hand, 
$gr_{\rho^{\ast}}(Z_pU)\cong Z_p[T_1]$ where $T_1$ is a free system  of generators $L_p(U,U_i^{\ast})$, so  we obtain that the cardinalities of $T$ and $T_1$ are equal and $L_p(U,U_i)\cong L_p(U,U_i^{\ast})$.

This completes the proof of the lemma and the proof of Theprem VI.\bigskip

{\bf 6.3.} We will need in the proof of theorem XII the following theorem which is  a refined version of statement vi) of Theorem VI.  \bigskip

{\bf Theorem 6.1.}  {\it  Let $H$ be a polycyclic group with Hirsch number $r$. Assume that there exists a $p$-series $(6.1)$ with unit intersection such that $L_p(H,H_i)$ is  abelian of rank $r$. Let $i_0$ be the index defined in statement}  vi) { \it of Theorem} VI. { \it Then there exists an index 
$i_1\geq i_0$ such that if $i\geq i_1$ then  the subgroups $Q=H_i$ and $N$  obtained  in Theorem VI  have the following series }

\begin{eqnarray} Q=Q^{(1)}\supseteq Q^{(2)}\supseteq\cdots\supseteq Q^{(k-1)}\supseteq Q^{(k)}=\nonumber\\
=N\supseteq Q^{(k+1)}\supseteq\cdots\supseteq Q^{(r-1)}\supseteq Q^{(r)}\supseteq Q^{(r+1)}=1\end{eqnarray}

{\it with following properties:}

 i) {\it Every factor $Q^{(j)}/Q^{(j+1)}\,\, (j=1,2,\cdots,r-1)$ is an infinite cyclic group.}

ii) {\it   Let  $L_p(Q^{(j)},Q^{(j)}_i)$ be the restricted Lie algebra of $Q^{(j)}$ associated to the $p$-series $Q^{(j)}_i=Q^{(j)}\bigcap H_i\,\, ( i=1,2,\cdots)$. Then the  algebra 

$L_p(Q^{(j)},Q^{(j)}_i)/L_p(Q^{(j+1)},Q^{(j+1)}_i)$ is free abelian of rank $1$.
If $q_j$ is an element of $Q^{(j)}$ which generates the quotient group $Q^{(j)}/Q^{(j+1)}$ then its  homogeneous component $\tilde{q}_j$  in $L_p(Q^{(j)},
Q^{(j)}_i)\subseteq L_p(H,H_i)$ generate the quotient algebra  $L_p(Q^{j)}_i)/ L_p(Q^{j+1)}_i)$. Hence the system of elements $q_1,q_2,\cdots,q_r$ generate the subgroup $Q$ and the system of homogeneous components  
$\tilde{q}_1,\tilde{q}_2, \cdots,\tilde{q}_r$  freely generate the free abelian algebra $L_p(Q,Q_i)$.}\bigskip

We need first two lemmas.\bigskip 

{\bf Lemma 6.3.} {\it Let $H$ be a group with a $p$-series $(6.1)$. Assume that there  exists a normal subgroup $Q$ with $(H\colon Q)=p^n$ and an index $i_0$ such that $H_i\subseteq Q\,\, \mbox{if}\,\,i\geq i_0  $ and let }
\begin{equation} H=S_1\supseteq S_2\supseteq\cdots\supseteq \S_n\supseteq S_{n+1}=Q\end{equation}

{\it be an invariant series with factors $S_j/S_{j+1}\,\, (j=1,2,\cdots,n)$ cyclic groups of order $p$; let $s_j$ be an element of $S_j$ which generates $S_j/S_{j+1}\,\, (j=1,2,\cdots,n)$. Let  $S_{j,i}=H_i\bigcap S_j\,\, (i=1,2,\cdots)$ for every given $1\leq j\leq n$. 
Then every quotient algebra $L_p(S_j, S_{j,i})/L_p(S_{j+1},S_{j+1,i})$ has dimension $1$ and is generated by the homogeneous component $\tilde{s}_j\,\, (j=1,2,\cdots,r-1)$;  
the system of these homogeneous components together with the ideal $L_p(Q,Q_i)$ generates the algebra $L_p(H,H_i)$}.\bigskip 

{\bf Proof.} We denote  $U=S_2$ and consider the series 

\begin{equation} U=S_2\supseteq S_3\supseteq\cdots\supseteq S_n\supseteq S_{n+1}=Q\end{equation}

 Let $\bar{H}=H/U$ and $\bar{H}_i=(H_iU)/U\,\, (i=1,2,\cdots)$ .  Since $\bar{H}_{i_0}=1$ we see that $\bigcap_{i=1}^{\infty}\bar{H}_i=1$ so the algebra   $L_p(\bar{H},\bar{H}_i)$ is non-zero. On the other hand, $\bar{H}$ is a cyclic group of order $p$ so the dimension of $L_p(\bar{H},\bar{H}_i)$ can not be greater than $1$. We obtain therefore that $dim(L_p(\bar{H},\bar{H}_i)=1$, and that $L_p(\bar{H},\bar{H}_i)$ is generated by the homogeneous component $\tilde{s}_1$. 
We see that the proof is now reduced to the subgroup $U=S_2$  and series 
$(6.14)$. Since $(U\colon Q)=p^{n-1}$ we obtain  after $n$ steps the system of elements  $\tilde{s}_1,\tilde{s}_2,\cdots, \tilde{s}_n$ which together with the subalgebra $L_p(Q,Q_i)$ generates $L_p(H,H_i)$.\bigskip

{\bf Corollary 6.2.} {\it The dimension of the algebra $L_p(H,H_i)/L_p(Q,Q_i)$ is $n$}.\bigskip

{\bf Corollary 6.3.} {\it The system of elements $\tilde{s_j}$ together with the ideal $L_p(Q,Q_i)$ generates the algebra $L_p(H,H_i)$.}\bigskip

We will need in the proof of Theorem XII one more corollary of Lemma $6.3.$  \bigskip

{\bf Corollary 6.4.}  {\it Let $V$ be a polycyclic group with Hirsch number $r$, $W$ be a finite subgroup of index $p$, $s$ be an element of $V\backslash W$, and let $s^p=a\in W$. Assume that there exists in $V$ a $p$-series $V_i\,\, (i=1,2,\cdots)$ with unit intersection such that the algebra $L_p(V,V_i)$ is free abelian of rank $r$. Let $\rho$ be the valuation in $KV$ defined by this series, $\widetilde{(s-1)}, \widetilde{(a-1)}$ be the homogeneous components of the element $s-1, a-1$ in the rings $gr(KV)$ and $gr(KW)$ respectively.

Then the polynomial ring $gr(KV)$ is a simple algebraic extension of the polynomial ring $KW$

\begin{equation} gr(KV)\cong gr(KW)[\widetilde{(s-1)}]\end{equation}

where the minimal polynomial of $\widetilde{(s-1)}$ is $t^p-\widetilde{(a-1)}$.}\bigskip

{\bf Proof.} Lemma $6.3.$ implies that the free abelian algebra $L_p(V,V_i)$ is an extension of the free abelian algebra $L_p(W,W_i)$ by the  one dimensional algebra generated by the element $\tilde{s}$, which is the homogeneous component of $s$, and $\tilde{s}^{[p]}=\tilde{a}$. We obtain from this that the ring $U_p(L_p(V,V_i))$ is a simple algebraic extension $U_p(L_p(V,V_i))\cong U_p(L_p(W,W_i))[\tilde{s}]$. The assertion now follows from the isomorphism  $U_p(L_p(V,V_i))\cong gr(KV)$ obtained in Proposition $2.7$.\bigskip

{\bf Lemma 6.4.} {\it Let $H$ be a free abelian group  of rank $r$. Assume that there exists a $p$-series }

\begin{equation} H=H_1\supseteq H_2\supseteq\cdots \end{equation} 

{\it with unit intersection such that the algebra $L_p(H,H_i)$ is free abelian of rank $r$. Then there exists a free system of generators $h_1,h_2,\cdots, h_r$ such that the homogeneous components $\tilde{h}_1,\tilde{h}_2,\cdots, \tilde{h}_r$ form a free system of generators for the free abelian algebra $L_p(H,H_i)$.} \bigskip

{\bf Proof.} Let $u_1,u_2,\cdots,u_r$ be a free system of generators for $L_p(H,H_i)$. Since the algebra $L_p(H,H_i)$ is graded we can assume that the elements $u_i\,\, (i=1,2,\cdots, r)$ are homogeneous, so every  element $u_i$ is the homogeneous component of an  element $v_i\in H$. We pick now a system of elements $v_i\in H$ such that $\tilde{v}_i=u_i\,\, (i=1,2,\cdots,n)$; let $V$ be the subgroup generated by these elements. 

The subgroup $V$ contains a $p$-series $V_i=V\bigcap H_i\,\, (i=1,2,\cdots)$ such that the Lie algebra $L_p(V,V_i)$ coincides with the algebra $L_p(H,H_i)$.  We claim that this implies that the index $(H\colon V)$ is prime to $p$. In fact, if we assume that $p|(H\colon V)$ then we can find a subgroup $U\supseteq V$ with index $(H\colon U)=p$ and $L_p(U,U_i)=L_p(H,H_i)$,  but  Corollary $6.2.$ implies that the dimension of the algebra $L_p(H,H_i)/L_p(U,U_i)$ is $1$. We see that $p$ is not a divisor of $(H\colon U)=p$.

We obtain therefore that there exists in $H$ a free system of generators $h_i\,\, (i=1,2,\cdots)$ such that $v_i=h_i^{m_i}\,\, (i=1,2,\cdots,r)$ where $m_1,m_2,\cdots, m_r$ are integers prime to $p$.  The  system of elements 
$h_1,h_2,\cdots, h_r$ generates $H$ and the homogeneous components of these elements generate $L_p(H,H_i)$. This system of generators satisfies the conclusion of the assertion. \bigskip 

{\bf Proof of Theorem 6.1.} Let   $i_0$ be as in statement vi) of Theorem VI, $H_{i_0}=Q_0$, $N_0\subseteq Q_0$  be the torsion free nilpotent normal subgroup with $Q_0/N_0$ free abelian. Let $Z_0$ be the center of $N_0$, $l$ be its rank,  and $\bar{H}=H/Z_0$. Statement ii) of Theorem VI implies  and the algebra $L_p(\bar{H},\bar{H}_i)$ is  abelian of rank $r-l$. Since $\bar{H}_{i_0}=H_{i_0}/Z_0$ is torsion free we obtain once again from statement ii)  that  $\bigcap_{i=1}^{\infty}\bar{H}_i=1$.  We can now apply induction by the Hirsch number of $H$ and to assume   that there exists an index $i_1\geq i_0$ such that if $i\geq i_1$ is given and $Z=H_i\bigcap Z_0$ then  the  subgroup $\bar{Q}=\bar{H}_i=H_i/Z$  of 
$\bar{H}$ contains the following  series with infinite cyclic factors which satisfies the conclusions of the assertion

\begin{eqnarray} \bar{Q}=\bar{Q}^{(1)}\supseteq \bar{Q}^{(2)}\supseteq\cdots\supseteq \bar{Q}^{(k-1)}\supseteq \bar{Q}^{(k)}=
\bar{N}\supseteq \nonumber\\
\supseteq \bar{Q}^{(k+1)} \supseteq\cdots\supseteq \bar{Q}^{(n-l)}\supseteq Q^{(n-l+1)}=1\end{eqnarray}

We obtain from this the following series in $H$

\begin{eqnarray} Q=Q^{(1)}\supseteq Q^{(2)}\supseteq\cdots\supseteq Q^{(k-1)}\supseteq Q^{(k)}=N\supseteq\nonumber\\ \supseteq Q^{(k+1)})\supseteq\cdots\supseteq Q^{(n-l)}\supseteq Q^{(n-l+1)}=Z\end{eqnarray}

where  every   factor
 $ Q^{(j)}/Q^{(j+1)}$ is an  infinite cyclic groups generated by an element $q_j$, and every  quotient algebra $L_p(Q^{(j)})/L_p(Q^{(j+1)})$ is  free abelian of rank $1$  generated by the element $\tilde{q}_j$. We apply now Lemma $6.4.$ to obtain a free system of generators whose homogeneous components freely generate the subalgebra $L_p(Z,Z_i)$   and the assertion follows.\bigskip

 Let $H$ be a polycyclic group with Hirsch number $r$. Assume that there exists a $p$-series $(6.1)$ with unit intersection such that $L_p(H,H_i)$ is  abelian of rank $r$. We will need in the proof of Theorem XII some special system of generators for the algebra $L_p(H,H_i)$ which is obtained in the following way. 
We obtain from Theorem VI and  Theorem  $6.1.$   that there exist  a poly-\{infinite cyclic\} normal subgroup $Q$  with quotient group $G=H/Q$ of index  $p^n$,  and   a system of generators 
$q_1,q_2,\cdots,q_r$  such that the free abelian algebra $L_p(Q,Q_i)$ is freely generated by the system of elements  $\tilde{q}_1,\tilde{q}_2,\cdots,\tilde{q}_r$. Lemma $6.3.$ implies that there exists  a system of elements $s_1,s_2,\cdots,s_n$ which generates the quotient group   $H/Q$ and the images of the elements $\tilde{s}_1,\tilde{s}_2,\cdots,\tilde{s}_n$ in $L_p(H,H_i)$ generate the quotient algebra $L_p(H,H_i)/L_p(Q,Q_i)\cong L_p(G,G_i)$. We have the following corollary of the results of this subsection.\bigskip

{\bf Corollary 6.5. and Definition 6.1.} {\it Assume that the conditions of Theorem $6.1.$ hold. Then the  system of elements  
$$q_1,q_2,\cdots,q_r; s_1,s_2,\cdots,s_n$$ 
generates the group $H$; the 
system of their homogeneous components   $$\tilde{q}_1,\tilde{q}_2,\cdots,\tilde{q}_r; \tilde{s}_1,\tilde{s}_2,\cdots,\tilde{s}_n$$ 
generates the algebra $L_p(H,H_i)$. This system of generators for $H$ is called a special system of generators.}

\setcounter{section}{7}
\setcounter{equation}{0}

\section*{\center \S  7}

{\bf 7.1.} {\bf Lemma 7.1.} {\it Let $H^{(j)}\,\, (j\in J)$ be a system of groups.  Assume that every group $H^{(j)}$ contains  invariant subgroups $S^{(j)}\supseteq U^{(j)}$, and $U^{(j)}$ contains  a p-series $U^{(j)}_i\,\, (i=1,2,\cdots)$ with unit intersection whose terms  are invariant in $S^{(j)}$. 
Let} $H=\prod_{j\in J}H^{(j)},
S=\prod_{j\in J}S^{(j)}, U=\prod_{j\in J}U^{(j)}$.

{\it There exists in $U$ a $p$-series $U_i\,\, (i=1,2,\cdots)$ such that all the subgroups $U_i$ are $S$-invariant, $U_i\bigcap U^{(j)}=U^{(j)}_i\,\, (i=1,2,\cdots)$ and the associated graded algebra $L_p(U,U_i)$ is isomorphic to the direct sum of the Lie algebras} $ L_p(U^{(j)},U^{(j)}_i)$

\begin{equation} L_p(U,U_i)\cong\bigoplus_{j\in J} L_p(U^{(j)},U^{(j)}_i)\end{equation}

{\bf Proof.}  The $p$-series in $U^{(j)}$ defines a weight function 
$f^{(j)}$ in $U^{(j)}$, and hence a weight function $f$ on the set $\bigcup_{j\in J}U^{(j)}$. We define now a weight function $f(x)$  on the group $U$. Let 

\begin{equation}x=x_{j_1}x_{j_2}\cdots x_{j_n}\end{equation}

 be an element of $U$. Then it weight in $U$ is defined as 

\begin{equation} f(x)=
\mbox{min}\{f(x_{j_1}), f(x_{j_2}),\cdots,f(x_{j_n})\}\end{equation}

A straightforward verification shows  that $f([x,y])\geq f(x)+f(y)$, $f(x)^p\geq p(f(x))$,\,$ f(xy^{-1})\geq f(x)-f(y)$; hence the  obtained  weight function $f$ on $U$ defines a $p$-series $U_i\,\, (i=1,2,\cdots)$ in $ U$ and we obtain a Lie algebra $L_p(U,U_i)$ corresponding to this $p$-series.  The definition of $f(x)$ implies that  its restriction on  every $U^{(j)}$ coincides with $f^{(j)}$, hence we obtain a natural imbedding of the Lie algebra $L_p(U^{(j)}, U^{(j)}_i)$ into 
$L_p(U,U_i)$. It is clear that for every $s\in S$ $f(s^{-1}us)=f(u)$ so every subgroup $U_i$ is $S$-invariant. 

We consider now the homogeneous component of the element $(7.2)$. Assume that $f(x)=k$. We can assume that   $f(x_{j_1})=f(x_{j_2})=\cdots=f(x_{j_l})=k$ and the weights of the remaining components is greater than $k$.   Hence, the homogeneous component $\tilde{x}$ of $x$ coincides with the homogeneous component of the element $x_1x_2 \cdots x_l$. We   apply 
 Lemma $2.2.$ and  obtain     that the homogeneous component of the  element $x $ is equal 

\begin{equation}\tilde{x}=\tilde{x}_{j_1}+\tilde{x}_{j_2}+\cdots+
\tilde{x}_{j_l}\end{equation}

where $\tilde{x}_{j_{\alpha}}\in L_p(U^{(j_{\alpha})}, U^{(j_{\alpha})}_i)\,\, (\alpha=1,2,\cdots,l)$.

 Hence the algebra $L_p(U,U_i)$ is generated by the system of subalgebras 
$L_p(U^{(j)}, U^{(j)}_i)\,\, (j\in J)$; we see that in order to prove relation $(7.1)$ it is enough to verify that   representation $(7.4)$ is unique. 

Assume that there exists  another representation. A routine argument shows that we can assume now that the set $J$ is finite and that after a suitable numeration of it this second representation has a form 

\begin{equation} \tilde{x}=\tilde{x}_1+\tilde{x}_2+\cdots\tilde{x}_m\end{equation}
 with $\tilde{x}_j\in L_p(U^{(j)}, U^{(j)}_i)\,\, (j=1,2,\cdots, m)$. Since all the elements in  the right side of  $(7.5)$ are homogeneous we obtain that their degrees and the degrees  must be equal to $k$. 
We obtain now from $(7.4)$ and $(7.5)$  the following two representation for $x$ 

\begin{equation} x=x_{j_1}x_{j_2}\cdots x_{j_l}y_1\end{equation}
and 

\begin{equation} x=x_1x_2\cdots x_m y_2\end{equation}

where $f(y_{\beta})>k\, \, (\beta=1,2)$  or, equivalently,   $y_{\beta}\in \prod_{j\in J}U^{(j)}_{k+1}\,\, (\beta=1,2)$. 
We compare now representations $(7.6)$ and $ (7.7)$ for $x$ and obtain that they coincide modulo the subgroup $\prod_{j\in J}U^{(j)}_{k+1}$; hence $l=m$ and after a renumeration we obtain that $ x_{j_{\alpha}}\in (x_j U^{(j) }_{k+1})\,\, (j=1,2,\cdots,l)$. Hence the homogeneous components of these elements in $L_p(U,U_i)$ coincide which proves the uniqueness of representation $(7.4)$. This completes the proof.\bigskip

{\bf Lemma 7.2.} {\it Let $W=S\wr G$ be the discrete wreath product of a group $S$ by a group $G$,  $S^{\ast}$ be the base group of $W$. 
 Assume that $S$ contains a normal subgroup $U$ with an $S$-invariant $p$-series 
\begin{equation}U=U_1\supseteq U_2\supseteq\cdots\end{equation}

 with unit intersection   and let 
$U_i^{\ast}=\prod_{g\in G} g^{-1}U_ig$.
  Then there  exists in $U^{\ast}$ a $W$-invariant $p$-series with unit intersection 
\begin{equation} U^{\ast}=U^{\ast}_1\supseteq U^{\ast}_2\supseteq\cdots\end{equation}  

  with  Lie agebra $L_p(U^{\ast},U^{\ast}_i)$  isomorphic to the direct sum of  copies of the algebra  $L_p(U,U_i)$ }

\begin{equation}   L_p(U^{\ast},U^{\ast}_i)\cong \bigoplus_{g\in G} L_p(U,U_i)\end{equation}

{\bf Proof.}  We apply Lemma $7.1.$ and obtain a  weight function $f$  on $U^{\ast}$ and  a $p$-series in $U^{\ast}$; the terms of this  $p$-series are invariant subgroups in $S^{\ast}$. We will prove now that they  are $G$-invariant and hence $W$-invariant.

Indeed, let $u$ be a non-unit element of $U^{\ast}$,  $u=u_{j_1}u_{j_2}\cdots u_{j_k}$ with $ u_{j_{\alpha}}\in U_{j_{\alpha}}=g_{\alpha}^{-1}Ug_{\alpha} \,\, (\alpha=1,2,\cdots,k)$. We can assume that  $f(u)=f(u_{j_1})=f(u_{j_2})=\cdots f(u_{j_l})=l$  and the rest of the factors, $u_{j_{l+1}},u_{j_{l+2}},\cdots, u_{j_k}$ have weights greater than $l$.   

 Let $g$ be an arbitrary element of $G$. Then $g^{-1}ug=
(g^{-1}u_{j_1}g)(g^{-1}u_{j_2}g)\cdots(g^{-1}u_{j_k}g)$ where 
$g^{-1}u_{j_{\alpha}}g\in g^{-1}(g_{\alpha}^{-1}Ug_{\alpha})g \,\, (\alpha=1,2,\cdots,k)$ and we obtain from this that $f (g^{-1}ug)=f(u)$; we see that the weight function in $U^{\ast}$ is $G$-invariant and so is the $p$-series defined by it. 

The rest of the statements follow immediately.\bigskip
 
We apply now this lemma to the case when $U=S$ and the group $G$ is finite. Lemma $7.1.$ yields a $W$-invariant $p$-series $S^{\ast}_i\,\, (i=1,2,\cdots)$  in the base group $S^{\ast}$.  We will need the following fact about this $p$-series.\bigskip

{\bf Corollary 7.1.} {\it Assume that $S$ contains a $p$-series with unit intersection and that the topology defined by this series is equivalent to the $p$-topology, and that the group $G$ is finite. Then the topology defined by the series $S_i^{\ast}\,\, (i=1,2,\cdots)$ is equivalent to the $p$-topology in $S^{\ast}$.}\bigskip

{\bf Proof.} For a given $n$ we can find a number $i(n)$ such that $S_{i(n)}\subseteq M_n(S)$ which implies that $S^{\ast}_{i(n)}\subseteq M_n(S^{\ast})$ and the assertion follows.\bigskip

 {\bf Proposition 7.1.} {\it Let $H$ be a group, $U$ be a normal subgroup which contains a $p$-series $(7.8)$ with unit intersection. Assume that there exists a normal subgroup $S\supseteq U$ of finite index $n$  in $H$ such that $h^{-1}U_ih=U_i\,\, (i=1,2,\cdots)$ for all $h\in S$. 

 Then there exists an $H$-invariant $p$-series 

\begin{equation}U=V_1\supseteq V_2\supseteq\cdots\end{equation}

 with unit intersection 
 such that the algebra $L_p(U,V_i)$ is isomorphic to a subalgebra of the direct sum of $n$ copies of the algebra $L_p(U,U_i)$.}\bigskip

{\bf Proof.} Let $G=H/S$.  We consider the usual imbedding of $H$ into the wreath product $W=S\wr G$  where $S$ is isomorphically imbedded into the base group
$S^{\ast}$  of $W$; this base group is the direct product of $(G:1)$ copies of  $S$. 

We apply now Lemma $7.2.$ to the  subgroup $U^{\ast}=\prod_{g\in G} g^{-1}Ug$ and obtain    that there exists a   $W$-invariant series $(7.9)$.

The series $V_i=U^{\ast}_i\bigcap U\,\, (i=1,2,\cdots)$ has unit intersection and the subalgebra $L_p(U,V_i)$ associated to this series is isomorphic to a subagebra of $L_p(U^{\ast},U^{\ast}_i) $ which in its turn is isomorphic to the sum of $n$ copies of $ L_p(U,U_i)$. This completes the proof. \bigskip

{\bf Corollary 7.2.} i) {\it If the algebra $L_p(U,U_i)$ is free abelian, (abelian), (abelian of finite rank)   then so is $L_p(U,V_i$.}

ii) {\it If $H$ is a polycyclic group with Hirsch number $r$ and $L_p(H,H_i)$ is free abelian,  (abelian of rank  $r$)  then so is $L_p(U,U^{\ast}_i)$.}\bigskip

{\bf Proof.} We will prove statement ii); the proof of statement i) is obtained by obvious simplification of the argument. The group $U$ is a subgroup of the group 
 $U^{\ast}=\underbrace{ U\times U\cdots\times U}_{n }$. The algebra $L_p(U^{\ast},U_i^{\ast})$ has rank  $rn$ because it is a direct sum of $n$ copies of $L_p(U,U_i)$. Since the rank of $L_p(U^{\ast},U_i^{\ast})$ coincides with the Hirsch number of $U^{\ast}$  we obtain from  statement i) of Theorem VI  that the rank of the algebra $L_p(U,V_i)$ must be  equal to the Hirsch number of $U$ which is equal $r$. 

Now assume that $L_p(U,U_i)$ is free abelian. Then the direct sum of $n$ copies of it is free abelian. Since $L_p(U,V_i)$ is a subalgebra of this direct sum it must be free abelian.

This completes the proof. \bigskip

{\bf 7.2.}  {\bf Proposition 7.2.} {\it Let $H$ be a group which contains a $p$-series} 
\begin{equation} H=H_1\supseteq H_2\supseteq\cdots   \end{equation}

{\it with unit intersection. Assume that the topology defined by this series is equivalent to the  $p$-topology in $H$. Let
$\Phi$ be a group of automorphisms of $H$.}

  i)  {\it  Assume that the restricted Lie algebra $L_p(H,H_i)$ is generated by the first $l$ factors $H_i/H_{i+1}\,\, (i=1,2,\cdots,l)$ where  the subgroups $H_i\,\, (i=1,2,\cdots,l)$ are $\Phi$-invariant. Then all the subgroups $H_i\,\, (i=1,2,\cdots )$ are $\Phi$-invariant.}

ii)   {\it If}  

\begin{equation} [\Phi,H_i]\subseteq H_{i+1}\,\, (i=1,2,\cdots,l)\end{equation}

 {\it i.e. $\Phi$ centralizes the first $l$ factors   $H_i/H_{i+1}\,\,$ then it centralizes all the factors $H_i/H_{i+1}$. }\bigskip 

{\bf Proof.}   We will prove first   statement i)  for the case  when series $(7.12)$ has a finite length, say $k$. Hence $H_k=1$ but  $H_{k-1}\not=1$. We  begin by proving   that $H_{k-1}$ is $\Phi$-invariant.  Let $h\in H_{k-1}$. Lemma $2. 4.$ implies that   $\tilde{h}$ can be  expressed in $L_p(H,H_i)$ as  a sum of the Lie monomials 

\begin{equation}[\tilde{h}_{{\alpha}_1},\tilde{h}_{{\alpha}_2},\cdots,\tilde{h}_{{\alpha}_s}]^{[p]^{n_{\alpha}}}\end{equation} 

where the homogeneous elements  $\tilde{h}_{\alpha_1},
\tilde{h}_{\alpha_2}, \cdots,\tilde{h}_{\alpha_s} $ are taken from the first $l$ factors $H_i/H_{i+1}$ and the weight of every such monomial is $k-1$.    Since $H_k=1$ we see  that  the element  

\begin{equation}[ h_{{\alpha_1}},h_{{\alpha_2}},
\cdots, h_{{\alpha_s}}]^{p^{n_{\alpha}}}\end{equation}

of the group $H$ is a coset representative for the Lie monomial $(7.15)$, and hence $h$ is a product of the elements $(7.14)$ with  weight   $k-1$. 

Now let $\phi$ be an arbitrary automorphism from $\phi$. We obtain that the image of the element $(7.15)$ is

\begin{equation} \phi([h_{{\alpha}_1},h_{{\alpha}_2},\cdots, h_{{\alpha}_s}]^{p^{n_{\alpha}}})= [\phi(h_{{\alpha}_1}),\phi(h_{{\alpha}_2}),\cdots, \phi(h_{{\alpha}_s})]^{p^{n_{\alpha}}}\end{equation}
  
 Let $w(x)$ denote the weight of an element $x\in H$. Since element $(7.15)$ is the representative of the monomial $(7.14)$ we obtain from Lemma $2.4.$ that its weight is 

\begin{equation}k-1=p^{n_{\alpha}}\sum_{i=1}^s  w(h_{\alpha_i})\end{equation}

Since $\phi(H_i)=H_i\,\, (i=1,2,\cdots,l)$ we obtain that 
$w(h_{\alpha_i})=w(\phi(h_{\alpha_i}))\,\, (i=1,2,\cdots,s)$
and then 

$$w([\phi(h_{{\alpha}_1}),\phi(h_{{\alpha}_2}),\cdots, \phi(h_{{\alpha}_s})]^{p^{n_{\alpha}}})\geq p^{n_{\alpha}}\sum_{i=1}^s  w(h_{\alpha_i})=k-1$$

This together with $(7.16)$  implies that $w(\phi([h_{{\alpha}_1},h_{{\alpha}_2},\cdots, h_{{\alpha}_s}]^{p^{n_{\alpha}}}))\geq k-1$, which means that $\phi([h_{\alpha_1},
h_{\alpha_2},\cdots, h_{\alpha_s}]^{p^{n_{\alpha}}})\in H_{k-1}$. Hence $\phi(h)$ is a product of elements from $H_{k-1}$ and we obtain that  $\phi(h)\in H_{k-1}$.
Since $h$ was an arbitrary element of $H_{k-1}$ we obtain that $H_{k-1}$ is $\Phi$-invariant.

  The group  $\bar{H}=H/H_{k-1}$ has a $p$-series   $\bar{H}_i=H_i/H_{k-1}\,\, (i=1,2,\cdots,k-1)$ whose length is shorter than the length of the series
 $H_i\,\, (i=1,2,\cdots,k)$ and the Lie algebra $L_p(\bar{H},\bar{H}_i)$ is  generated by the first $l$ factors $\bar{H}_i/\bar{H}_{i+1}\,\, (i=1,2,\cdots,l)$ because $L_p(\bar{H},\bar{H}_i)$ is a homomorphic image of $L_p(H,H_i)$. 
  The action of $\Phi$ on $H$ defines in a natural way the action of $\Phi$ on  $\bar{H}$ and we  assume by induction that all the subgroups $\bar{H}_i\,\, (i=1,2,\cdots,k-1)$ are $\Phi$-invariant. This implies that their inverse images $H_i\,\, (i=1,2,\cdots,k)$ are $\Phi$-invariant; this completes the proof of statement i) for the special case when the series $(7.12)$ has a finite length.

We consider now the general case. To prove that the  term 
$H_i\,\, (i>l)$ of series $(7.12)$ is $\Phi$-invariant we pick $n>i$ and consider the quotient group  
 $G=H/M_n(H)$; let $G_i$ be the image of $H_i$ under the natural homomorphism $H\longrightarrow G$. Proposition $2.1.$ implies that the epimorphism $H\longrightarrow G$ defines an epimorphism $L_p(H,H_i)\longrightarrow L_p(G,G_i)$. Since the algebra $L_p(H,H_i)$ is generated by its first $l$ factors we conclude that the algebra $L_p(G,G_i)$ is generated by the first $l$ factors $G_i/G_{i+1}\,\, (i=1,2,\cdots,l)$.    There exists $m\geq i$ such that $H_m\subseteq M_n(H)$ and hence $G_m=1$. We see that the   series $G_i\,\, (i=1,2,\cdots)$ has finite length in $G$ and the algebra $L_p(G,G_i)$ is generated by the factors $G_i/G_{i+1}\cong ({H_i}/M_n(H))/(H_{i+1}/M_n(H))\cong H_i/H_{i+1}\,\, (i=1,2,\cdots,l)$. The group of automorphisms $\Phi$ acts in a natural way on  $G$, we denote this group by $\Psi$, and obtain from the proven special case     that  an arbitrary  subgroup  $G_i$ is  $\Psi$-invariant; hence its inverse image $H_i$ is $\Phi$-invariant. 

The proof of  statement i) is complete.
  
ii) We prove now the second statement.   Since all the subgroups $H_i\,\, (i=1,2,\cdots)$ are $\Phi$-invariant the action of the group $\Phi$ on $H$ defines in a natural way its action   on the graded Lie algebra $L_p(H,H_i)$  by Corollary $2.3.$; since $\Phi$  centralizes the factors $H_i/H_{i+1}\,\, (i=1,2,\cdots,l)$ which generate the algebra $L_p(H,H_i)$  it centralizes all the factors and the proof is complete. \bigskip

{\bf Proposition 7.3.} {\it  Let $H$ be a finitely generated group,  $\Phi$ be  a group of automorphisms of $H$. Assume that $H$  contains a  $p$-series $(7.12)$  with unit intersection such that the topology defined by  this $p$-series is equivalent to  the  $p$-topology and the Lie algebra $L_p(H,H_i)$ is finitely generated, or equivalently, is generated by a finite number of factors $H_i/H_{i+1}\,\, (i=1,2,\cdots, l)$. Then }

i) {\it  There exists a normal subgroup  $\Phi_1$  of finite index in $\Phi$ such that   
  all the subgroups $H_i\,\, (i=1,2,\cdots,)$ are $\Phi_1$-invariant and the factors  $H_i/H_{i+1}\,\, (i=1,2,\cdots)$ are cenralized by $\Phi_1$}.

ii) {\it If the order of every automorphism $\phi\in \Phi$ in the quotient group $H/H^{\prime}H^p$ is a power of $p$ then the index of $\Phi_1$ is also a power of $p$.}\bigskip

{\bf Proof.} We pick an arbitrary $m>l$  and consider the quotient group $\bar{H}=H/M_m(H)$.  The action of the group $\Phi$ on $H$ defines in a natural way a group $\bar{\Phi}$ of automorphisms of $\bar{H}$ and an epimorphism 
$\Phi\longrightarrow \bar{\Phi}$.  

Since $\bar{\Phi}$ is finite   we obtain a normal subgroup $\Phi_1\subseteq \Phi$ which acts on $\bar{H}$ trivially,  hence $\Phi_1$ acts trivially on every subgroup $\bar{H}_i=H_i/M_m(H)\,\, (i=1,2,\cdots,m)$.

Since the subgroups $\bar{H}_i\,\, (i=1,2,\cdots,m)$ are $\Phi_1$-invariant subgroups of $\bar{H}$   their inverse images $H_i\,\, (i=1,2,\cdots,m)$ are $\Phi_1$-invariant in $H$. Proposition 
$7.2$ now implies that all the subgroups $H_i\,\, (i=1,2,\cdots)$ are $\Phi_1$-invariant. Finally, $\Phi_1$ acts trivially on every factor $\bar{H}_i/\bar{H}_{i+1}\cong H_i/H_{i+1}\,\, (i=1,2,\cdots,l)$ so statement ii) of Proposition $7.2$  implies that $\Phi_1$ centralizes all the factors 
$H_i/H_{i+1}$. This completes the proof of statement i).

 We prove now statement ii). The subgroup $\Phi_1$ is the kernel of the map $\Phi\longrightarrow \bar{\Phi}$; if    the order of every 
$\phi\in \Phi$ in the quotient group $\bar{H}/\bar{H}^{\prime}\bar{H}\cong H/H^{\prime}H^p$ is a power of $p$ then $\bar{\Phi}$ is a $p$-group by Burnside's Theorem, so the index of $\Phi_1$       is a power of $p$ if the conditions of statement ii) hold, and statement ii) follows. \bigskip

{\bf Proposition 7.4.} {\it  Let $H$ be a finitely generated group. Assume that it contains a $p$-series $(7.12) $  with unit intersection such that the topology defined by  this $p$-series is equivalent to  the $p$-topology and the Lie algebra $L_p(H,H_i)$ is finitely generated.

There exists a $p$-series }

\begin{equation}H=V_1\supseteq V_2\supseteq\cdots  \end{equation} 

{\it of characteristic subgroups with unit intersection whose associated graded algebra  $L_p(V,V_i)$  is isomorphically imbeded into the direct sum of a finite number of isomorphic copies of the algebra  $L_p(H,H_i)$.}\bigskip

 {\bf Proof.} Let $\Phi$ be the automorphism group of $H$. Proposition $7.3.$ implies that there exists a normal subgroup $\Phi_1$ of finite index in $\Phi$ 
such that $\Phi_1(H_i)=H_i\,\, (i=1,2,\cdots)$. We consider now the holomorph of $H$,  that is the group $Hol(H)$ which is a   split extension of the group $H$ by  its  automorphism  group $\Phi$. Let $Q$ be its  subgroup generated by $H$ and $\Phi_1$. Then $Q$ is a split extension of $H$ by $\Phi_1$, and it  is a normal subgroup of finite index in $Hol(H)$; in fact the index of $Q$ in 
$Hol(H)$ is equal to the index $(\Phi\colon \Phi_1)$. We apply now Proposition $7.1.$ and obtain that there exists in $H$ a $p$-series $(7.18)$ whose terms are $\Phi$-invariant and the algebra $L_p(V,V_i)$    isomorphically imbeded into the direct sum of a finite number of isomorphic copies of the algebra  $L_p(H,H_i)$.   This completes the proof.\bigskip

{\bf Corollary 7.3.} {\it The algebra $L_p(V,V_i)$ has the following additional properties.}

 i) {\it  If $L_p(H,H_i)$ is abelian (free abelian) then so is  $L_p(V,V_i)$; if $L_p(H,H_i)$ is abelian of finite rank then so is $L_p(V,V_i)$.}

ii) {\it  Assume that the group $H$ in Proposition $7.4.$ is polycyclic with Hirsch number $r$ and the algebra $L_p(H,H_i)$ is abelian (free abelian) of rank $r$. Then the algebra $L_p(V,V_i)$ is  abelian  (free abelian) of rank $r$.}\bigskip

The first statement follows from Proposition $7.4.$ The second statement is obtained by the same argument as     Corollary $7.2.$\bigskip

{\bf Theorem VIII.} {\it Let $H$ be a finitely generated  group which has  a  $p$-series $H_i\,\, (i=1,2,\cdots)$ with unit intersection and with the associated restricted Lie algebra $L_p(H,H_i)$ abelian (free abelian) of finite rank. Assume that the topology defined by this $p$-series is equivalent to the  $p$-topology. Then there exists a $p$-series 

\begin{equation}H=U_1\supseteq U_2\supseteq\cdots\end{equation}

 whose terms $U_i\,\, (i=1,2,\cdots)$ are characteristic subgroups and the Lie algebra $L_p(H,U_i)$ is abelian (free abelian) of finite rank.}\bigskip

{\bf Proof.}  The assertion follows from Proposition $7.4.$ and Corollary $7.3.$\bigskip 

{\bf Theorem IX.} {\it  Let $H$ be a  torsion free polycyclic group with Hirsch number $r$.  Assume that there exists a  $p$-series $(7.12)$ with unit intersection and associated graded  Lie algebra $L_p(H,H_i)$  free abelian (abelian) of  rank $r$.  Then there exists a $p$-series

 \begin{equation}H=U_1\supseteq U_2\cdots\end{equation}

 whose terms $U_i\,\, (i=1,2,\cdots)$ are characteristic subgroups and the Lie algebra $L_p(H,U_i)$ is free abelian (abelian) of  rank $r$.}\bigskip 

{\bf Proof.} Theorem VI implies that the topology defined by series $(7.12)$ is equivalent to the $p$-topology. The assertion now follows from Theorem VIII and Corollary $7.3$.\bigskip

{\bf 7.3.} {\bf Theorem X.} {\it $H$ be a torsion free   polycyclic group with Hirsch number $r$ which contains a $p$-series  $(7.12)$  
 with unit intersection. Assume  that the Lie algebra $L_p(H,H_i)$ is free abelian of rank $r$. Let $\Phi$ be a group of automorphisms of $H$ such that the order of every automorphism $\phi\in  \Phi$ on  the quotient group $H/H^{\prime}H^p$ is a power of $p$. }

 { \it Then there exists a $p$-series }

\begin{equation} H=H_1^{\ast}\supseteq H_2^{\ast}\supseteq\cdots\end{equation}

 {\it with unit intersection such that  all the subgroups $H_i^{\ast}\,\, (i=1,2,\cdots)$ are $\Phi$-invariant,   $\Phi$ centralizes all the factors $H_i/H_{i+1}\,\, (i=1,2,\cdots)$ and  the algebra $L_p(H,H_i^{\ast})$  is free abelian of rank $r$. }  \bigskip 

We recall  that series $(7.21)$ defines a weight function $f$  in $H$, this weight function defines a valuation $v$ in the group ring $KH$;  conversely, the function $f$   completely defines series $(7.21)$. The construction of the function $f$ and valuation $v$ in the proof of Theorem X will  yield the following fact about the function $f$. \bigskip

{\bf Corollary 7.4.} {\it Let $r$ be the Hirsch number of $H$. Then 
$f(h)\geq  2r\,\, (h\in h)$. } 

Corollary $7.4.$ will not be used in the  proofs of other results.\bigskip

{\bf Proof of Theorem X and Corollary 7.4.} The proof will be given in $4$ steps.

{\it Step 1.} Theorem IX implies that we can assume that all the subgroups $H_i\,\, (i=1,2,\cdots)$ are characteristic. Series $(7.12)$ defines a valuation $\rho$ in the group ring $Z_pH$ with  associated graded ring $gr(Z_pH)$.  We extend now the  valuation $\rho$ to the Laurent polynomial ring $Z_pH[t,t^{-1}]$ assuming that $\rho(t)=1$; let $V$ be the valuation ring of  $Z_pH[t,t^{-1}]$ and $\bar{X}$ be the image of a subset $X\subset V$ under the homomorphism $V\longrightarrow V/(t)$. 
For an arbitrary element $h_j\in H$ we denote $\rho(h_j-1)=n_j\,\, (j\in J)$;  Proposition $2.9.$ implies that the $Z_p$-linear combinations of all the elements $\overline {(h_j-1)t^{-n_j}}\,\, (h_j\in H, j\in J)$ form a subalgebra isomorphic to $L_p(H,H_i)$ and all these linear combinations  form the set of  homogeneous elements of $L_p(H,H_i)$, and that $V/(t)\cong U_p(L_p(H))$.
We will construct at this step some  special free system of generators for this algebra  and a weight function $f$ on this system of generators. 
 
 Proposition $7.3.$ implies that there exists a normal subgroup $\Phi_1$ in $\Phi$ which centralizes all the factors $H_i/H_{i+1}$  and the index $(\Phi\colon \Phi_1)$ is  power of $p$. Since all the subgroups $H_i\,\, (i=1,2,\cdots)$ are characteristic the group $\Phi$ acts as a group of automorphisms of the algebra $L_p(H,H_i)$ and every homogeneous component $H_i/H_{i+1}$ is $\Phi$-invariant; since  the subgroup $\Phi_1$ acts trivially we obtain that    the finite $p$-group $G=\Phi/\Phi_1$ acts as a group of automorphisms of the algebra $L_p(H,H_i)$ and of the vector spaces $H_i/H_{i+1}\,\,$  so   these vector spaces become $G$-modules.  

 Lemma $2.11$ implies that the free abelian algebra $L_p(H,H_i)$ is generated by the vector space   $Q\cong L_p(H,H_i)/L^{[p]}_p(H,H_i)$  which is in fact a subspace of $L_p(H,H_i)$; this vector space has dimension $r$ because the rank of $L_p(H,H_i)$ is $r$ and the group $G$ acts in a natural way on this vector space. The grading in $L_p(H,H_i)$ defines a grading in the vector subspace $Q$. Since the homogeneous components of $L_p(H,H_i)$ are $G$-invariant   we obtain that  the homogeneous components  of $Q$  are also $G$-invariant. 
  Since $Q$ has finite dimension $r$ there will be a finite number $k\leq r$ of non-zero homogeneous components $Q_i\,\, (i=1,2,\cdots,k)$ 

\begin{equation} Q=\bigoplus_{i=1}^k  Q_i\end{equation}

 It is worth remarking that every  element $x\in Q_i\,\, (i=1,2,\cdots,k)$   is a homogeneous element in $L_p(H,H_i)$.

We pick now an arbitrary homogeneous component $Q_i$ and consider the    series 
\begin{eqnarray} Q_i\supseteq \omega(KG) \bullet Q_i\supseteq\omega^2(KG)\bullet Q_i\supseteq\cdots\nonumber\\
\cdots\supseteq \omega^{s-1}(KG)\bullet Q_i  \supseteq \omega^s(KG)\bullet Q_i=0\end{eqnarray}

The inclusions in this series are strict because if two consecutive terms had coincided, say   $\omega^n(KH)\bullet Q_i= \omega^{n+1}(KH)\bullet Q_i\not=0$,  then we would have gotten  that $\omega^j(KH)\bullet Q_i=\omega^{j+1}(KH)\bullet Q_i$ for all $j\geq n$ which is impossible because the ideal $\omega (KG)$ is nilpotent. Since the inclusions  are strict and the dimension of $Q=r$  we obtain that in $(7.23)$  $s\leq r$.

We pick  now in  $\omega^{n}(Z_pG)\bullet Q_i$ a system of elements  $\bar{T}_{i,n}$ which  give   a basis of the quotient module  
$(\omega^{n}(Z_pG)\bullet Q_i)/ (\omega^{n+1}(Z_pG)\bullet Q_i )  \,\, (n=0,1,\cdots, s-1)$.  It is also important that all the elements in the system $\bar{T}_{i,n}$ are homogeneous elements of $L_p(H,H_i)$.
The system of elements $\bar{T}_i=\bigcup_{n=1}^{s-1}\bar{T}_{i,n}$ forms a basis for $Q_i$.  We obtained    a basis $\bar{T}_i$ for every vector subspace $Q_i\,\, (i=1,2,\cdots,k)$ and 
  obtain then a basis $\bar{T}$ for $Q$ taking 
$\bar{T}=\bigcup _{i=1}^k \bar{T}_i$, where all the elements of $\bar{T}$ are homogeneous and $\bar{T}$ freely generate the algebra $L_p(H,H_i)$.

 Let $m$ be an arbitrary integer greater than $2r$. We define a weight function $f$ on $\bar{T}_{i,n}$ as follows 

\begin{equation} f(x)=m+2n+1\,\, \mbox{if}\,\, x\in \bar{T}_{i,n}\,\,(n=0,1,\cdots,s-1)  \end{equation}

Since $\bar{T}$ is a disjoint union of the systems $\bar{T}_{i,n}$ we obtain a weight function $f$ on $\bar{T}$. It follows immeadiately that

\begin{equation}m+1\leq  f(x)\leq m+2r+1\,\, (x\in \bar{T})\end{equation}

We complete this step by reminding  that it was pointed out in the beginning of the proof   that $V/(t)\cong Z_p[\bar{T}]$ and that every element $\bar{t}_j\in \bar{T}$ has in fact a form $\bar{t}_j=\overline{(h_j-1)}t^{-n_j}\,\, (j=1,2,\cdots,r)$ where $n_j=\rho(h_j-1)$.

{\it Step 2.} We extend now the weight function $f$ which was defined on the system of generators $\bar{T}$ to a valuation $\bar{v}$ of the algebra 
$Z_p[\bar{T}]\cong \bar{V}$; this valuation will be used on Step $3$  for a construction of a valuation $v$ in $V$ and in $Z_pH$.  We will show now that the group $G$ centralizes the  valuation $\bar{v}$.

     Let  $x\in \bar{T}$.
The definition of $\bar{T}$ implies that there exists a unique pair $i,n$ such that $x\in \bar{T}_{i,n}$ so $\bar{v}(x)=m+2n+1$.
We see now from $(7.23)$  that if  $x\in \bar{T}_{i,n}$ the the image of the element $(g-1)\bullet x$ in $Q$ belongs to $\omega^{n+1}\bullet Q$. We obtain from this that   either   there exists $0\not=u$  which is a linear combination of elements from $\bigcup_{k=n+1}^{s-1} \bar{T}_{i,k}$ and $ y\in Q^p$ such that 

\begin{equation} g\bullet x-x=(g-1)\bullet x=u+y\end{equation}

or 

\begin{equation} g\bullet x-x=(g-1)\bullet x=y\end{equation}

We will now show that $\bar{v}(u)$ and $\bar{v}(y)$ are greater than $m+2n+2$. This will imply that the element of $G$ centralize the elements of $\bar{T}$ with respect to the valuation $\bar{v}$.

In fact, we obtain from $(7.24)$ that the $\bar{v}$-values of the elements from $\bar{T}_{i,l}$ are greater than or equal than $m+2n+3$ if $l\geq (n+1)$. Further, we see from relation  $(7.25)$ that the values of the function $f(x)$ on the system  $\bar{T}$  are greater than or equal $m+1$; this together with  the   definition of the function $\bar{v}$  implies  that  the $\bar{v}$-values of all the elements of $Q$ are greater than or equal $m+1$. Since the element $y$ belongs to the subalgebra $Q^p$ 
the  value of the element $y$  is greater than or equal to $p(m+1)\geq 2(m+1)>m+2r+1$. 
We obtain from this that in both cases, when  relation $(7.26)$ or $(7.27)$ hold,  

\begin{equation}\bar{v}(g\bullet x- x)> \bar{v}(x)+1\end{equation}

and our claim is proven.

{\it Step 3.} Let $\bar{t}_j=\overline{(h-1)}t^{-n_j }$ be an arbitrary element of $ \bar{T}$.  We denote now $(h_j-1)t^{-n_j}\,\, (j=1,2,\cdots,r)$ and obtain a system of elements $<t_1,t_2,\cdots,t_r>=T\in Z_pH$ such that its image in $V/(t)$ is $\bar{T}$.

Consider now  the system of elements $<t,T>$ in $V$. Since the algebra $V/(t)$   is isomorphic to the polynomial agebra $Z_p[T]$ the system  $<t,T>$ is an independent polycentral system in $V$. We  extend now the weight function $f(x)$ to the system $<t,T>$ by defining $f(t)=M$ where $M$ is an arbitrary natural number greater than $2( m+2r+1)$ and $f(t_j)=f(\bar{t}_j).$ Since the  values of $f$ on the subsystem 
$T$  are less than or equal to $ m+2r+1$ we obtain that $f(t)>2f(x)$ for every $x\in T$ and we can apply Theorem II and obtain  that this  weight function extends to a valuation $v$ in $V$ with graded ring isomorphic to the polynomial ring $Z_p[t,T]$. 

We will prove now  that the group $G$ centralizes the graded ring $gr_v(V)$. 
Since $\phi(t)=t\,\, (t\in T)$ by definition, we have to prove only that $G$ centralizes the elements from $T$. Let $t_j$ be an element from $T$.  We  have once again  $\bar{t}_j\in \bar{T}_{i,n}$ for some pair $ i,n$. If  relation $(7.26)$ holds in the ring $V/(t)\cong Z_p[T] $ then we obtain in $V$ 

\begin{equation} g\bullet t_j-t_j=(g-1)\bullet t_j=u+y+u_1\end{equation}
  
where $u$ and $y$ are the same as in $(7.26)$ and $u_1$ is an element from the ideal $(t)$. Since all the elements from the ideal $(t)$ have values greater than or equal $2(m+2r+1)$ we obtain from $(7.26)$ and $(7.29 )$ that 
$v(( g-1)\bullet t_j) >v(t_j)+1$; the same relation is obtained in the case when  
$(7.27)$ holds. This shows that $G$ 
centralizes the system of elements $t, T$; hence it centralizes the ring $gr_v(V)$. 

{\it Step 4.} We restrict now the valuation function $v$ to the subring $Z_pH\subseteq V$.  The graded ring $gr_v(Z_pH)$ is a subring of $gr_v(V)$, so 
$gr_v(Z_pH)$  is centralized by $\Phi$.  Further, we denote

\begin{equation}H_i^{\ast}=\{h\in H|v(h-1)\geq i\}\,\, (i=1,2,\cdots)\end{equation}

and obtain o in $H$ a new $p$-series with unit intersection 

\begin{equation}H=H^{\ast}_1\supseteq H^{\ast}_2\supseteq\cdots\end{equation}

Proposition $2.10.$ implies that the homogeneous components  $\widetilde{(h-1)} \,\, (h\in H)$ generate in $gr_v(Z_p(H))\subseteq gr_v(V)$  a subalgebra isomorphic to the algebra $L_p(H^{\ast},H^{\ast}_i)$. Since $gr(V)$ is centralized by 
$\Phi$  the elements of this  subalgebra  are centralized by the group $\Phi$. Since $gr_v(V)\cong Z_p[t,T]$ we obtain that the algebra $L_p(H,H_i^{\ast}$ is abelian and it contains no nilpotent elements by  Corollary $2.6.$, so it is free abelian.

We prove now  that the rank of  $L_p(H,H_i^{\ast})$ is $r$.  The system of homogeneous components $t$ and   $ \tilde{t}_j=\widetilde{(h_j-1)}t^{-n_j}\,\, (j=1,2,\cdots,r)$ generate a ring isomorphic to the polynomial ring $Z_p [t,T]$. This implies that the system of elements 
$\tilde{t}_j t^{-n_j}=\widetilde{(h_j-1)}\,\, (j=1,2,\cdots,r)$ is algebraically independent in $gr(Z_pH)$ over $Z_p$.  Since $U_p(L_p(H,H_i^{\ast})\cong gr(Z_pH)$ we obtain that $U_p(L_p(H,H_i^{\ast})$ is a polynomial ring with $r_1\geq r$ variables. 
On the other hand, the rank of the algebra $L_p(H,H_i^{\ast})$  does not exceed $r$ by Proposition $ 2.5. $  so $r_1$ can not be greater than $r$. We obtain therefore that this rank is $r$. 

This completes the proof of Theorem X. \bigskip

 {\bf Lemma 7.3.} {\it Let $G$ be a finite $p$-group, $(G\colon 1)=p^n$. There exists a $p$-series}

\begin{equation} G=G_1\supseteq G_2\supseteq\cdots\supseteq G_{n-1}\supseteq G_n=1\end{equation}

{\it such that the restricted Lie algebra $L_p(G,G_i)$ is abelian of dimension $n$ and exponent $p$.}\bigskip

{\bf Proof.} Let $U$ be a central subgroup of order $p$ in $G$ generated by an element $u$. We can assume that there exists a $p$-series in the quotient group $\bar{G}=G/U$

\begin{equation} \bar{G}=\bar{G}_1\supseteq \bar{G}_2\supseteq\cdots\supseteq  \bar{G}_{n-2}\supseteq \bar{G}_{n-1}=1\end{equation}

such that the algebra $L_p(\bar{G},\bar{G}_i)$ is free abelian of dimension  $n-1$ and of exponent $p$. Let $m$ me the maximum of the weights of non-unit elements of $\bar{G}$. We pick now an integer $M>pm$ and for 
$u^n\in U\, \, (1\leq n\leq p-1)$  define the weight  of $u^n$ in $G$ as $\omega(u)=nM$,    $\omega(1)=\infty$. We define then for an element $g\not\in U $ its  weight   $\omega(g)$ to be equal to the weight of the coset $\bar{g}=gU$ in $\bar{G}$. The weight function $\omega$ is now defined for all the elements $g\in G$ and $\omega(gu^k)=\omega(g)$ for an arbitrary $g\in G$ and a natural $k$. 

Let $g_{i_1},g_{i_2}$ be two elements of $G$. If the commutator $[g_{i_1},g_{i_2}]$ does not belong to $U$ then 

\begin{equation} \omega([g_{i_1},g_{i_2}])>  \omega(g_{i_1})+\omega(g_{i_2}) \end{equation}

because this relation holds in $\bar{G}$ for the elements $\bar{g}_{i_1},\bar{g}_{i_1}$. On the other hand, if $[g_{i_1},g_{i_2}]\in U$ then relation $(7.34)$ holds because the elements of $U$ have weights greater than $pm\geq \omega(g_{i_1})+\omega(g_{i_2})$.

The same argument shows that 

\begin{equation} \omega(g^p)> p \omega(g)\end{equation}

Relations $(7.34)$ and $(7.35)$ show that the weight function $\omega$ defines a $p$-series with associated graded algebra abelian of exponent $p$. Since $1$ is the only element with infinite weight the intersection of all the terms of this series is $1$ and the dimension of the associated Lie algebra is $n$. 

This completes the proof.\bigskip

{\bf Proposition 7.5.} {\it Let $H$ be a polycyclic with Hirsch number $r$, $U$ be a normal subgroup such that the quotient group $G=H/U$ is a finite $p$-group, $(G\colon U)=p^n$. Assume that there exists a $p$-series 

\begin{equation} U=U_1\supseteq U_2\supseteq\cdots\end{equation} 

with unit intersection with the associated graded   algebra $L_p(U,U_i)$ is free abelian (abelian) of rank $r$.  }

 {\it Then there  exists a $p$-series }

\begin{equation} H=H_1\supseteq H_2\supseteq\cdots\end{equation}

{\it with unit intersection such that the algebra $L_p(H,H_i)$ is abelian of rank $r$ and  the algebra 
$L_p(U,U_i^{\ast})$ associated to the $p$-series $U_i^{\ast}=U\bigcap H_i\,\, (i=1,2,\cdots)$ is free abelian of rank $r$ (abelian of rank $r$).  }\bigskip

{\bf Proof.} {\it Step 1.} We consider first the case when the algebra $L_p(U,U_i)$ is free abelian of rank $r$.  We pick in the group $G$ a weight function $\omega(G)$ obtained in Lemma $7.3.$ with the associated restricted Lie algebra abelian of dimension $n$.
Consider the wreath product $W=U\wr G$ and obtain from 
Lemma  $7.2.$  that there exists  in the base group $U^{\ast}$
a  $W$-invariant $p$-series  $(7.9)$  such that the Lie algebra $L_p(U^{\ast},U^{\ast}_i)$ is isomoprphic to the direct sum of $p^n$ copies of $L_p(U,U_i)$; so it is free abelian of rank $p^nr$.    Theorem X  implies that we can get in 
$U^{\ast}$ a $G$-invariant $p$-series $V_i\,\, (i=1,2,\cdots)$ such that the group $G$ centralizes the algebra $L_p(U^{\ast},V_i)$ and the rank of $L_p(U^{\ast},V_i)$ is $p^nr$. Let $\rho$ be the weight function defined by this series and  let $m$ be the maximum of the weights $\omega(g)\,\, (g\in G)$.  We pick an  arbitrary integer $M$ greater than $pm$. We can assume (see Corollary 
$2.7.$) that the values of the weight function $\rho$ are multiples of $M$.

 We define now a weight function $\Omega$ on $W$ as follow. 

\begin{equation} \Omega(u)=\rho(u)\,\, \mbox{if}\,\, u\in U^{\ast}; \Omega(gu)=\omega(g)\,\, \mbox{for}\,\, 1\not=g\in G; u\in U^{\ast}\end{equation}
 
The restrictions of $\Omega$ on the groups $G$ and $U^{\ast}$ define algebras $L_p(G,G_i)$ and $L_p(U^{\ast},V_i)$ respectively. This fact together with  the definition of $\Omega$ implies that for every two elements $u_1,u_2$ we have 

\begin{equation}\Omega[u_1,u_2]> \Omega(u_1)+\Omega(u_2)\end{equation} 

in the following $2$ cases:

$1)$ If $u_1,u_2\in U^{\ast}$; $2)$ If $u_1,u_2\in G$. 

We consider now the third case when when $u_1=g_1v_1, u_2=g_2v_2$ where $g_1,g_2$ are non-unit elements of $G$, $v_1,v_2$ are non-unit elements of $U$. In this case $\Omega(u_i)=g_i\,\, (i=1,2)$ and this case is reduced to  case $2$. 

It remains to prove relation 
$(7.39)$ in the case when $u_1\in G $ and $u_2\in U^{\ast}$.  Since $G$ centralizes the algebra $L_p(U^{\ast},V_i)$ we have in this case $\rho(u_1^{-1}u_2^{-1}u_1u_2)> \rho(u_2)+1$. Since the value of $\rho$ are multiples of $M>m$ we obtain that 

 \begin{equation} \rho(u_1^{-1}u_2^{-1}u_1u_2)\geq \rho(u_2)+M>\rho(u_1)+\rho(u_2)\end{equation}

which proves $(7.39)$.

The same argument proves the relation 

\begin{equation} \Omega(x^p)\geq p\Omega(x)\,\, \mbox{for}\,\, (x\in W)\end{equation}

 We see that the weight function $\Omega$  defines a $p$-series 

\begin{equation} W=W_1\supseteq W_2\supseteq\cdots\end{equation} 

with unit intersection and $W_i\bigcap U^{\ast}=V_i; W_i\bigcap G=G_i\,\, (i=1,2,\cdots)$ and the algebra $L_p(W,W_i)$ is a split extension of the free abelian algebra $L_p(U^{\ast}, V_i)$ of rank $p^nr$ by a finite abelian algebra $L_p(G,G_i)$.  We obtain from this together with relation $(7.39)$  that the algebra $L_p(W,W_i)$  is a direct sum of the algebra $L_p(U^{\ast},V_i)$ and  the algebra $L_p(G,G_i)$, so it  is abelian of rank $p^nr$ which is equal to the Hirsch number of $W$. Since $H$ is a subgroup of $W$  statement i) of Theorem VI implies that the rank of the algebra $L_p(H,H_i)$ must coincide with its Hirsch number of $H$. Since the restriction of the function $\Omega$ on the group $U^{\ast}$ defines the free abelian algebra  $L_p(U^{\ast}, U^{\ast}_i)$ of rank $p^nr$ equal to the Hirsch rank of $U^{\ast}$ its restriction on $U$ defines a free abelian algebra $L_p(U,V_i)$, its rank  must be equal $r$ by statement i) of Theorem VI.

This completes the proof for the case when the algebra $L_p(U,U_i)$ is free abelian of rank $r$. 

{\it Step 2.} Now assume that $L_p(U,U_i)$ is abelian of rank $r$. Apply Theorem $6.1.$ and find $U_i=Q$ such that the algebra $L_p(Q,Q_i)$ associated to the $p$-series $Q_i=U_i\bigcap Q\,\, (i=1,2,\cdots)$ is free abelian of rank $r$. We find then in $U$ a characteristic subgroup $R\subseteq Q$ with index $(U\colon R)=p^n$. The algebra $L_p(R,R_i)$ associated to the $p$-series $R_i=R\bigcap Q_i=R\bigcap U_i\,\, (i=1,2,\cdots)$ is free abelian of rank $r$ because it is a subalgebra of $L_p(Q,Q_i)$. Since the index $(H\colon R)$ is a power of $p$ the assertion follows from the case which was considered at step $1$. 
\bigskip

{\bf Theorem XI.} {\it Let $H$ be a polycyclic group with Hirsch number $r$. Assume that there exists a $p$-series $(1.1)$ with unit intesection such that the algebra $L_p(H,H_i)$ is finitely generated.  Then the center $Z$ of $L_p(H,H_i)$ has rank $r$ iff the  following two conditions hold}

i) {\it  $H$ contains normal subgroups $Q\supseteq N$ which satisfy the conditions of Theorem VII.}

ii) {\it The topology defined by series $(1.1)$ is equivalent to the $p$-topology. }\bigskip

{\bf Proof.}   If conditions i) and ii) hold then the assertion follows from Theorem VII and statement viii) of Theorem VI.

 We prove now the necessity of these conditions.
Assume that the rank of $Z$ is $r$. Theorem V implies that there exists a normal subgroup $U$ of index $p^n$ such that the algebra $L_p(U,U_i)$ is free abelian subalgebra of $Z$ of finite index in $L_p(H,H_i)$.  The index of $L_p(U,U_i)$ in $Z$ is also finite so the rank of $L_p(U,U_i)$ is $r$. Proposition $7.5.$  
implies that there exists in $H$ a $p$-series with unit intesection and associated graded  Lie algebra abelian of rank $r$. The necessity of the conditions i) and ii)  now follows from statements v) and vi) of Theorem VI.

 \bigskip

\setcounter{section}{8}
\setcounter{equation}{0}

\section*{\center \S 8. Proof of Theorem XII.}

{\bf 8.1.} {\bf Proposition 8.1.} {\it Let $R$ be an algebra over a field $K$ with a non-negative discrete pseudovaluation $\rho$, and associated graded ring $gr(R)$, $R\ast H$ be a suitable  skew group ring of $R$ with infinite cyclic group $H$. Assume that there exists a natural number $k$ such that for every $x\in R$

\begin{equation} \rho(hxh^{-1}-x)>\rho(x)+k\end{equation}

Then the pseudovaluation $\rho$ extends to a pseudovaluation $\rho_1$ of $R\ast H$ such that $\rho_1(h-1)=k$. The graded ring $gr_{\rho_1}(R\ast H)$ is isomorphic to the polynomial ring 
$gr(R)_{\rho}[t]$ where $t=\widetilde{(h-1)}$ is the homogeneous component of the element $h-1$. The extension $\rho_1$ is the only extension of $\rho$ with 
value of $h-1$ equal $k$ and associated graded ring isomorphic to $gr_{\rho}(R)[t]$. }

{\bf Proof.} Let $h$ be the generator of $H$. We extend first $\rho$ to the skew polynomial subring $R[h]$. Every element  $x\in R[h]$ has a unique representation

\begin{equation}x=\sum_{i=0}^n\lambda_i(h-1)^i\,\, (\alpha_i\in R;\,\,i=0,1,\cdots,n)\end{equation}

We  define now  the value of element $(8.2)$   by 

\begin{equation}\rho_1(x)=\min_i\{\rho(\lambda_i+ki)\}\end{equation} and $\rho_1(0)=\infty$. 

We will now prove that $\rho_1$ is a pseudovaluation on in $R[h]$. Assume that 
$\rho_1(x)=\rho(\alpha_{i_0})+i_0k$,  and let $0\not=y=\sum_{j=1}^m \beta_j(h-1)^j$ be an element of $R\ast H$ with $\rho_1(y)=\rho(b_{j_0})+j_0k$. We see immediately that

\begin{equation}\rho_1(x+y)\geq \mbox{min}\{\rho_1(x),\rho_1(y)\}\end{equation}

We prove now the relation 
\begin{equation}\rho_1(xy)\geq\rho_1(x)+\rho_1(y)\end{equation}

We have $x=x_1+x_2,y=y_1+y_2$ where all the summands in the representations of $x_1,y_1$ have values $\rho_1(x_1)$ and $\rho_1(y_1)$ respectively and $\rho_1(x_2)>\rho_1(x), \rho_1(y_2)>\rho_1(y)$. Relation $(8.5)$ will follow if we prove   that 

\begin{equation}\rho_1(x_1y_1)\geq \rho_1(x_1)+\rho_1(y_1)\end{equation}

Let $\alpha (h-1)^{i}, \beta (h-1)^{j}$ be two arbitrary  terms  in the representations of $x_1,y_1$.  We have $(h-1)\beta-\beta (h-1)=h\beta-\beta h=
(h\beta h^{-1}-\beta)$ and we obtain from this and $(8.1)$ that

\begin{equation}(h-1)\beta=\beta(h-1) +u\end{equation} 

where $\rho(u)>\rho(\beta)+k$. We will use now an  induction argument to show that for every natural $i$ 

 \begin{equation} (h-1)^i\beta=\beta(h-1)^i+u_i\end{equation}  

where $u_i$ is an element of $R\ast H$ with $\rho_1(u_i)>\rho(\beta)+ki$

In fact, assume that it has already been proven that 
\begin{equation} (h-1)^{i-1}\beta=\beta (h-1)^{i-1}+u_{i-1}\end{equation}

where $u_{i-1}$ is an element of $R\ast H$ such that $\rho_1(u_{i-1})>\rho(\beta)+k(i-1)$.

We obtain from this

\begin{eqnarray} (h-1)^i\beta=(h-1)\beta(h-1)^{i-1}+(h-1)u_{i-1}=\nonumber\\
=(\beta(h-1)+u)(h-1)^{i-1}+  (h-1)u_{i-1}= \nonumber\\
=\beta(h-1)^i+u(h-1)^{i-1}+(h-1)u_{i-1}
\end{eqnarray}

We denote now $u_i=u(h-1)^{i-1}+(h-1)u_{i-1}$ and obtain $(8.8)$.

 We have now from $(8.8)$ 

 \begin{eqnarray}\alpha(h-1)^i\beta(h-1)^j=\alpha(\beta(h-1)^{i+j}+u_i(h-1)^j)=\nonumber\\
=\alpha\beta(h-1)^{i+j}+\alpha u_i (h-1)^j=\alpha\beta(h-1)^{i+j}+v\end{eqnarray} 
where $v=\alpha u_i(h-1)^j$ and 

\begin{equation}\rho_1(v)>\rho(\alpha)+\rho(\beta)+k(i+j)\end{equation}

Definition $(8.3)$ implies that the    summand $\alpha\beta (h-1)^{i+j}$ has value $\rho(\alpha \beta) +k(i+j)\geq \rho(\alpha)+\rho(\beta)+k(i+j)$. This and $(8.12)$ imply that   the right side of $(8.11)$ has value greater than or equal 
$\rho(\alpha)+\rho(\beta)+k(i+j) $ and we obtain from this and $(8.11)$ 

\begin{equation}\rho_1(\alpha(h-1)^i\beta(h-1)^j)\geq \rho_1(\alpha(h-1)^i)+\rho_1(\beta(h-1)^j)\end{equation} 

Since this relation holds for arbitrary summands in the representations of $x_1,y_1$ we obtain that $\rho_1(x_1y_1)\geq \rho(x_1)+\rho_1(y_1)$; this  together with $(8.4)$ proves that $\rho_1$ is a pseudovaluation in $R[h]$ which extends the valuation $\rho$ of $R$. 

The element $h-1$ has value $k$ and relation $(8.7)$ implies that the homogeneous component $t=\widetilde{(h-1)}$ of $h-1$ commutes with all the homogeneous components $\tilde{\beta}$ of elements $(\beta\in R)$ which implies that  $t$ commutes with the subring 
$gr_{\rho}(R)\cong gr_{\rho_1} (R)$. Further, if $\rho_1(x)=l$ then $x=u+x_2$ 
where $\rho_1(x_2)>l$,

\begin{equation} u=\sum_{i=0}^n\alpha_i(h-1)^i\end{equation}

 and $\rho(\alpha_i)+ki=l$; we see that the  homogeneous component of $x$ is equal to the homogeneous component of $x_1$.  We will show now that this  homogeneous component  has a unique representation 

\begin{equation} \tilde{u}=\sum_{i=0}^n\tilde{\alpha}_it^i\end{equation}

This will imply that $gr_{\rho_1}(R[h])$ is isomorphic to the polynomial ring $gr_{\rho}(R)[t]$.

We observe first that the definition of $\rho_1$ implies that the weight of every summand $\alpha_i(h-1)^i $ in $(8.14)$ is equal
to $\rho_1(\alpha)+ki=l$ and Lemma $2.3$ implies that the homogeneous component of this summand is $\tilde{\alpha}_it^i$. Since the weight of $x_1$ is $l$ we conclude now once again from Lemma $2.3$  that $\tilde{x}_1$ does have representation $(8.15)$, which proves that $gr_{\rho_1} (R[h])\cong gr(R)[t]$. 

To extend $\rho$ to the ring $R\ast H$ we use the fact that for every element $y\in R\ast H$ we can find $x\in R[h]$ and an integer $m$ such that  $y=xh^m$. We define now $\rho_1(y)=\rho_1(x)$ and a straightforward argument shows that that we obtain a pseudovaluation in $R\ast H$ and that the graded ring 
$gr_{\rho_1}(R\ast H)$ is isomorphic to $gr_{\rho_1}(R[h])$.

 To prove that $gr_{\rho_1}(R\ast H)\cong gr_{\rho} (R)[t]$ we  have to prove   the uniqueness of representation $(8.15)$; once again, it is enough to do this for the subring $R[h]$. 

If representation  $(8.15)$ is not unique then  a standard argument yields that there exist non-zero elements $\lambda_j\in R\,\, (j=1,2,\cdots,n)$ such that 

\begin{equation}\sum_{j=0}^m\tilde{\lambda}_jt^j=0 \end{equation}

where

\begin{equation}\rho(\lambda_j)+kj=l\,\, (j=1,2,\cdots,m)\end{equation} 

 Equation $(8.16)$ means that 

\begin{equation} \rho_1 (\sum_{j=0}^m \lambda_j(h-1)^j)\geq (l+1)\end{equation}

This contradicts the definition $\rho_1$ and we proved the uniqueness of representation $(8.15)$.

It remains to prove that $\rho_1$ is the only extension of $\rho$ such that $\rho(h-1)=k$ and $gr_{\rho_1}(R\ast H)\cong gr_{\rho}(R)[t]$. The  proof of this fact is obtained by the same argument, with obvious simplifications, which   will be used in step $2$ of  the proof of a similar statement in Proposition $8.2.$ and  we omit it. \bigskip

{\bf Proposition 8.2.}  {\it Let $R_1$  be an algebra of characteristic $p$, $R$ be a subalgebra of $R_1$.   Assume that there exists an invertible element 
$g\in R_1$  such that $g^{-1}Rg=R, g^{p^m}=r_0\in R$ and    the elements }

\begin{equation} 1,g,g^2,\cdots, g^{p^m}-1\end{equation}

{\it form a basis of the left $R$-module $R_1$. Assume also that there exists   
 a non-negative discrete pseudovaluation $\rho$, with commutative  associated   graded ring $gr(R)$
such that $\rho(r_0-1)=kp^m$ and }

 \begin{equation} \rho(g^{-1}rg-r)>\rho(r)+k\end{equation} 

{\it for every $r\in R$. } 

{\it Then there exists an extension of $\rho$ to a pseudovaluation $\rho_1$ of  the ring $R_1$ such that $\rho(g-1)=k$, the graded ring $gr_{\rho_1}(R_1)$ is isomorphic to the algebraic extension $gr_{\rho}(R)[\theta]$ where $\theta=\widetilde{(g-1)}$ is the homogeneous component of $g-1$, the minimal polynomial of $\theta$ is $t^{p^m}-\widetilde{(r_0-1)}$.

 The pseudovaluation  $\rho_1$ is the unique extension such that the value of $r_0-1$ is $k$ and  the  asssociated graded ring is isomorphic to the quotient ring 
$(gr_{\rho}(R)[t])/(t^{p^m}-\widetilde{(r_0-1)})$.}\bigskip

{\bf Proof.} {\it Step 1.}  Let $\phi$ be the inner automorphism $x\longrightarrow g^{-1}xg$ of $R_1$. Since the subring $R$ is $\phi$-invariant we can consider the skew group ring $R\ast H$ of $R$ with infinite group $H$ where $h^{-1}rh=\phi(r)$ for the generator $h$ of $H$ and an  arbitrary $r\in R$. We apply Proposition $8.1$ to obtain a  pseudovaluation $\tau$ of $R\ast H$ which coincides with $\rho$ on $R$,  $\tau(h-1)=k$, and $gr_{\tau}(R\ast H)\cong R[t]$. We consider now the  homomorphism $\psi\colon R\ast H\longrightarrow R_1$ whose kernel is the principal  ideal  $A$ generated by the element $(h-1)^{p^m}-(r_0-1)$.  Proposition $2.2.$ implies that we obtain  a  pseudovaluation  $\rho_1$ of the ring $R_1$. We will show that this pseudovaluation satisfies all the conclusions of the assertion.

We will  prove first the relation 

\begin{equation} gr_{\rho_1}(R_1)\cong R/(t^{p^m}-\widetilde{(r_0-1)})\end{equation}

The homogeneous component of the element $(h-1)^{p^m}-(r_0-1)$ in 
$gr_{\tau}(R\ast H)$ is $t^{p^m}-\widetilde{(r_0-1)}$. So the element $t^{p^m}-\widetilde{(r_0-1)}$  belongs to the kernel of the homomorphism $gr_{\tau}(R\ast H)\longrightarrow gr_{\rho_1}(R_1)$.  To prove  relation $(8.21)$ it is enough to verify that if  $x\in ker(\phi)=(h^{p^m}-r_0)$ then  $\tilde{x}\in (t^{p^m}-\widetilde{(r_0-1)})$.  Assume that $\tau(x)=l$ and let 

 \begin{equation} x=\sum_{i,j}x_i(h^{p^m}-r_0)x_j\end{equation}

Since the graded ring $gr_{\tau}(R\ast H)\cong R[t]$ is commutative we obtain from 
$(8.22)$ that 

\begin{equation} x=\sum_{i,j}(h^{p^m}-r_0)x_ix_j+y\end{equation}

where $\tau(y)>l$. We see that the homogeneous component of $x$ coincides with 
the homogeneous component of the element $(h^{p^m}-r_0)a$ where  $a=\sum_{i,j}x_ix_j$; we can assume in fact that $x= (h^{p^m}-r_0)a$. We have already observed that  the homogeneous component of $(h^{p^m}-r_0)$ is equal to $\widetilde{(h^{p^m}-r_0)}=(\widetilde{(h^{p^m}-1)}- \widetilde{(r_0-1)})=(t^{p^m}-\widetilde{(r_0-1}))$. Since this element can not be a zero divisor in the ring 
$gr_{\rho}(R)[t]$ we obtain that the homogeneous component of $x$ is a product of homogeneous components of the element $h^{p^m}-r_0$
and  of the element $a$: 

\begin{equation}\tilde{x}=  (t^{p^m}-\widetilde{(r_0-1}))\tilde{a}\end{equation} 
which means that $\tilde{x}$ belongs to the ideal generated by $ (t^{p^m}-(r_0-1))$. This  proves  $(8.21)$.

{\it Step 2.} We will prove now the uniqueness of the extension. If $\rho_2$ is a second extension of the pseudovaluation $\rho$ with associated graded ring isomorphic to 
$(gr_{\rho}(R)[t])/(t^{p^m}-\widetilde{(r_0-1)})$
then $\rho_2(g-1)$ must be less  than or equal  to $k$ because $(g-1)^{p^m}=r_0-1$. If $\rho_2(g-1)<k$ then the homogeneous component $\widetilde{(g-1)}$ is nilpotent and the minimal polynomial of this component can not be   
$t^{p^m}-\widetilde{(r_0-1)}$.   We see that $\rho_2(g-1)=k$. Further  for every $\lambda\in R$  and $(g-1)^i\,\, (1\leq i\leq p^m-1)$ we  have $\rho_2(\lambda(g-1)^i)\geq \rho(\lambda)+ki$; if we had in the last equation a strict inequality    we would have 
$\tilde{\lambda} \widetilde{ (g-1)^i)}=0$ which is impossible because we assumed that the minimal polynomial of $g-1$ is $t^{p^m}-\widetilde{(r_0-1)}$. We obtain therefore that $\rho_2(\lambda(g-1)^i)= \rho(\lambda)+ki$.

Let 

\begin{equation} u=\sum_{i=1}^n\lambda_i(g-1)^i\,\, (i\leq p^m-1)\end{equation}

be an arbitrary element of $R_1$, and  let  $ \mbox{min}_{1\leq i\leq n}\,\,  
\{\rho(\lambda_i)+ik)\}=l$; we can assume that 
$\rho(\lambda_i)+ik=l\,\, \mbox{if}\,\, 1\leq i\leq n_1$ and $\rho(\lambda_i)+ki>l\,\, \mbox{if}\,\, n_1<i\leq n$. Let $u_1= \sum_{i=1}^{n_1}\lambda_i(g-1)^i$. Clearly $\rho_1(u)=\rho_1(u_1)$,   $\rho_2(u)=\rho_2(u_1)$ and  $\rho_1(\lambda_i(g-1)^i)=\rho_2(\lambda_i(g-1)^i)=l\,\, (1\leq i\leq n_1)$.  The  element $\sum_{i=1}^{n_1}\tilde{\lambda}_i\widetilde{(g-1)^i}=\sum_{i=1}^{n_1}\tilde{\lambda}_it^i$ of $gr_{\rho_1}(R_1)$ is nonzero and  we conclude now from Lemma $2.3$   that the $\rho_1$-value of $u$ must be $l$. The same argument shows thar $\rho_2(u)=l$. We proved that the pseudovaluation $\rho_2$ must coincide  with $\rho_1$.
This completes the proof.\bigskip

We will need in the proof of Theorem XII a corollary of Propositions $8.1.$ and $8.2.$\bigskip 

{\bf Corollary 8.1.} {\it   Assume that the conditions of Proposition $8.2.$ hold. Let   $A$  be  an ideal of $R$ such that $g^{-1}Ag=A$ and $A_1=AR_1$ be the ideal in $R_1$ generated by $A$. Let $x$ be an arbitrary element of $A_1$. Then the homogeneous component $\tilde{x}$ of $x$ in $gr_{\rho_1}(A_1)$ has a unique representation 

\begin{equation} \tilde{x}=\sum_{i=1}^m \tilde{\mu}_i \widetilde{(g-1)}^{n_i}  \end{equation}

where $\tilde{\mu}_i$ is the  homogeneous component of an element $\mu_i\in A$,$\widetilde{(g-1)}$  is the homogeneous component of $g-1$,  
$0\leq n_i\leq p^n-1$. So the ideal $gr_{\rho_1}(A_1)$ of the ring $gr(R_1)$ is isomorphic to the ring   $A[\theta]\cong A[t]/(t^{p^m}-\widetilde{(r_0-1)})$.}\bigskip

{\bf Proof.} The element $x$ has a representation 

\begin{equation} x=\mu_0+\mu_1(g-1)+\cdots \mu_{p^m-1}(g-1)^{p^m-1}\,\, (\mu_i\in A; \,\, i=1,2,\cdots,p^m-1) \end{equation}

Assume that $\rho_1(x)=n$. Then $x=x_1+x_2$ where $\rho_1(x_2)>n$, $\rho_1(x_1)=n$ and 

\begin{equation} x_1=\sum_{i=1}^m \mu_i \widetilde{(g-1)
}^{n_i}\end{equation}
 
where all the numbers $n_i$ are distinct,  $0\leq n_i\leq p^m-1$ and $\rho(\mu_i)+kn_i=n$. The homogeneous component of the element $\mu_i (g-1)^{n_i}$ is $\tilde{u}_i \widetilde{(g-1)}^i$. Lemma $2.3.$ now implies that either the homogeneous component of $x_1$ is equal to the right side of $(8.28)$, and in this case it is the homogeneous component of $x$,  or that 

\begin{equation} \sum_{i=1}^m \tilde{\mu}_i \widetilde{(g-1)}^{n_i}=0\end{equation}  

Relation $(8.29)$ is  impossible because the minimal polynomial of $\widetilde{g-1}$ is $t^{p^m}-\widetilde{(r_0-1)}$ whereas all  the coefficients $\tilde{\mu}_i\,\, (i=1,2,\cdots,m)$ in $(8.29)$ belong to $gr(A)$. We obtain from this that the homogeneous component of $x$ is given by $(8.26)$. The uniqueness of this representation follows from Proposition $(8.2)$. 

This completes the proof.\bigskip

We will need one more corollary of Proposition $8.2.$

The ring $R_1$ in Proposition $8.2.$ and Corollary $8.1.$ is isomorphic to a cross product $R\ast G$ where $G$ is the cyclic group of order $p^n$. Assume that conditions of Corollary $8.1.$ hold and let  $\phi$ be the homomorphism $R\longrightarrow R/A$ and $\phi_1$ be the homomorphism 
$R\longrightarrow \bar{R}_1=R/A_1$ where $A_1=AR$, and $\bar{X}$ be the image of a subset $X\in R_1$ under this homomorphism.  Then the restriction of the homomorphism $\phi_1$ on $R$ is $\phi$,  and  $\phi(R)\cong \bar{R}$,    the ring $\bar{R}_1$ is isomorphic to a suitable  cross product $ \bar{R}\ast \bar{G}$ where the group $\bar{G}$ is isomorphic to $G$. This implies  that 
the    ring $\bar{R}_1$ is a free left module with basis $1,\overline{(g-1)},\cdots, \overline {(g-1)}^{p^n-1}$ over the subring $\bar{R}=R/A$ and 
$\overline{(g-1)}^{p^n}=\overline{(r_0-1)}.$ 

We obtain now  from Proposition $2.3.$  that the epimorphism $\phi_1$ together with the  pseudovaluation $\rho_1$ in $R_1$ defines a pseudovaluation $\bar{\rho}_1$ in $\bar{R}_1$ and we have an epimorphism    $\tilde{\phi}_1\colon  gr_{\rho_1}(R_1) \longrightarrow  gr_{\bar{\rho_1}}(\bar{R}_1) $ with kernel $gr_{\rho_1}(A_1)$. Similarly, the  homomorphism $\phi\colon R\longrightarrow  \bar{R}$  and pseudovaluation $\rho$ in $R$ define a pseudovaluation $\bar{\rho}$ in $\bar{R}$ and the  restriction of $\tilde{\phi}_1$ on $gr(R)$ defines an epimorphism 
$\tilde{\phi}\colon gr_{\rho}(R)\longrightarrow gr_{\bar{\rho}} (\bar{R})$. The ring $gr_{\bar{\rho_1}}(\bar{R}_1)$ is generated by the  subring $\tilde{\phi}_1(gr_{\rho}(R)\cong    gr_{\bar{\rho}} (\bar{R})$ and the element $\tilde{\phi}(\theta)$ where $\theta$ is the homogeneous component of the element  $ g-1$ in 
$gr_{\rho_1}(R_1)$. We obtain now from Corollary $8.1.$  that 
$gr_{\bar{\rho_1}}(\bar{R}_1)
$ is isomorphic to the quotient ring $R[\theta]/A[\theta]\cong (R/A)[\theta]\cong \bar{R}[\theta]$ which implies in particular 
that $ \tilde{\phi}(\theta)\not=0$. We obtain from this and from Proposition $2.3.$ that the homomorphism $\tilde{\phi}$   maps the homogeneous component  $\theta=\widetilde {(g-1)}$  in $gr_{\rho_1}(R_1)$ on  the homogeneous component of the element $\overline{g-1}$ in $gr_{\bar{\rho_1}}(\bar{R}_1)$; we denote this homogeneous component by $\bar{\theta}$.

We proved  in these notations the following fact. \bigskip

{\bf Corollary 8.2.} {\it  Let $\theta$ and $\bar{\theta}$ be the homogeneous component of the elements $g-1$ in $gr_{\rho_1}(R_1)$ and in  $gr_{\bar{\rho}_1}(\bar{R}_1)$ respectively. Then }

\begin{equation} \tilde{\phi}(gr_{\rho_1}(R_1))\cong gr_{\bar{\rho_1}}(\bar{R}_1)\cong \bar{R}[\bar{\theta}]\end{equation}

{\it where $\bar{\theta}=\tilde{\phi}(\theta)$ is the homogeneous component of the element $\phi(g-1)$ in $\bar{R}_1$ and the minimal polynomial of $\bar{\theta}$ is $t^{p^m}-\widetilde{(r_0-1)}$.}\bigskip

{\bf 8.2.} Let  $R$ be a ring, $\phi$ be an automorphism of $R$,  $A$ be an $\phi$-invariant ideal in $R$.    Assume that there exists an $\phi$-invariant non-negative pseudovaluation $\rho$ of $R$, and let $\bar{\rho}$ be the pseudovaluation of  the ring 
$\bar{R}=R/A$ obtained from the homomorphism $R\longrightarrow \bar{R}$ (see Proposition $2.3.$)
Let $a_i\,\, (i\in I)$ be a system of elements in $A$ whose homogeneous components $\tilde{a}_i\,\, (i\in I)$ generate the ideal $gr(A)\subseteq gr(R)$. Proposition $2.3.$ implies  that the homomorphism $ R\longrightarrow \bar{R}$ defines in a natural way a homomorphism $gr(R)\longrightarrow gr(\bar{R})$. Let $b_j\,\, (j\in J)$ be a system  of element in $R$ whose homogeneous components $\tilde{b}_j\,\, (j\in J)$ generate $gr(R)$ modulo the ideal $gr(A)$. We can assume that the images of  $\tilde{b}_j\,\, (j\in J)$ in $gr(R)$ are non-zero. 
This implies that for every $j\in J$  the value $\bar{\rho}(\bar{b}_j)$ is equal to the weight $\rho(b_j)$. \bigskip

{\bf Proposition 8.3.} {\it Let $R$ be a ring, $\phi$ be an automorphism of
 $R$,   $\rho$ be a non-negative  $\phi$-invariant  valuation in $R$, $A$ be a $\phi$-invariant ideal of $R$.   Assume that there exists  systems  of  elements $a_i\,\, (i\in )$ whose homogeneous components $\tilde{a}_i\,\, (i\in I)$ generate the ideal $gr(A)$ of $gr(R)$, and a system of elements $b_j\,\, (j\in J)$ whose homogeneous components  $\tilde{b}_j\,\, (j\in J)$ generate the ring $gr(R)$ modulo the ideal $gr(A)$, and  natural number $k$ such that} 

\begin{equation} \rho (\phi(a_i) -a_i) > k+\rho(a_i)\,\, (i\in I)\end{equation}

\begin{equation}  \bar{\rho}(\phi(\bar{b}_j) -\bar{b}_j) > k+
\bar{\rho}(\bar{b}_j)\,\, (j\in J)\end{equation}

{\it Assume also  that for every $i\in I, j\in J$}

\begin{equation}\rho(a_i)> k+\rho(b_j)\end{equation}

{\it Then for every $r\in R$}

\begin{equation} \rho(\phi(r)-r)> k+\rho(r)\end{equation}

We need first the following  lemma.\bigskip

{\bf Lemma 8.3.} {\it Let $R$ be a ring with a non-negative valuation $\rho$, 

$x_i\,\, (i\in I)$ be a system of elements whose homogeneous components $\tilde{x}_i\,\, (i\in I)$ generate the ring $gr(R)$. Let $r\in R$ be an element with $\rho(r)=n$  and $m>n$ be a natural number.

Then 

\begin{equation} r=\sum_{j=1}^s \pi_j+y\end{equation} 

where $\rho(y)\geq m$ and every $\pi_j\,\, (j=1,2,\cdots,s)$ is a monomial with value $n\leq \rho(\pi_j)\leq m-1$ on the set of elements  
$x_i\,\, (i\in I)$.}\bigskip

{\bf Proof.} We find monomials $\pi_j\,\, (j=1,2,\cdots, s_1)$ with values $\rho(\pi_j)=n\,\, (j=1,2,\cdots,s_1)$ such that

 \begin{equation} r=\sum_{j=1}^{s_1}\pi_j+r_1\end{equation}

where $\rho(r_1)>n$. We can obtain a similar  representation for the element 
$r_1$ and representation $(8.35)$ is obtained in a finite  number of steps. 

{\bf Proof of Proposition 8.3}. We see that the  homogeneous components of elements $a_i, b_j\,\, (i\in I,j\in J)$ generate the ring $gr(R)$.   Let $\rho(r)=n$.  We pick the number $m=n+k$ and obtain from Lemma $8.3.$ that  it is enough to prove the assertion for  an arbitrary  monomial  $\pi=\pi_j$ in the representation $(8.35)$ of $r$ where $n\leq \rho(\pi)= n_1\leq m-1$.  We assume that  

\begin{equation}\pi=x_1x_2\cdots x_l\end{equation}

where the elements $x_{\alpha}\,\, (\alpha=1,2,\cdots,l)$ are taken from the set of generators $a_i,b_j\,\, (i\in I;j\in J)$ and $n\leq \rho(\pi)=n_1\leq (m-1)$. The conditions of the assertion yield that it is true when the element $\pi$ is equal to one of the generators, $a_i$ or $b_j$. We can assume that it is true when $\pi$ is a product of $l-1$ generators. We apply now the identity

\begin{equation} \phi(uv)-uv=\phi(u)(\phi(v)-v)+(\phi(u)-u)v\end{equation}

to the elements $u=x_1x_2\cdots x_{l-1}$ and $y=x_l$   and  obtain now for the first summand  in the right side of $(8. 38)$

\begin{eqnarray}\rho(\phi(u)(\phi(v)-v))\geq \rho(\phi(u))+\rho (\phi(v)-v)=\nonumber\\
=\rho(u)+ \rho(\phi(v)-v)>\rho(u)+\rho(v)+k=\nonumber\\
=\rho(uv)+k=\rho(\pi)+k=n_1+k\end{eqnarray}

We obtain in the same way 

\begin{equation}\rho ((\phi(u)-u)v)>n_1+k\end{equation}

and hence 

\begin{equation} \rho(\phi(\pi)-\pi)>n_1+k\end{equation}

Since $\pi$ was an arbitrary summand in the representation $(8.35)$ of $r$ the assertion follows.\bigskip

{\bf 8.3. }    {\bf Theorem XII.} {\it Let $H$ be a torsion free polycyclic group with Hirsch number $r$, $U$ be a normal subgroup with Hirsch number $k$ and torsion free quotient group $\bar{H}=H/U$.   }

{\it Assume that the following $3$ conditions hold.} 

1) {\it There exists in $U$ a $p$-series}
 \begin{equation} U=U_1\supseteq U_2\supseteq \end{equation} 

{ \it with associated restricted graded Lie algebra $L_p(U,U_i)$ free abelian of rank $k$.}

2) {\it  There exists in the group $\bar{H}=H/U$ a $p$-series} 

\begin{equation} \bar{H}=\bar{H}_1\supseteq \bar{H}_2\supseteq\cdots\end{equation} 
{\it with associated restricted  Lie algebra $L_p(\bar{H},\bar{H}_i)$ free abelian of rank $r-k$.} 

3) {\it For for every subgroup $R=gp(h,U)$ generated by $U$ and an element $h\in H$ the quotient group $\bar{R}=H/U^{\prime}U^p$ is a residually \{finite $p$-group\} or equvalently $[ U, h^{p^k}]\subseteq U^{\prime}U^p$. }

{\it Then there exists a $p$-series} 

\begin{equation} H=H_1\supseteq H_2\supseteq\cdots\end{equation}

{ \it  with unit intersection and  associated restricted Lie algebra $L_p(H,H_i)$ free abelian of rank $r$ such that   $\bar{H}_i=(H_iU)/U\,\, (i=1,2,\cdots)$.  }\bigskip

The  last statement of Theorem XII together with Proposition $2.1.$ implies immediately the following   corollary.\bigskip 

{\bf Corollary 8.3.} {\it  The natural homomorphism $\phi\colon H\longrightarrow \bar{H}$ defines a homomorphism  
$\tilde{\phi}\colon L_p(H,H_i)\longrightarrow L_p(\bar{H},\bar{H}_i)$ of graded algebras.}\bigskip

{\bf Proof.} {\it Step 1.} Apply  Theorem $6.1.$ and Corollary $6.2.$    to obtain  in the group $\bar{H}$ a poly\{infinite cyclic\}  normal subgroup $\bar{Q}$ with the the quotient group $G=\bar{H}/\bar{Q}$ a finite $p$-group of order $p^n$,  and a special system of generators 

\begin{equation}\bar{s}_{\alpha},\bar{q}_j\,\, (\alpha=1,2,\cdots,n;j=1,2,\cdots,r)\end{equation}

 where the free abelian  algebra $L_p(\bar{Q},\bar{Q}_i)$ is freely generated by the system of elements  $\tilde{q}_j\,\, (j=1,2,\cdots, r-k)$.     Let $m$ be the maximum of the weights of the elements from this system of generators. We pick an arbitrary natural $M>pm$ and apply Theorem X and Corollary $7.4.$ to get an in $U$  an $H$-invariant $p$-series $U_i\,\, (i=1,2,\cdots,k)$ with unit intersection  such that the weights  of all the elements from $U$ are multiples of   $M$ and the algebra $L_p(U,U_i)$ is centralized by the group $H$.

Let  $Q$ be the inverse image of $\bar{Q}$ in $H$ and 

\begin{equation}s_{\alpha},q_j\,\, (\alpha=1,2,\cdots,n;j=1,2,\cdots,r)\end{equation}

be a system of coset representatives for the elements of system $(8.45)$.
 The group $H$ is obtained from $U$ by a chain of $r-k$ infinite cyclic extensions, and then by a chain of  $n$ cyclic extensions of order $p$; every infinite cyclic extension is generated by some  element $q_j$, every cyclic extension of order $p$ is generated by an element $s_{\alpha} $.  Let $V\supseteq W$ be two subgroups of this series with the quotient group $V/W$ infinite cyclic or cyclic of order $p$,  $\bar{V}=V/U,\bar{W}=W/U$. We  define by $\phi$ the natural homomorphism $\bar{H}=H\longrightarrow H/U$   and  use the same notation  for the restrictions of this homomorphism on the subgroups  $V$ and $W$. 

We consider now the group ring $K\bar{H}$ where $K$ is an arbitrary field of characteristic $p$. Series $(8.43)$ defines a filtration and a valuation $\bar{\rho}$ in the group ring $K\bar{H}$;  we   keep the same notation $\bar{\rho}$ for  the restrictions of this valuation on $K\bar{V}$ and $K\bar{W}$.

 Assume that it has already been proven that the subgroup $W$ contains a $p$-series with unit intersection

\begin{equation} W=W_1\supseteq W_2\supseteq\cdots\end{equation} 

 and with associated algebra $L_p(W,W_i)$ free abelian with rank equal to the Hirsch number of $W$, and that  every  subgroup $ \bar{W}_i=\bar{H}_i\bigcap \bar{W}$  is isomorphic to  $ (W_iU)/U$.  This  assumption   implies that the homomorphism $\phi\colon W\longrightarrow \bar{W}$
 defines the related homomorphism of algebras $  L_p(W,W_i)\longrightarrow L_p(\bar{W},\bar{W}_i)$,   and   a homomorphism $U_p(L_p(W,W_i))\longrightarrow U_p(L_p(\bar{W},\bar{W}_i))$; we will use for the last two   homomorphisms the  notation $\tilde{\phi}$. Further, Proposition  $2.7.$ imply  that 
$gr_{\rho}(KW)\cong U_p(L_p(W,W_i)$, and $gr_{\rho}(K\bar{W})\cong U_p(L_p(\bar{W},\bar{W}_i))$ so  we obtain also a  homomorphism 
$gr_{\rho}(KW)\longrightarrow  gr_{\bar{\rho}}(K\bar{W})$. Once again, we will use for this homomorphism the same notation $\tilde{\phi}$.

We will  prove    that there exists a $p$-series 

\begin{equation} V=V_1\supseteq V_2\supseteq\cdots\end{equation}

such that $V_i\bigcap W=W_i\,\, (i=1,2,\cdots)$, the algebra $L_p(V,V_i)$ is free abelian with rank equal to the Hirsch number of $V$ and 

\begin{equation} (V_iU)/U\cong \bar{V}_i\,\, (i=1,2,\cdots)\end{equation}    

This will imply that after a finite number of steps we will get a required $p$-series $(8.44)$ in $H$. We will give the proof for   for the existence of series $(8.48)$ in the  case when the quotient group $V/W$ is cyclic of order $p$; the case when this quotient group is infinite cyclic is obtained by the same argument with obvious  simplifications.

{\it Step 2.}   We will extend at this step   the valuation $\rho$ in  $KW$ to a valuation $\rho_1$ in $KV$.

Let $s\in V$ be the an element from system $(8.46)$  which generates the quotient group $V/W$; its  image $\bar{s}$ in $H/U$ belongs to $\bar{V}$ and generates the quotient group $\bar{V}/\bar{W}\cong V/W$, so $s^p=a\in W$ and $\bar{s}^p=\bar{a}\in \bar{W}$. Let  $\bar{\rho}(\bar{s}-1)=\delta$.

Let $\bar{\tau}$ be homogeneous component of the element $\overline{s-1}=\bar{s}-1$ in 
$gr_{\bar{\rho}}(K\bar{V})$. Then we obtain from Proposition $8.2.$ 

\begin{equation} gr_{\bar{\rho}}(K\bar{V})\cong gr_{\bar{\rho}}(K\bar{W})[\bar{\tau}]\end{equation} 

where the minimal polynomial of $\bar{\tau}$ is 
$t^p-\widetilde{(\bar{a}-1)}$; here $\widetilde{(\bar{a}-1)}$ denotes the homogeneous component of the element $\bar{a}-1$ in $gr(K\bar{V})$.

Since the algebra $L_p(\bar{H},\bar{H}_i)$ is free abelian with rank equal to the Hirsch number of $\bar{H}$ we obtain from statement i) of Theorem VI that the  $p$-series $\bar{V}_i=\bar{H}_i\bigcap \bar{V}\,\, (i=1,2,\cdots)$ has  unit intersection and the algebra $L_p(\bar{V},\bar{V}_i)$ is free abelian of rank equal to the Hirsch number of $\bar{V}$; we have also  $ \bar{V}_i\bigcap \bar{W}=\bar{H}_i\bigcap \bar{W}=\bar{W}_i\,\, (i=1,2,\cdots)$.

 Let $u$ be an arbitrary element whose homogeneous component belongs to the system of   generators of the algebra $L_p(U,U_i)$.  Since $s$ centralizes the series $U_i\,\, (i=1,2,\cdots)$ we obtain that the weight of the element $s^{-1}us-u$ is greater than the weight of $u$ (see subsection $2.8.$); we obtain from this 

\begin{equation} \rho(s^{-1}(u-1)s-(u-1))> \rho(u-1)\end{equation}  

It is important that the values of the elements of $L_p(U,U_i)$ are multiples of $M$; we recall also that $M>\delta=\rho(\bar{s}-1)$  and obtain from $(8.51)$  

\begin{equation} \rho(s^{-1}us-u)\geq \rho(u-1)+M>\rho(u-1)+\delta\end{equation}

We pick now in $\bar{W}$ an  element $\bar{w}$ whose homogeneous component  belongs to the system of generators of the algebra $L_p(\bar{W},\bar{W}_i)$.  Since the algebra $L_p(\bar{V},\bar{V_i})$ is commutative the homogeneous components of $\bar{s}$ and $\bar{w}$ commute  and we obtain in the same way as in $(8.51)$ and $(8.52)$

\begin{equation} \bar{\rho}(\bar{s}^{-1}\bar{w}\bar{s}-\bar{w})>\delta+\bar{\rho}(\bar{w}-1)\end{equation}

Since the  elements of $U$ have weights greater than $M$ all the elements of $\omega(KU)$  and  all the elements of the ideal $\omega(KU)KW$ have $\rho$-values greater than $M$. Since $\omega(KU)KW$  is the kernel of the homomorphism $KW\longrightarrow K\bar{W}$ we  obtain now from Proposition $8.3$ that for every element $x\in KW$

\begin{equation} \rho(s^{-1}xs-x)>\rho(x-1)+\delta\end{equation}

We recall now that  $s^p=a$,  $\bar{s}^p=\bar{a}$ and the weight of the element $\bar{a}=\bar{s}^p\in \bar{W}$  is equal to $p\delta$ because the algebra $L_p(\bar{V},\bar{V}_i)$ is free abelian. This implies immediately that  the element $a\in W$ has weight less than or equal to $p\delta$.  We will now prove that this weight is equal to $p\delta$. 

In fact,  if the weight of $a$ were  less than $p\delta$, say $\delta_1<p\delta$,    then its image $\bar{a}$ would also have weight  $\delta_1$ because all the elements of the kernel $\omega(KU)KW$ have weights greater than $p\delta>\delta_1$.   This contradiction shows that   the weight of $a$ is $p\delta$; we obtain from this that   $\rho(a-1)=p\delta$.

We can apply now  Proposition $8.2$   to the group $V$ and its normal subgroup $W$. We obtain  that there exists a unique  extension  of the valuation  $\rho$ to a pseudovaluation $\rho_1$  of the algebra $KV$, such that $\rho(s-1)=\delta$, the graded ring $gr(KV)$ is commutative and  we have for the homogeneous components $\widetilde{s-1}$ and $\widetilde{a-1}$ of the elements $s-1$ and $a-1$.

\begin{equation}( \widetilde{s-1})^p=\widetilde{(a-1)}\end{equation}

  We consider now once again the homomorphism $\phi\colon V\longrightarrow \bar{V}=V/U$ and the related homomorphism $KV\longrightarrow K\bar{V}$ of group rings. We obtain  a pseudovaluation $\bar{\rho}_1$ of $K\bar{V}$; if $A_i$ is a filtration defined by $\rho_1$ in $KV$ then $\phi(A_i)=\bar{A}_i\,\, (i=1,2,\cdots)$ is the filtration in $K\bar{V}$ defined by $\bar{\rho}_1$. The restriction of $\bar{\rho}_1$ on $K\bar{W}$ coincides with the valuation $\bar{\rho}$ of $K\bar{W}$.

The $p$-series  $\bar{V}_i^{\prime}\,\, (i=1,2,\cdots)$  in $\bar{V}$ defined by the pseudovaluation $\bar{\rho}_1(\bar{V})$ is in fact obtained as $\bar{V}_i^{\prime}=\phi(V_i)\,\, (i=1,2,\cdots)$; this means  that 

\begin{equation} \bar{V}_i^{\prime}=(V_iU/U)\,\, (i=1,2,\cdots)\end{equation}

Corollary $8.2.$ implies that $gr_{\bar{\rho}_1}(K\bar{V})\cong gr_{\bar{\rho}}(K\bar{W})[\bar{\theta}]$ where $\bar{\theta}$ is the homogeneous component of $\bar{s}-1$ in  $gr_{\bar{\rho_1}}(K\bar{V})$   and it has the same minimal polynomial  $t^p-\widetilde{\phi(a-1)}$ as the element $\bar{\tau}$ in $(8.47)$.  We see that $gr_{\bar{\rho_1}}(K\bar{V})\cong gr_{\bar{\rho}}(K\bar{V})$. 

We obtained two extensions, $\bar{\rho}$ and  $\bar{\rho_1}$,  of the valuation $\bar{\rho}(K\bar{W})$ on $K\bar{V}$ with the same associated graded ring. Proposition $8.2.$ implies that $\bar{\rho_1}=\bar{\rho}$.  We obtain from this that $\bar{V}_i^{\prime}=\bar{V}_i\,\, (i=1,2,\cdots)$ for the $p$-series $\bar{V}_i\,\, (i=1,2,\cdots)$ defined by $\rho$. This together with $(8.56)$ implies that relation $(8. 49)$ holds  and the proof is complete.

\setcounter{section}{9}
\setcounter{equation}{0}

\section*{\center \S 9.  Examples.}

   We construct now two  following examples.

{\it Example 1. } Let $H$ be an infinite cyclic group. We will construct a $p$-series $H_i\,\, (i=1,2,\cdots)$ such that the algebra $L_p(H,H_i)$ is not finitely generated. 
We define a weight function $w(h)$ on $H$ as follows. The elements of $H\backslash H^p$ have weight $1$, the elements of $H^p\backslash H^{p^2}$ have 
weight $p+1$, the elements of $H^{p^2}\backslash H^{p^3}$ have weight $  p(p+1)+1$; if   the weight of the elements from 
$H^{p^n}\backslash  H^{p^{n+1}}$ is $w_n$ then the weight of the elements from   $H^{p^{n+1}}\backslash  H^{p^{n+2}}$ is $pw_n+1$.

This weight function defines a $p$-series $H_i\,\, (i=1,2,\cdots)$ in $H$ and we obtain that if $h\in H$ then $\tilde{h}^p=0$, and $L_p(H,H_i)$ is a restricted   infinite dimensional abelian Lie  algebra of exponent $p$ over $Z_p$.

{\it Example 2. } We will construct now an example of a free abelian group of rank $2$  which contains  a $p$ series $H_i\,\, (i=1,2,\cdots)$ such that the algebra $L_p(H,H_i)$ is free abelian of  rank $1$.

We consider the ring of polynomials  $Z_p[t]$ over the field $Z_p$. Let $R$ be the ring of fractions of $Z_p[t]$ with respect to the complement of the ideal $(t)$. Then $R$ has an ideal $(t)$ with the quotient ring $R/(t)$ isomorphic to $Z_p$, and the powers of this ideal define a $t$-adic valuation $\rho$ in $R$. The graded ring $gr_{\rho}(R)$ is isomorphic to the polynomial ring $Z_p[t]$. We pick the  polynomials $p_1[t]=1+t, p_2[t]=1+t+t^p$. Since every polynomial in $Z_p[t]$ has a unique representation as a product of irreducible polynomials and $p_2[t]$ is not divisible by $1+t$ we obtain that these   polynomials  freely generate a free abelian subsemigroup in the ring $Z_p[t]$. Every element  $p_i[t]\,\, (i=1,2)$ is  invertible in the ring $R$, and we conclude easily,  once again from the uniqueness of factorization in 
$Z_p[t]$, that the elements $u=p_1[t]$ and $v=p_2[t]$  freely generate a free abelian subgroup $H$ of the group of units of $R$. 

  Proposition $2.10.$  implies  that the valuation $\rho$ defines a  $p$-series $H_i\,\, (i=1,2,\cdots)$ in $H$, and the homogeneous components 
$\widetilde{(h-1)} \,\, (h\in H)$ generate in the algebra  $gr_{\rho}(R)\cong Z_p[t]$
a subalgebra  isomorphic to the algebra $L_p(H,H_i)$. This algebra is abelian; we obtain from Corollary $2.6.$ that it   contains no nilpotent elements because the graded ring $gr_{\rho}(R)$ is a domain, so it is free abelian.

We show now that   that the rank of $L_p(H,H_i)$ is $1$.

  Since  the homogeneous components of $p_1[t]-1$  and $ p_2[t]-1$ in 
$gr_{\rho}(R)$ are equal to $t$  we obtain  that $\tilde{u}=\tilde{v}=t$.  
Let $h$ be a non-unit element of $H$, $h=u^n v^m$. The homogeneous components of the elements $u^n-1$ and $v^m-1$ are $nt$ and $mt$ respectively. Since  
$u^nv^m-1=(u^n-1)+(v^m-1)+ (u^n-1)(v^m-1)$ we obtain   that if $n\not=-m$ then the homogeneous component of the element $h-1$ is $(n+m)t$. 

Consider the case when  $n=-m$.    We have to calculate in $R$  the  homogeneous component of the element 

\begin{equation} u^nv^{-n}-1=
(1+t)^n(1+t+t^p)^{-n}-1=((1+t)^n-(1+t+t^p)^n)(1+t+t^p)^{-n}\end{equation}   

The element $(1+t+t^p)$ has $\rho$-value zero  so the homogeneous component of  $(1+t+t^p)^{-n}$ is zero, because of this  the homogeneous component of element $(9.1)$ coincides with the  homogeneous component of $(1+t)^n-(1+t+t^p)^n$.  We consider first the subcase when  $(n,p)=1$ and  replace the term $(1+t)+t^p)^n$ by its binomial expansion and obtain that the homogeneous component of the element $(1+t)^n-(1+t+t^p)^n$ coincides with the homogeneous component of  $-n(1+t)^{n-1}t^p$  which is is equal to $-n t^p$.

In the general case, we have $n=p^k n_1$ where $(n_1,p)=1$. Since the    homogeneous component of the element $u^{n_1}v^{n_1}-1$ is $-n_1t^p$  the homogeneous component of the element $u^nv^{-n}-1$ is $(-n_1t^p)^{p^k}=-n_1t^{p^{k+1}}$.  

We obtained  that the homogeneous components $\tilde{h}\,\, (h\in H)$ generate in $gr(R)$ the free abelian   restricted  Lie subalgebra  with the  generator $t$. This  proves that the algebra $L_p(H,H_i)$ is free abelian of rank $1$.

\newpage

\begin{center} {\bf References} \\ \end{center}
 \noindent

1. A. Yu. Bachturin, Identical Relations in Lie Algebras, VNU Science Press, Utrecht, 1987.\\

\noindent

2. N. Bourbaki,  Algebra,  II, Chap.4-7, Springer, Heidelberg,1990.\\

\noindent

3. P. M. Cohn, Skew Fields, ( Theory of General Division Rings), Cambridge University Press, Cambridge, 1977.\\

\noindent
 
4. B. H. Hartley, Topics in nilpotent groups, (in Group Theory, Essays for Philip Hall), Academic Press, New York, 1984, pp. 61-120.\\

\noindent

5. N. Jacobson, Lie Algebras, Wiley, New York/London, 1989. \\

\noindent

6. M. I. Kargapolov and Yu.I. Merzlyakov, Fundamentals of the Theory of  Groups, Nauka, Moscow, 1972 (in Russian) (English translation: Springer-Verlag, 1979, NY). \\

\noindent

7. M. Lazard
, Sur les groupes nilpotents et les anneaux de Lie, Ann. Sci. Ecole Norm. Sup. (3) (1954), 101-190.\\

\noindent

8. A. I. Lichtman, Valuation methods in division rings, Journal of Algebra, 177 (1993), 870-898.\\

\noindent

9. A. I. Lichtman, Valuation methods in group rings and in skew fields, I, J. Algebra 257 (2002), 106-167.\\

\noindent

10. A. I. Lichtman, Restricted Lie algebras of polycyclic groups, Journal of Algebra and Its Applications, 5 (2006), 571-627.\\

\noindent

11. A. I. Lichtman, Restricted Lie algebras of polycyclic groups, III, (in preparetion)\\

\noindent

12. A. I. Lichtman, Matrix rings and linear groups over a field of fractions of enveloping algebras and group rings, I, J. Algebra 88 (1984), 1-37.\\

\noindent

13. D. S. Passman, The Algebraic Sructure of Group Rings, Wiley, New York, 1971.\\

\noindent

14. D. S. Passman, Universal fields of fractions for polycyclic group algebras, Glasgow Math. J. 23 (1982), 103-113. \\

\noindent

15.  D. G. Quillen, On the associated graded ring of a group ring, J. Algebra 10, (1968), 711-718.\\

\noindent

16. D. Segal, Polycyclic Groups, Cambridge Univ. Press, 1983.\\

\noindent

\end{document}